 \documentclass[10pt,reqno]{amsart}
 \usepackage{amsmath}
 \usepackage{amssymb}
 \usepackage{multicol}
 \usepackage{amsfonts}
\newcommand{\R}{{\mathbb R}}
 \newcommand{\N}{{\mathbb N}}
 \newcommand{\Z}{{\mathbb  Z}}
 
 \newcommand{\C}{{\mathbb  C}}

 \newcommand{\hf}{\frac{1}{2}}

 \newcommand{\intR}{\int_{-\infty}^{\infty}}

 \renewcommand{\Re}{{\mathrm {Re}}\ }
 \renewcommand{\Im}{{\mathrm {Im}}\ }
 \newcommand{\generic}{\left( \begin{smallmatrix} a & b \\ c & d
 \end{smallmatrix} \right)}

 \renewcommand{\vec}{\boldsymbol}

\renewcommand{\sc}[1]{\langle#1\rangle}

 \renewcommand{\l}{\left}
 \renewcommand{\r}{\right}

 \renewcommand{\include}{\input}

\begin{document}
 \theoremstyle{plain}
 \newtheorem{theorem}{Theorem}[section]
 \newtheorem{proposition}[theorem]{Proposition}
 \newtheorem{lemma}[theorem]{Lemma}
 \newtheorem{corollary}[theorem]{Corollary}

 \theoremstyle{definition}
 \newtheorem{definition}[theorem]{Definition}

\def\sgn{\mathop{\rm{sgn}}\nolimits}
\def\Ind{\mathop{\rm{Ind}}\nolimits}
\def\g{\mathfrak{g}}
\textwidth 14cm \textheight 23cm

 \numberwithin{equation}{section}
\title[PSEUDODIFFERENTIAL ANALYSIS  ON PROJECTIVE SPACE] {PROJECTIVE PSEUDODIFFERENTIAL ANALYSIS
AND  HARMONIC ANALYSIS}
\author{Michael Pevzner,\ Andr\'e Unterberger}
\address{Math\'ematiques (UMR 6056)\\
Universit\'e de Reims\\ Moulin de la Housse, B.P.1039\\ F-51687
REIMS Cedex 2}
 \email{pevzner@univ-reims.fr,andre.unterberger@univ-reims.fr}
\begin{abstract}
We consider pseudodifferential operators on functions on
\,$\R^{n+1}$\, which commute with the Euler operator, and can thus
be restricted to spaces of functions homogeneous of some given
degree. Their symbols can be regarded as functions on a reduced
phase space, isomorphic to the homogeneous space
\,$G_n/H_n=SL(n+1,\R)/GL(n,\R)$\,, and the resulting calculus is a
pseudodifferential analysis of operators acting on spaces of
appropriate sections of line bundles over the projective space
\,$P_n(\R)$\,: these spaces are the representation spaces of the
maximal degenerate series \,$(\pi_{i\lambda,\varepsilon})$\, of
\,$G_n$\,. This new approach to the quantization of \,$G_n/H_n$\,,
already considered by other authors, has several advantages: as an
example, it makes it possible to give a very explicit version of the
continuous part from the decomposition of \,$L^2(G_n/H_n)$\, under
the quasiregular action of \,$G_n$\,. We also consider interesting
special symbols, which arise from the consideration of the
resolvents of certain infinitesimal operators of the representation
\,$\pi_{i\lambda,\varepsilon}$\,.
\end{abstract}
 \date{\today}

\maketitle
\tableofcontents

\setcounter{section}{0}
\section*{Introduction}

This paper is devoted to harmonic analysis on the homogeneous space
\,$G_n/H_n\\=SL(n+1,\,\R)/GL(n,\R)$\, and more precisely to a study,
on this example, of the interaction between harmonic analysis and
pseudodifferential analysis. We here combine two ideas, both of
which stem from a long--standing tradition in mathematics or
physics: the  {\em superselection method\/}  in pseudodifferential
analysis, and the
{\em square--root method\/} in the analysis of certain second--order differential operators.\\

Recalling here the definition of a not over--specialized version of
pseudodifferential analysis on \,$\R^{n+1}$\, will save space later
(the dimension \,$n+1$\, is of course meant for coherence with the
sequel). A {\em symbolic calculus\/} of operators on \,$\R^{n+1}$\,
is a linear way to associate linear operators \,${\mathrm{Op}}(H)$\,
on functions of $n+1$ variables to functions \,$H$\, of \,$2n+2$
variables: we are not interested in  an axiomatization of the
concept, but in the following parameter--dependent special case.
Given \,$\kappa\in \R$\,, consider the defining equation, in which
the integration with respect to  \,$dy$\, is to  be carried first:
\begin{equation}\label{1Op}
({\mathrm{Op}}_{\kappa}(H)\,u)(x)=\int_{\R^{n+1}}\int_{\R^{n+1}}
H((1-\kappa)\,x+\kappa\, y,\,\eta)\, e^{2i\pi\langle
x-y,\,\eta\rangle}\ u(y)\ d\eta\,dy\,.
\end{equation}
In the case when \,$\kappa=\hf$\,, this is the Weyl calculus of
operators, or Weyl pseudodifferential analysis. In the case when
\,$\kappa=0$\,, this is the {\em standard pseudodifferential
calculus\/}, or {\em convolution--first\/} calculus, a terminology
which the reader will feel justified after he has examined the case
when the symbol \,$H$\, decomposes as
\,$H(x,\eta)=h_1(x)\,h_2(\eta)$\,: whether this is the case or not,
one can reduce the double sign of integration to a simple one with
the help of the Fourier transform of \,$u$\,. One cannot introduce
the standard calculus without considering, at the same time, the
{\em antistandard calculus\/}, which is the case \,$\kappa=1$\, of
the formula above, would it be only for the fact that the adjoint of
the operator \,${\mathrm{Op}}_0(H)$\, is the operator
\,${\mathrm{Op}}_1(\bar{H})$\,.\\

What we call the superselection method originates from the
physicists'  superselection rule: we want to devote our interest to
a special class of operators, to wit those which commute with some
fixed differential operator \,$M$\, with symbol \,$m$\,, with a
given self--adjoint realization. It {\em may\/} happen that the
symbols \,$H$\, of such operators are exactly the ones satisfying
the Poisson bracket equation \ $\{m\,,\,H\}=0$\,: this will be the
case when \,$M$\, is the infinitesimal operator of a one--parameter
unitary group lying in the covariance group of the symbolic
calculus, a group which contains the metaplectic group when
\,$\kappa=\hf$\,, and the group \,$SL(n+1,\,\R)$\, of
transformations of \,$L^2(\R^{n+1})$\, in all other cases. In
significant instances, reducing the phase space \,$\R^{n+1}\times
\R^{n+1}$\, so as to account for this equation cannot fail to lead
to a geometrically interesting structure: it does, indeed
\cite{klg}, in the case when \,$n=3$\, and \,$M$\, is the d'Alembert
operator \ $\frac{\partial^2}{\partial x_0^2}-\sum_{j=1}^3
\frac{\partial^2}{\partial x_j^2}$\,, since it leads in a canonical
way to the basic geometric concepts which occur in the development
of special relativity. Moreover, given any positive number
\,$\mu$\,, the reduced phase space is then the natural choice for a
pseudodifferential analysis of operators acting on the space of
solutions of the equation \ $Mu+\mu\,u=0$\,. In the case just
alluded to, this spectral equation is the Klein--Gordon equation,
and the resulting Klein--Gordon pseudodifferential analysis --- a
reduction of the Weyl calculus --- was developed in ({\em loc.cit.\/}). \\

We here consider a reduction of the standard--antistandard
pseudodifferential analysis. It goes along the lines of the above
general scheme, if one chooses for \,$M$\, the Euler operator \
$(2i\pi)^{-1}\, (\sum_{k=1}^{n+1} x_k\,\frac{\partial}{\partial
x_k}+\frac{n+1}{2})$\,. In other words, we consider operators  on
functions of \,$n+1$\, variables which preserve the homogeneity ---
and the parity as well --- of functions: these operators then give
rise to operators acting on functions defined on the projective
space \,$P_n(\R)$\,. The reduced phase space turns out to be the
homogeneous space \,${\mathcal X}_n^{\bullet}=G_n/H_n$\,, and the
corresponding pseudodifferential
analysis may also be referred to as a quantization of \,$G_n/H_n$\,.\\

Such a quantization was studied --- in a way independent from the
superselection method --- in \cite{U3} in the case when \,$n=1$\,
and, in the general case, in a series of papers \cite{van2, vdm}. In
all cases, the development has to include a description of the
decomposition of \,$L^2(G_n/H_n)$\, under the quasiregular
representation of \,$G_n$\, in this space, which amounts to a
description of the spectral  decomposition of the basic
\,$G_n$--invariant differential operator \,$\Delta_n$\, on
\,${\mathcal X}_n^{\bullet}=G_n/H_n$\,. In the one--dimensional
case, one could dispense with this task, relying instead on results
relative to the harmonic analysis of hyperboloids \cite{stri}. In
general, in the references just given, the problem was dealt with by
the use of the \,$H_n$--spherical distribution method. We wish to
give some
idea, in the present introduction, of the {\em square--root\/} method adopted here.\\

It consists in replacing a second--order operator, here
\,$\Delta_n$\,, by an equivalent, or almost equivalent --- it
usually provides more information --- first-order operator or system
of operators. This idea resurfaces in a variety of domains and
disguises. Without any attempt at completeness, let us recall the
following well--known, or not so well--known, instances. The first
one that springs to mind is Dirac's replacement of the second--order
Klein--Gordon equation by his system of four first--order equations;
by the way, this can be followed up, again, in the domain of
pseudodifferential analysis, leading to the construction of the {\em
Dirac symbolic calculus\/} of operators \cite{dir}. Another circle
of ideas, quite close to the one which we will adhere to in the
present work, is familiar to harmonic analysts and deals with such
objects as the Weyl group and Harish--Chandra's isomorphism. We
prefer to come to it in terms most readers will probably not be
quite as familiar with, starting from the Lax--Phillips scattering
theory for the automorphic wave equation \cite{lap}. Automorphic
functions are functions in the upper half--plane invariant under the
action, by fractional--linear transformations, of some arithmetic
group: when at the same time generalized eigenfunctions of the
non-Euclidean Laplacian \,$\Delta$\,, they are called
non--holomorphic modular forms. In the Lax--Phillips scattering
theory, pairs of automorphic functions are made to appear as the set
of Cauchy data on some hyperboloid for the d'Alembert equation in
the three--dimensional forward light--cone. In \cite{aumod}, it was
shown that the space of such pairs can be identified with functions
on \,$\R^2$\,, in such a way that, under the transfer, the operator
\,$\Delta-\frac{1}{4}$\, becomes the square of the first--order
operator \ $(2i)^{-1}\,(x_1\,\frac{\partial}{\partial x_1}
+x_2\,\frac{\partial}{\partial x_2}+1)$\,. A fully similar idea,
replacing \,$\Delta$\, by \,$\Delta_n$\,, will work here. Let us
just mention {\em en passant\/} that the concept of automorphic
distribution that arose from this transfer made automorphic
pseudodifferential
analysis \cite{birk} possible.\\

More details follow: the superselection rule present here calls for
the consideration of symbols \,$H=H(x,\,\xi)$\, invariant under the
ever--present action \ $t.\,(x,\,\xi) =(tx,\,t^{-1}\xi)$\, of the
group \,$\R^{\times}$\,; then, \,${\mathcal X}_n^{\bullet}$\, may be
realized as the hypersurface, in the corresponding quotient, of
equation \ $\sc{x,\,\xi}=1$\,. Next, we consider in
\,$\R^{n+1}\times \R^{n+1}$\, the operator \ $\square= \sum
_{k=1}^{n+1} \frac{\partial^2}{\partial x_j\,\partial{\xi}_j}$\,.
This operator will turn out, in a moment, to be a fundamental one in
connection with pseudodifferential analysis. For the time being,
its importance stems from what follows: within the domain
\,$\Omega^+$\, defined by \ $\sc{x,\,\xi}>0$\,, set \
$\tau=\log\,\sc{x,\,\xi}$\,, which provides an identification of the
quotient of \,$\Omega^+$\, by the group \,$\R^{\times}$\, with the
product \,${\mathcal X}_n^{\bullet}\,\times \R$\,. Then, under the
transformation \ $H \mapsto H_1=e^{\frac{n\tau}{2}}\,H$\,, the
equation \ $\square \, H=0$\, is equivalent to the wave equation \
$\frac{\partial^2 H_1}{\partial
\tau^2}+(\Delta_n-\frac{n^2}{4})\,H_1=0$\,. This explains a
fundamental property of one part at least (the continuous one) of
the decomposition of \,$L^2({\mathcal
X}_n^{\bullet})=L^2(G_n/H_n)$\,, to wit the fact that the
generalized eigenvalues always come by pairs \,$(\rho,\,-n-\rho)$\,.
Though one can trace this to several possible sources, the
following, a continuation of the Lax--Phillips point of view, seems
to us especially striking: solutions of the wave equation above can
be characterized by their first two traces on \,${\mathcal
X}_n^{\bullet}$\,, not just one. On the other hand, constructing
\,$\R^{\times}$--invariant solutions, in \,$\R^{n+1}\times
\R^{n+1}$\,, of the equation \ $\square \, H=0$\, can be achieved by
means of a Fourier transformation, starting from functions
\,$\Phi$\, on \,$\Sigma_n^{\bullet}$\,, the quotient by the now
familiar action of \,$\R^{\times}$\,  of the cone \,$\Sigma_n$\, of
equation \ $\sc{y,\,\eta}=0$\,. It is then not surprising that there
exists an involution \,${\mathcal K}$\, of this latter space of
functions such that two \,${\mathcal K}$--related functions
\,$\Phi$\, and \,${\mathcal K}\,\Phi$\, should always lead to
solutions of the above wave equation with the same first trace on
\,${\mathcal X}_n^{\bullet}$\,. Analyzing the involution
\,${\mathcal K}$\, leads without too much difficulty to a full
understanding of the continuous part of
\,$L^2(G_n/H_n)$\, (section 3).\\

Next (section 2), we decompose functions on \,$\R^{n+1}$\, into
their homogeneous components: since the Euler operator commutes with
the linear action of \,$G_n$\,, this action decomposes as a
continuous sum \ $(\pi_{i\lambda,\varepsilon})$\, with \,$\lambda\in
\R$\, and \,$\varepsilon=0$\, or $1$\,, a ``series'' of irreducible
unitary representations in \,$L^2(\R^n)$\, also arising from the
general theory \cite{kna} as a maximal degenerate series of
representations of the group \,$G_n$\,; then, the  Fourier
transformation decomposes as the family of intertwining operators
relative to this series. The decomposition, along the general lines,
of the standard--antistandard pseudodifferential analysis, leads for
every pair \,$(i\lambda,\varepsilon)$\, to the definition of two
linear maps \ ${\mathrm{Op}}_{i\lambda,\varepsilon}$\, and \
${\mathrm{Op}}^{\vee}_{i\lambda,\varepsilon}$\, from functions on
\,${\mathcal X}_n^{\bullet}$\, to linear operators in the space of
the representation \ $\pi_{i\lambda,\varepsilon}$\,: these two
symbolic calculi are, of course, exactly the ones
used in the above--given references concerning the quantization of the space \,$G_n/H_n$\,.\\

The occurrence of the operator \,$\square$\, above is remarkable
since the equation \ $\square \,H=0$\, exactly means that the
operator \ ${\mathrm{Op}}_{\kappa}(H)$\, does not depend on
\,$\kappa$\,. In particular, the operator, on \,$\R^{n+1}$\,, with
symbol \,$H$\,, is the same, whether one considers \,$H$\, as a
standard or antistandard symbol. Now, it is easy to connect the
symbol, in the \ ${\mathrm{Op}}_{i\lambda,\varepsilon}$-- calculus,
of the associated operator, to \,$H$\, viewed as a standard symbol,
while the symbol in the \
${\mathrm{Op}}^{\vee}_{i\lambda,\varepsilon}$--calculus of an
operator is easily connected to the antistandard symbol of the
operator on \,$\R^{n+1}$\, it comes from. In this way, one finds
(section 2) a simple proof of a formula, first given in \cite{vdm}
when \,$n>1$\,, connecting the \
${\mathrm{Op}}_{i\lambda,\varepsilon}$\, and \
${\mathrm{Op}}^{\vee}_{i\lambda,\varepsilon}$\,-- symbols of the
same operator. However, the method only works for symbols lying in
the continuous part of the decomposition of
\,$L^2(G_n/H_n)$\,.\\

To prevent some possible misunderstanding, let us emphasize two
points, both related to the fact that the covariance group of the
pseudodifferential analysis under discussion is
\,$G_n=SL(n+1,\R)$\,, not \,$O(n)$\,. There does not exist on the
projective space \,$P_n(\R)$\, any measure invariant under
\,$G_n$\,: still, a representation such as
\,$\pi_{i\lambda,\varepsilon}$\, is unitary because it really acts,
rather than on functions on \,$P_n(\R)$\,, on sections of some line
bundle; only, this fact is sometimes blurred by the use on
\,$\R^n\subset P_n(\R)$\, of affine coordinates. Next, there exists
on \,$P_n(\R)$\,, viewed if so wished as the usual quotient of the
sphere \,$S^n$\,, a vast family of pseudodifferential analyses
covariant under the action of the orthogonal group: the \
${\mathrm{Op}}_{i\lambda,\varepsilon}$\,-- calculus is covariant
under some specific action of \,$G_n$\,, which makes it almost
unique, since any two such calculi have to be related under a
transformation of functions on \,$G_n/H_n$\, expressing itself, in
spectral--theoretic terms, as
a function of the operator \,$\Delta_n$\,.\\

The problem of analyzing the {\em sharp composition\/} --- the
terminology is the one in use in pseudodifferential analysis --- of
symbols, by which is meant the bilinear operation that corresponds
to the composition of the associated operators, is a difficult one
which, in the case when \,$n=1$\,, was partly solved (for symbols
lying in the discrete part of the decomposition of \,$L^2(G_1/H_1)$)
in \cite{U3}, where the Rankin--Cohen brackets were shown for the
first time to have a significant role in pseudodifferential
analysis. We do not solve it here for general \,$n$\,, but we seize
this opportunity to show that the integral formula for the sharp
composition of symbols --- which is trivial to obtain --- is very
far from revealing the more interesting aspects of the operation
under consideration. Our main point, in section 4, will be to do
away, on the basis of it, with two related popular misconceptions:
one of them consists in pushing too far the concept that the inverse
of the parameter \,$\lambda$\, that specifies an irreducible
representation of \,$G_n$\, within its series might be interpreted
as a ``Planck's constant''; the other one consists in believing that
the composition of symbols can be, in some reasonable sense,
approximated by a
series of bidifferential operators.\\

Some functions on \,$G_n/H_n$\,, while not in \,$L^2(G_n/H_n)$\,
when \,$n\geq 2$\, (in the one--dimensional case, these functions
lie in the discrete part of the decomposition of \,$L^2(G_1/H_1)$)
are very interesting  to consider in view of the role they play in
the symbolic calculus: for they provide the symbols of certain
operators in the algebra generated by resolvents of elements of the
(complexified) space of infinitesimal operators of the
representation \,$\pi_{i\lambda,\varepsilon}$\,. These symbols are
introduced in Section 5, where it is also shown that the
above--mentioned formula, linking the two species of symbols of the
same operator, continues to hold in this new context. The analysis
of individual operators obtained in this way --- which played an
essential role, when \,$n=1$\,, in \cite{U3} ---
can, up to some point, be reduced to the one--dimensional case.\\

To conclude, let us make it clear that, though the present paper
certainly provides more familiarity with the \
${\mathrm{Op}}_{i\lambda,\varepsilon}$-- calculus, we are still far
from having reached a point where this could be considered as a
genuine pseudodifferential analysis in the sense demanded, say, by
possible applications to partial differential equations:
developments in this direction may prove surprisingly new, in
particular in view of the fact that the representations
\,$\pi_{i\lambda,\varepsilon}$\, are not
square--integrable.\vskip15pt


\section{Pseudodifferential analysis, from \,$\R^{n+1}$\, to \,$P_n(\R)$}

The projective space \,$P_n(\R)$\, is the quotient of
\,$\R^{n+1}\backslash\{0\}$\, by the equivalence that identifies two
vectors when proportional: we denote as \,$x\mapsto x^{\bullet}$\,
the projection map. The vector \,$x=(x_1,\dots,x_{n+1})$\, is called
a set of homogeneous coordinates of \,$x^{\bullet}$\,: we shall also
represent \,$x^{\bullet}$\, by the vector
\,$s=x_{n+1}^{-1}\,(x_1,\dots,x_n)$\, in the case when
\,$x_{n+1}\neq 0$\,. The group \,$G_n=SL(n+1,\R)$\, acts on
\,$\R^{n+1}$\, in the linear way, which defines an action on
\,$P_n(\R)$\, too, denoted as \ $(g,\,s)\mapsto [g]\,s$\, in
inhomogeneous coordinates (note that the use of inhomogeneous
coordinates makes this action look like a singular one, which it is not). \\

We first decompose the Hilbert space \,$L^2(\R^{n+1})$\, under the
action \,$(g,\,v)\mapsto v\,\circ\,g^{-1}$\, of \,$G_n$\,. This
action preserves the parity of functions, and we denote as
\,$L^2_{\varepsilon}(\R^{n+1})$\,, with \,$\varepsilon=0$\, ({\em
resp.\/} $1$) the subspace of \,$L^2(\R^{n+1})$\, consisting of even
({\em resp.\/} odd) functions. Given \,$v=v_0+v_1\in
L^2_0(\R^{n+1})\oplus L^2_1(\R^{n+1})$\,, decompose it as
\begin{equation}\label{2declamb}
v=\sum_{\varepsilon=0,\,1} \intR
v_{i\lambda,\varepsilon}\,d\lambda\,,
\end{equation}
where the function
\begin{equation}\label{2lambdapart}
v_{i\lambda,\varepsilon}(x)=\frac{1}{2\pi}\,\int_0^{\infty}
t^{\frac{n-1}{2}+i\lambda}
\,v_{\varepsilon}(tx)\,dt=\frac{1}{4\pi}\,\intR
|t|_{\varepsilon}^{\frac{n-1}{2}+i\lambda} \,v(tx)\,dt
\end{equation}
is homogeneous of degree and parity
\,$(-\frac{n+1}{2}-i\lambda,\,\varepsilon)$\,, a phrasing that we
shall adopt for brevity: we here set
\begin{equation}
|t|_{\varepsilon}^{\alpha}=|t|^{\alpha}\,({\mathrm{sign}}\,t)^{\varepsilon}\qquad
{\mathrm{for}} \ t\in \R\backslash\{0\}\,,\ \alpha\in\C\,.
\end{equation}
The function \,$v_{i\lambda,\varepsilon}$\, is of course
characterized by the function \,$v_{i\lambda,\varepsilon}^{\flat}$\,
on \,$\R^n$\, such that
\begin{equation}\label{2defflat}
v_{i\lambda,\varepsilon}^{\flat}(s)=v_{i\lambda,\varepsilon}(s,1)
\end{equation}
since, with \,$x=(x_*,\,x_{n+1})$\,, one has
\begin{equation}\label{2charac}
v_{i\lambda,\varepsilon}(x)=|x_{n+1}|_{\varepsilon}^{-\frac{n+1}{2}-i\lambda}\,
v_{i\lambda,\varepsilon}^{\flat}(\frac{x_*}{x_{n+1}})\,.
\end{equation}
Applying the equation
\begin{equation}
\intR
|(v_{i\lambda,\varepsilon}(x)|^2\,d\lambda=\frac{1}{2\pi}\,\int_0^{\infty}
t^n\,|v_{\varepsilon}(tx)|^2\,dt\,,
\end{equation}
valid for almost every \,$x\in \R^{n+1}\backslash \{0\}$\,, with
\,$x=(s,1)$\,, and integrating the result with respect to \,$ds$\,,
we obtain
\begin{equation}
\Vert\,v\,\Vert^2_{L^2(\R^{n+1})}=4\pi\,\sum_{\varepsilon=0,\,1}
\intR \Vert\,v_{i\lambda,\varepsilon}^{\flat}\,\Vert^2_{L^2(\R^n)} \
d\lambda\,.
\end{equation}\\

Next, given \,$g\in G_n$\, of the form \ $g=\l(\begin{smallmatrix}  M  & \vec{p} \\
\vec{q}^{\dag} &m \end{smallmatrix}\r)$\,, where \,$\vec{p}\in
\R^n$\, is a column vector and \,$\vec{q}^{\dag}$\, is the transpose
of the column vector \,$\vec{q}\in \R^n$\, , one has for every
\,$s\in \R^n$\,, as a consequence of (\ref{2declamb}), the equation
\begin{equation}
(v\circ g)\,(s,1)=\sum_{\varepsilon=0,\,1} \intR
v_{i\lambda,\varepsilon}(Ms+\vec{p},\,\langle
\vec{q},\,s\rangle+m)\,d\lambda
\end{equation}
or, using the homogeneity,
\begin{equation}
(v\circ g)\,(s,1)=\sum_{\varepsilon=0,\,1} \intR |\langle
\vec{q},\,s\rangle+m|_{\varepsilon}^{-\frac{n+1}{2}-i\lambda} \
v_{i\lambda,\varepsilon}(\frac{Ms+\vec{p}}{\langle
\vec{q},\,s\rangle+m}\,,1)\,d\lambda\,.
\end{equation}
It follows that
\begin{equation}\label{2transform}
(v\circ g)_{i\lambda,\varepsilon}^{\flat}\,(s)=|\langle
\vec{q},\,s\rangle+m|_{\varepsilon} ^{-\frac{n+1}{2}-i\lambda}\
v_{i\lambda,\varepsilon}^{\flat}\,(\frac{Ms+\vec{p}} {\langle
\vec{q},\,s\rangle+m})\,.
\end{equation}
Set
\begin{equation}\label{2defpi}
\pi_{i\lambda,\varepsilon}(g^{-1})\,v_{i\lambda,\varepsilon}^{\flat}=
(v\circ g)_{i\lambda,\varepsilon}^{\flat}\,:
\end{equation}
since
\begin{equation}
\biggr|\,\frac{D\,([g]s)}{Ds}\,\biggr|=|\langle
\vec{q},\,s\rangle+m|^{-n-1}\,,
\end{equation}
the representation \,$\pi_{i\lambda,\varepsilon}$\, of \,$G_n$\, in
\,$L^2(\R^n)$\, so introduced
is unitary.\\

Together with \,$\pi_{i\lambda,\varepsilon}$\,, we consider the
contragredient representation
\,$\pi_{i\lambda,\varepsilon}^{\sharp}$\, defined by the equation
(in which \,$g\mapsto g'$\, denotes the matrix transposition)
\begin{equation}
\pi_{i\lambda,\varepsilon}^{\sharp}(g')=\pi_{i\lambda,\varepsilon}(g^{-1})\,.
\end{equation}
Although the formal definition of the intertwining operator
\,$\theta_{i\lambda,\varepsilon}$\, from the representation
\,$\pi_{i\lambda,\varepsilon}$\, to the representation
\,$\pi_{-i\lambda,\varepsilon}^{\sharp}$\, is a consequence of the
general theory, a better understanding of its properties can be
obtained from its definition in terms of the usual Fourier
transformation
on \,$L^2(\R^{n+1})$\, ({\em cf.\/} \cite[p.\,28]{aumod} for the case when \,$n=1$).\\

Applying the Fourier transformation, normalized as
\begin{equation}
({\mathcal F}\,v)(x)=\int_{\R^{n+1}} v(y)\,e^{-2i\pi\langle
x,\,y\rangle}\,dy\,,
\end{equation}
to both sides of (\ref{2declamb}), and noting that the Fourier
transformation sends functions homogeneous of degree
\,$-\frac{n+1}{2}-i\lambda$\, to functions homogeneous of degree
\,$-\frac{n+1}{2}+i\lambda$\, with the same parity, we obtain
\begin{equation}
{\mathcal F}\,v=\sum_{\varepsilon=0,\,1} \intR {\mathcal
F}\,v_{i\lambda,\varepsilon}\,d\lambda= \sum_{\varepsilon=0,\,1}
\intR ({\mathcal F}\,v)_{-i\lambda,\varepsilon}\,d\lambda\,.
\end{equation}
We may thus define
\begin{equation}\label{2deftheta}
\theta_{i\lambda,\varepsilon}\,v^{\flat}_{i\lambda,\varepsilon}=({\mathcal
F}\,v) ^{\flat}_{-i\lambda,\varepsilon}\,.
\end{equation}
Checking that the operator \,$\theta_{i\lambda,\varepsilon}$\, has
the required intertwining property is easy: indeed, given \,$g\in
G_n$\,, one has on one hand, applying (\ref{2defpi}) and the
definition just given,
\begin{equation}
\theta_{i\lambda,\varepsilon}\,\pi_{i\lambda,\varepsilon}(g^{-1})\,v^{\flat}_{i\lambda,\varepsilon}
=\theta_{i\lambda,\varepsilon}\,(v\circ
g)^{\flat}_{i\lambda,\varepsilon}= ({\mathcal F}\,(v\circ
g))^{\flat}_{-i\lambda,\varepsilon}= \l[\,({\mathcal
F}\,v)\,\circ\,{g'}^{-1}\,\r]^{\flat}_{-i\lambda,\varepsilon}\,,
\end{equation}
on the other hand
\begin{equation}
\pi_{-i\lambda,\varepsilon}(g')\,\theta_{i\lambda,\varepsilon}\,v^{\flat}_{i\lambda,\varepsilon}=
\pi_{-i\lambda,\varepsilon}(g')\,({\mathcal
F}\,v)^{\flat}_{-i\lambda,\varepsilon}= \l[\,({\mathcal
F}\,v)\,\circ\,{g'}^{-1}\,\r]^{\flat}_{-i\lambda,\varepsilon}\,.
\end{equation}\\

Note that, even though (\ref{2deftheta})  defines the function
\,$\theta_{i\lambda,\varepsilon}\,v^{\flat}_{i\lambda,\varepsilon}$\,
for almost all \,$\lambda$\, only, it is easy to introduce classes
of functions \,$v$\, such that, for every \,$s\in \R^n$\,, \
$v^{\flat}_{i\lambda,\varepsilon}(s)$\, has a well--defined meaning
for every \,$\lambda$\,. A simple, useful example is provided by the
space \,${\mathcal S}_{\mathrm{flat}}(\R^{n+1})$\, consisting of all
{\em flat\/} functions in \,${\mathcal S}(\R^{n+1})$\,, {\em
i.e.,\/} functions in this latter space every derivative of which is
bounded, near \,$0$\,, by a constant times
such power of \,$|x|$\, as one may wish.\\

We must now compute the operator \,$\theta_{i\lambda,\varepsilon}$\,
as an integral operator in terms of the inhomogeneous coordinates on
\,$P_n(\R)$\,. The computations that follow have only formal value,
but we plead not guilty on this account: for we only need to compare
the perfectly valid definition (\ref{2deftheta}) of the intertwining
operator to the formal one taken from the general theory. The part
homogeneous of degree \,$-\frac{n+1}{2}+i\lambda$\, and of parity
\,$\varepsilon$\, of the function \,$x\mapsto e^{-2i\pi\langle
x,\,y\rangle}$\, is given by the (divergent) integral, taken from
(\ref{2lambdapart}), to be interpreted as defining a Fourier
transform,
\begin{multline}
\frac{1}{4\pi}\,\intR |t|_{\varepsilon}^{\frac{n-1}{2}-i\lambda}\,
e^{-2i\pi\,t\langle x,\,y\rangle}\ dt=\\
\frac{1}{4\pi}\,(-i)^{\varepsilon}\,\pi^{-\frac{n}{2}+i\lambda}\,
\frac{\Gamma(\frac{1+n}{4}-\frac{i\lambda}{2}+\frac{\varepsilon}{2})}
{\Gamma(\frac{1-n}{4}+\frac{i\lambda}{2}+\frac{\varepsilon}{2})}\
|\langle x,\,y\rangle|_{\varepsilon}^{-\frac{n+1}{2}+i\lambda}\,:
\end{multline}
hence
\begin{equation}
({\mathcal F}\,v)^{\flat}_{-i\lambda,\varepsilon}(x)=
\frac{1}{4\pi}\,(-i)^{\varepsilon}\,\pi^{-\frac{n}{2}+i\lambda}\,
\frac{\Gamma(\frac{1+n}{4}-\frac{i\lambda}{2}+\frac{\varepsilon}{2})}
{\Gamma(\frac{1-n}{4}+\frac{i\lambda}{2}+\frac{\varepsilon}{2})}\
\int_{\R^{n+1}}\,|\langle
x,\,y\rangle|_{\varepsilon}^{-\frac{n+1}{2}+i\lambda}\,v(y)\,dy
\end{equation}
and, for \,$\sigma\in \R$\,, setting \,$y_*=(y_1,\dots,y_n)$\,,
\begin{multline}
(\theta_{i\lambda,\varepsilon}\,v^{\flat}_{i\lambda,\varepsilon})(\sigma)=
({\mathcal F}\,v)^{\flat}_{-i\lambda,\varepsilon}(\sigma,\,1)=
(-i)^{\varepsilon}\,\pi^{-\frac{n}{2}+i\lambda}\,
\frac{\Gamma(\frac{1+n}{4}-\frac{i\lambda}{2}+\frac{\varepsilon}{2})}
{\Gamma(\frac{1-n}{4}+\frac{i\lambda}{2}+\frac{\varepsilon}{2})}
\ \times\\
 \frac{1}{4\pi}\int_{\R^n}dy_*\intR |\langle \sigma,\,y_*\rangle+y_{n+1}|
_{\varepsilon}^{-\frac{n+1}{2}+i\lambda}\
v(y_*,\,y_{n+1})\,dy_{n+1}\,:
\end{multline}
set \,$y_*=y_{n+1}\,s$\,, transforming the second line into
\begin{multline}
 \frac{1}{4\pi}\int_{\R^n} |1+\langle s,\,\sigma\rangle|_{\varepsilon}^{-\frac{n+1}{2}+i\lambda}\,
 ds\intR |y_{n+1}|_{\varepsilon}^{\frac{n-1}{2}+i\lambda}
v(y_{n+1}\,(s,1))\,dy_{n+1}=\\
\int_{\R^n} |1+\langle
s,\,\sigma\rangle|_{\varepsilon}^{-\frac{n+1}{2}+i\lambda}\,
v^{\flat}_{i\lambda,\varepsilon}(s)\,ds\,,
\end{multline}
as seen after another application of (\ref{2lambdapart}). Thus, the
(formal) definition of the intertwining operator
\,$\theta_{i\lambda,\varepsilon}$\, is finally
\begin{equation}\label{2thetaugly}
(\theta_{i\lambda,\varepsilon}\,u)(\sigma)=C_{i\lambda,\varepsilon}\,
\int_{\R^n} |1+\langle
s,\,\sigma\rangle|_{\varepsilon}^{-\frac{n+1}{2}+i\lambda}\,
u(s)\,ds\,,
\end{equation}
with
\begin{equation}\label{2Cilambda}
C_{i\lambda,\varepsilon}=(-i)^{\varepsilon}\,\pi^{-\frac{n}{2}+i\lambda}\,
\frac{\Gamma(\frac{1+n}{4}-\frac{i\lambda}{2}+\frac{\varepsilon}{2})}
{\Gamma(\frac{1-n}{4}+\frac{i\lambda}{2}+\frac{\varepsilon}{2})}\,,
\end{equation}
just the classical expression of the intertwining operator as
obtained from the general theory \cite{kna}: the extra phase factor
\,$(-i)^{\varepsilon}\,\pi^{i\lambda}$\,, not necessary for
unitarity, could be dispensed with in the present context, but is
important \cite{aumod} in
modular form theory, where it plays a role in the functional equations.\\

We are now in a position to introduce the
\,$(\lambda,\,\varepsilon)$\,-- dependent pseudodifferential
analysis of operators on \,$L^2(\R^n)$\, to be considered in the
present paper. Starting with an operator \,$A$\, on
\,$L^2(\R^{n+1})$\, commuting with the transformations \,$x\mapsto
tx,\ t\neq 0$\,, of the argument, so that \,$A$\, preserves the
parity of functions and transforms homogeneous functions into
functions homogeneous of the same degree, we restrict the operator
\,$A$\, to the space of functions homogeneous of a given degree and
parity \,$(-\frac{n+1}{2}-i\lambda,\,\varepsilon)$\,, identified
with the help of the map
\,$h\mapsto h_{i\lambda,\varepsilon}^{\flat}$\, to a space of functions on the projective space.\\

The following three spaces play a role here:\\
(i) the space \,$\Omega=\{(x,\,\xi)\in
\R^{n+1}\,\times\,\R^{n+1}\colon
\langle x,\,\xi\rangle \neq 0\}$\,;\\
(ii) the quotient \,$\Omega^{\bullet}$\, of \,$\Omega$\, under the
equivalence that identifies \,$(x,\,\xi)$\, to
\,$(tx,\,t^{-1}\xi)$\, for every \,$t\neq 0$\,; the image of
\,$(x,\,\xi)\in \Omega$\, under the associated canonical projection
is denoted as \,$(x,\,\xi)^{\bullet}$\,, not to be confused with
\,$(x^{\bullet},\,\xi^{\bullet})$\,;\\
(iii) finally, the subset \,${\mathcal X}_n^{\bullet}$\, of
\,$\Omega^{\bullet}$\, consisting of all points
\,$(x,\,\xi)^{\bullet}$\, such that \ $\langle x,\,\xi\rangle=1$\,;
\,${\mathcal X}_n^{\bullet}$\, can be identified with the quotient
of the hypersurface of \,$\Omega$\, of equation
\ $\langle x,\,\xi\rangle =1$\,  under the same equivalence as in (ii).\\

The space \,${\mathcal X}_n^{\bullet}$\, can also be identified with
the subset of \ $P_n(\R)\,\times\,P_n(\R)$\, consisting of all
points \,$(x^{\bullet},\,\xi^{\bullet})$\, such that \ $\langle
x,\,\xi\rangle \neq 0$\,, under the embedding \
$(x^{\bullet},\,\xi^{\bullet})\mapsto \l(x\,;\,\frac{\xi}{\langle
x,\,\xi\rangle}\r)^{\bullet}$\, of this latter space into
\,$\Omega^{\bullet}$\,. In terms of the (almost always defined only)
inhomogeneous coordinates \,$(s,\,\sigma)$\,, this embedding takes
the form
\begin{equation}\label{2ssigma}
(s,\,\sigma)\mapsto (s,1\,;\,\frac{\sigma}{1+\langle
s,\,\sigma\rangle},\frac{1}{1+\langle s,\,\sigma\rangle})\,.
\end{equation}
Finally, the space \,${\mathcal X}_n^{\bullet}$\, can be thought of
as the coset space \,$G_n/H_n$\,, where \,$G_n=SL(n+1,\R)$\, and
\,$H_n$\, is a subgroup of \,$G_n$\, isomorphic to \,$GL(n,\R)$\,,
to wit that made up by the linear transformations that respect the
splitting \ $\R^{n+1}=(\R^n\times\{0\})\, \oplus\,(\{0\}\times
\R)$\,. An invariant measure on \,${\mathcal X}_n^{\bullet}$\,
expresses itself, in terms of the coordinates above, as \
$|1+\langle s,\sigma\rangle|^{-n-1}\,ds\,d\sigma$\,.
Let us mention at once that taking quotients under this equivalence will be a fixture of what follows.\\

Recall from the beginning of the introduction that the standard
symbolic calculus \,${\mathrm{Op}}_0$\, and the antistandard
symbolic calculus \,${\mathrm{Op}}_1$\, on \,$\R^{n+1}$\, are
defined by the formulas
\begin{equation}\label{2stand}
({\mathrm{Op}}_0(H)\,v)(x)=\int_{\R^{n+1}}\,H(x;\,\xi)\,e^{2i\pi\langle
x,\,\xi\rangle}\,({\mathcal F}\,v)(\xi)\,d\xi
\end{equation}
and
\begin{equation}\label{2antistand}
({\mathcal
F}\,{\mathrm{Op}}_1(H)\,v)(\xi)=\int_{\R^{n+1}}\,H(x;\,\xi)\,e^{-2i\pi\langle
x,\,\xi\rangle}\,v(x)\,dx\,.
\end{equation}
In complete analogy, only replacing the Fourier transformation and
its integral kernel  by the operator
\,$\theta_{i\lambda,\varepsilon}$\,  and its integral kernel, one
introduces two species of symbols in the
\,$(\lambda,\varepsilon)$--dependent pseudodifferential calculus on
the projective space. The standard symbol \,$f$\, and the
antistandard symbol \,$h$\, of some operator \,$A$\, are the
functions such that \ $A={\mathrm{Op}}_{i\lambda,\varepsilon}(f)$\,
or
\ $A={\mathrm{Op}}_{i\lambda,\varepsilon}^{\vee}(h)$\, according to the definitions that follow.\\

\begin{definition}
The standard and antistandard symbolic calculi associated with the
pair \,$(\lambda,\varepsilon)$ are defined by the equations
\begin{equation}\label{2Op}
\left({\mathrm{Op}}_{i\lambda,\varepsilon}(f)u\right)(s)=(-1)^{\varepsilon}\,
C_{-i\lambda,\varepsilon}\,\int f(s,\sigma)\,|1+\langle
s,\sigma\rangle|_{\varepsilon}
^{-\frac{n+1}{2}-i\lambda}\left(\theta_{i\lambda,\varepsilon}\,u\right)(\sigma)\,d\sigma
\end{equation}
and
\begin{equation}\label{2Opvee}
(\theta_{i\lambda,\varepsilon}\,{\mathrm{Op}}_{i\lambda,\varepsilon}^{\vee}(h)\,u)(\sigma)=
C_{i\lambda,\varepsilon}\,\int h(s,\sigma)\,|1+\langle
s,\sigma\rangle|_
{\varepsilon}^{-\frac{n+1}{2}+i\lambda}\,u(s)\,ds\,.
\end{equation}\\
\end{definition}

Since the intertwining operator \,$\theta_{i\lambda,\varepsilon}$\,
is unitary, either defining map, after it has been divided by the
constant in front of the integral that defines it, sets up an
isometry between the Hilbert space \,$L^2({\mathcal
X}_n^{\bullet},\,|1+\langle s,\sigma\rangle|^{-n-1}\,ds\,d\sigma)$\,
and the space of Hilbert--Schmidt operators on \,$L^2(\R^n)$\,. It
is also immediate that the adjoint of the operator \
${\mathrm{Op}}_{i\lambda,\varepsilon}(f)$\, is the operator
\,${\mathrm{Op}}_{i\lambda,\varepsilon}^{\vee}(\bar{f})$\,. The
normalisation constants in front of the two integrals have been
introduced so that the (standard or antistandard) symbol of the
identity operator should be the constant \,$1$. An immediate, purely
formal, property of the
\,${\mathrm{Op}}_{i\lambda,\varepsilon}$--\,symbolic calculus is its
covariance under the representation \,$\pi_{i\lambda,\varepsilon}$\,
and the action defined by
\begin{equation}
g.\,(s,\,\sigma)=([g]\,s,\,[{g'}^{-1}]\,\sigma)
\end{equation}
of \,$G_n$\, in \,$P_n(\R)\times P_n(\R)$\,: this means that, for
every \,$g\in G_n$\,, one has the equation
\begin{equation}\label{2covariance}
    \pi_{i\lambda,\varepsilon}(g)\,{\mathrm{Op}}_{i\lambda,\varepsilon}(f)\,
    \pi_{i\lambda,\varepsilon}(g^{-1})=
    {\mathrm{Op}}_{i\lambda,\varepsilon}(f\circ g^{-1})\,.
\end{equation}
The same holds with the
\,${\mathrm{Op}}^{\vee}_{i\lambda,\varepsilon}$\,--symbolic calculus.\\

We now connect the
\,${\mathrm{Op}}_{i\lambda,\varepsilon}$\,--calculus ({\em resp.\/}
the \,${\mathrm{Op}}^{\vee}_{i\lambda,\varepsilon}$\,--calculus) to
the standard ({\em resp.\/} antistandard) calculus of operators on
functions on \,$\R^{n+1}$\,. It is necessary to consider symbols
\,$H=H(x,\,\xi)$\, invariant under transformations \
$(x,\,\xi)\mapsto (tx,\,t^{-1}\xi)\,,\ t\in \R^{\times}$\,: this
condition means that the associated operator (from the standard or
antistandard calculus) commutes with the transformations \ $x\mapsto
tx,\ t\neq 0$\,, of the argument. There is a slight difficulty, in
relation with the fact that the invariance of \,$H$\, does not
permit it to satisfy estimates (relative to its derivatives) of any
kind usual in pseudodifferential analysis: a very crude analysis,
however, will be sufficient for our purpose. To start with, if
\,$H$\, is bounded, \,${\mathrm{Op}}_0(H)$\, sends the space
\,${\mathcal S}(\R^{n+1})$\, into the space \,${\mathcal B}$\, of
continuous bounded functions, and \,${\mathrm{Op}}_1(H)$\, sends the
space \,${\mathcal S}(\R^{n+1})$\, into the image, under the Fourier
transformation, of \,${\mathcal B}$\,. This is not yet satisfactory
since, for the analysis to follow, we need to end up in the space
\,$L^2(\R^{n+1})$\,. To that effect, let us assume that the symbol
\,$H$\, is \,$C^{\infty}$\, and bounded, and that it remains in this
space after it has been applied any operator in the algebra
generated by the differential operators \
$x_j\,\frac{\partial}{\partial x_k}$\, or \
$\xi_j\,\frac{\partial}{\partial \xi_k}$\,: this condition is
compatible with the invariance of \,$H$\,. Let \,${\mathcal
S}_{\mathrm{flat}}(\R^{n+1})$\, be the already mentioned space of
all  rapidly decreasing \,$C^{\infty}$\, functions on
\,$\R^{n+1}$\,, flat at the origin. Then, an elementary integration
by parts shows that the operator \,${\mathrm{Op}}_0(H)$\, sends the
image \,${\mathcal F}\,{\mathcal S}_{\mathrm{flat}}(\R^{n+1})$\, of
the space \,${\mathcal S}_{\mathrm{flat}}(\R^{n+1})$\, under the
Fourier transformation into the space (contained in
\,$L^2(\R^{n+1})$) of continuous functions which remain bounded
after they have been multiplied by any polynomial in \,$x$\,. Also,
the operator \,${\mathrm{Op}}_1(H)$\, sends the space \,${\mathcal
S}_{\mathrm{flat}}(\R^{n+1})$\, into \,$L^2(\R^{n+1})$\,: of course,
both spaces \,${\mathcal S}_{\mathrm{flat}}(\R^{n+1})$\, and
\,${\mathcal F}\,{\mathcal S}_{\mathrm{flat}}(\R^{n+1})$\, are dense in \,$L^2(\R^{n+1})$\,.\\

\begin{proposition}
Let \,$H=H(x,\,\xi)$\, be a symbol in the space
\,$C^{\infty}(\R^{n+1}\times \R^{n+1})$\,, invariant under the
transformations \,$(x,\,\xi)\mapsto (tx,\,t^{-1}\xi)$\,, \,$t\neq
0$\,: assume that it is bounded and remains so after it has been
applied any operator in the algebra generated by the differential
operators \ $x_j\,\frac{\partial}{\partial x_k}$\, or \
$\xi_j\,\frac{\partial}{\partial \xi_k}$\,. Let
\,${\mathrm{Op}}_0(H)$\, be the pseudodifferential operator: \,$
{\mathcal F}\,{\mathcal S}_{\mathrm{flat}}(\R^{n+1})\to
L^2(\R^{n+1})$\, with standard symbol \,$H$\,. For every function
\,$v\in {\mathcal S}(\R^{n+1})$\, and every
\,$(\lambda,\,\varepsilon)$\,, one has
\begin{equation}\label{2OpHpres}
({\mathrm{Op}}_0(H)\,v)_{i\lambda,\varepsilon}={\mathrm{Op}}_0(H)\,v_{i\lambda,\varepsilon}\,.
\end{equation}
For any given pair \,$(\lambda,\,\varepsilon)$\,, the operator
\,$A_{i\lambda,\varepsilon}$\, on functions of \,$n$ variables
characterized by the property
\begin{equation}\label{2Alambeps}
({\mathrm{Op}}_0(H)\,v)_{i\lambda,\varepsilon}^{\flat}=
A_{i\lambda,\varepsilon}\,v_{i\lambda,\varepsilon}^{\flat}
\end{equation}
can be identified with the operator \
${\mathrm{Op}}_{i\lambda,\varepsilon}(f)$\, if one defines the
function \,$f$\, on \,${\mathcal X}_n^{\bullet}$\, in terms of \
$(H,\,\lambda,\varepsilon)$\, by
\begin{equation}\label{2compuf}
f((x,\,\xi)^{\bullet})=(-1)^{\varepsilon}\,C_{-i\lambda,\varepsilon}^{-1}\intR
H\l((x,\,t\xi)^{\bullet}\r)\ e^{2i\pi t} \
|t|_{\varepsilon}^{\frac{n-1}{2}+i\lambda}\,dt\,.
\end{equation}
In a similar way, the operator \,$B_{i\lambda,\varepsilon}$\,
characterized by the equation
\begin{equation}\label{2Alambepsbis}
({\mathrm{Op}}_1(H)\,v)_{i\lambda,\varepsilon}^{\flat}=
B_{i\lambda,\varepsilon}\,v_{i\lambda,\varepsilon}^{\flat}
\end{equation}
can be written as  \
$B_{i\lambda,\varepsilon}={\mathrm{Op}}_{i\lambda,\varepsilon}^{\vee}(h)$\,
with
\begin{equation}\label{2compuh}
h((x,\,\xi)^{\bullet})=C_{i\lambda,\varepsilon}^{-1}\intR
H\l((tx,\,\xi)^{\bullet}\r)\ e^{-2i\pi t} \
|t|_{\varepsilon}^{\frac{n-1}{2}-i\lambda}\,dt\,.
\end{equation}\\
\end{proposition}

\begin{proof}
That the symbol \,$H$\, is even just means that the operator
\,${\mathrm{Op}}(H)$\, preserves the parity of functions; that it is
invariant under the one--parameter group of transformations as
defined above with \,$t>0$\, means that the associated operator
preserves the space of homogeneous functions of any given degree.
The equation (\ref{2OpHpres}) follows: we now make the operator
\,$A_{i\lambda,\varepsilon}$\, explicit, letting
\,${\mathrm{Op}}_0(H)$\,  act on a function \,$v$\, in the space
\,${\mathcal F}\,{\mathcal S}_{\mathrm{flat}}(\R^{n+1})$\,. Setting
\ $\xi=(\xi_*,\,\xi_{n+1})$, noting that \ ${\mathcal
F}\,v_{i\lambda,\varepsilon}= ({\mathcal
F}\,v)_{-i\lambda,\varepsilon}$\, and using (\ref{2charac}) and
(\ref{2deftheta}), we find
\begin{multline}
({\mathcal F}\,v_{i\lambda,\varepsilon})(\xi_*,\,\xi_{n+1})=
({\mathcal F}\,v)_{-i\lambda,\varepsilon}(\xi_*,\,\xi_{n+1})=\\
|\xi_{n+1}|_{\varepsilon}^{-\frac{n+1}{2}+i\lambda}\, ({\mathcal
F}\,v)^{\flat}_{-i\lambda,\varepsilon}\,(\frac{\xi_*}{\xi_{n+1}})=
|\xi_{n+1}|_{\varepsilon}^{-\frac{n+1}{2}+i\lambda}\,(\theta_{i\lambda,\varepsilon}\,
v^{\flat}_{i\lambda,\varepsilon})\,(\frac{\xi_*}{\xi_{n+1}})\,:
\end{multline}
since, using the definition (\ref{2stand}) of
\,${\mathrm{Op}}_0(H)$\,,
\begin{align}
({\mathrm{Op}}_0(H)\,v)_{i\lambda,\varepsilon}^{\flat}\,(s)&=
({\mathrm{Op}}_0(H)\,v)_{i\lambda,\varepsilon}\,(s,1)\nonumber\\
&=\int_{\R^{n+1}} H(s,1\,;\,\xi)\,e^{2i\pi\,(\langle
s,\,\xi_*\rangle+\xi_{n+1})}\, ({\mathcal
F}\,v_{i\lambda,\varepsilon})(\xi)\,d\xi\,,
\end{align}
one sees from the equation (\ref{2Op}) that one has \
$A_{i\lambda,\varepsilon}= {\mathrm{Op}}_{i\lambda,\varepsilon}(f)$
\  provided that one defines
\begin{align}\label{2f}
f(s,\,\sigma)&=(-1)^{\varepsilon}\,C_{-i\lambda,\varepsilon}^{-1}\,
|1+\langle s,\sigma\rangle|_{\varepsilon}^{\frac{n+1}{2}+i\lambda}\,
\intR H(s,1\,;\,t\sigma,t)\ e^{2i\pi t\,(1+\langle
s,\sigma\rangle)}\,
|t|_{\varepsilon}^{\frac{n-1}{2}+i\lambda}\,dt\nonumber\\
&=(-1)^{\varepsilon}\,C_{-i\lambda,\varepsilon}^{-1} \intR
H(s,1\,;\,\frac{t\,\sigma}{1+\langle
s,\,\sigma\rangle},\frac{t}{1+\langle s,\,\sigma\rangle})\ e^{2i\pi
t}\, |t|_{\varepsilon}^{\frac{n-1}{2}+i\lambda}\,dt\,,
\end{align}
an expression which can be identified with (\ref{2compuf}).\\

To arrive at the computation of the function \,$h$\, such that \
$B_{i\lambda,\varepsilon}=
{\mathrm{Op}}^{\vee}_{i\lambda,\varepsilon}(h)$\, according to
(\ref{2Opvee}), we write, using in succession
(\ref{2Alambepsbis}),\,(\ref{2deftheta}),\,(\ref{2defflat}) and
(\ref{2OpHpres}), and starting this time from a function \ $v\in
{\mathcal S}_{\mathrm{flat}}(\R^{n+1})$\,,
\begin{multline}
(\theta_{i\lambda,\varepsilon}\,B_{i\lambda,\varepsilon}\,v_{i\lambda,\varepsilon}^{\flat})(\sigma)=
(\theta_{i\lambda,\varepsilon}\,({\mathrm{Op}}_1(H)\,v)_{i\lambda,\varepsilon}^{\flat})(\sigma)=
(({\mathcal F}\,{\mathrm{Op}}_1(H)\,v)_{-i\lambda,\varepsilon}^{\flat})(\sigma)=\\
(({\mathcal
F}\,{\mathrm{Op}}_1(H)\,v)_{-i\lambda,\varepsilon})(\sigma,1)=
({\mathcal
F}\,(({\mathrm{Op}}_1(H)\,v)_{i\lambda,\varepsilon}))(\sigma,1)=
({\mathcal
F}\,({\mathrm{Op}}_1(H)\,v_{i\lambda,\varepsilon}))(\sigma,1)\,:
\end{multline}
we then use the definition (\ref{2antistand}) of
\,${\mathrm{Op}}_1(H)$\,, expressing what precedes as
\begin{multline}
\int_{\R^{n+1}} H(x_*,x_{n+1};\,\sigma,1)\,e^{-2i\pi\,[\langle
x_*,\,\sigma\rangle
+x_{n+1}]}\,v_{i\lambda,\varepsilon}(x_*,\,x_{n+1})\,dx_*\,dx_{n+1}\\
=\int_{\R^{n+1}} H(\frac{ts}{1+\langle
s,\,\sigma\rangle},\frac{t}{1+\langle s,\,\sigma\rangle};\,
\sigma,1)\,e^{-2i\pi
t}\,|t|_{\varepsilon}^{\frac{n-1}{2}-i\lambda}\,v^{\flat}_{i\lambda,\varepsilon}(s)
\ |1+\langle
s,\,\sigma\rangle|_{\varepsilon}^{-\frac{n+1}{2}+i\lambda}\
ds\,dt\,.
\end{multline}
We must thus take, this time,
\begin{equation}\label{2h}
h(s,\,\sigma)=C_{i\lambda,\varepsilon}^{-1}\intR
H(\frac{ts}{1+\langle s,\,\sigma\rangle},\frac{t}{1+\langle
s,\,\sigma\rangle};\,\sigma,1)\ e^{-2i\pi t}\,
|t|_{\varepsilon}^{\frac{n-1}{2}-i\lambda}\,dt\,,
\end{equation}
which leads finally to (\ref{2compuh}).\\
\end{proof}

A fully equivalent, more expressive way to write the function
\,$f$\, or \,$h$\, on \,${\mathcal X}_n^{\bullet}$\, is as the
restriction to \,${\mathcal X}_n^{\bullet}$\,  of the image of
\,$H$\, under some operator expressed, in spectral--theoretic terms,
as a function of the pair \ $({\langle \xi,\,
\frac{\partial}{\partial \xi}\rangle,\,\delta})$\,, where the first
entry denotes the Euler operator \,$\langle \xi,\,
\frac{\partial}{\partial
\xi}\rangle=\sum_{k=1}^{n+1}\xi_k\,\frac{\partial}{\partial
\xi_k}$\,, and the second one is the usual indicator of parity of
functions {\em with respect to\/} \,$\xi$\, only. Note that the
dissymmetry between the variables \,$x,\,\xi$\, is only apparent
since, as \,$H$\, lives on \,$\Omega^{\bullet}$\,, it satisfies the
equation \ $\sum_{k=1}^{n+1} (x_k\,\frac{\partial}{\partial x_k}-
\xi_k\,\frac{\partial}{\partial \xi_k})\,H=0$\,; also, it is {\em
globally\/} even as
a function of \,$(x,\,\xi)$\,.\\

\begin{proposition}
Under the assumptions of the preceding proposition, observe that the
number \,$|\varepsilon-\delta|$\, is the number, equal to \,$0$ or
$1$, characterized by the congruence \ $|\varepsilon-\delta|\equiv
\varepsilon+\delta$\, mod $2$\,, and set
\begin{align}
D_1&=i^{-\varepsilon+|\varepsilon-\delta|}\,\pi^{-\langle
\xi,\,\frac{\partial}{\partial \xi} \rangle}\
\frac{\Gamma(\frac{1-n}{4}-\frac{i\lambda}{2}+\frac{\varepsilon}{2})}
{\Gamma(\frac{1+n}{4}+\frac{i\lambda}{2}+\frac{\varepsilon}{2})}\,
\frac{\Gamma(\frac{n+1}{4}+\hf\,\langle
\xi,\,\frac{\partial}{\partial \xi}\rangle+
\frac{i\lambda}{2}+\frac{|\varepsilon-\delta|}{2})}
{\Gamma(\frac{1-n}{4}-\hf\,\langle \xi,\,\frac{\partial}{\partial
\xi}\rangle
-\frac{i\lambda}{2}+\frac{|\varepsilon-\delta|}{2})}\,,\nonumber\\
D_2&=i^{\varepsilon-|\varepsilon-\delta|}\,\pi^{-\langle
\xi,\,\frac{\partial}{\partial \xi} \rangle}\
\frac{\Gamma(\frac{1-n}{4}+\frac{i\lambda}{2}+\frac{\varepsilon}{2})}
{\Gamma(\frac{1+n}{4}-\frac{i\lambda}{2}+\frac{\varepsilon}{2})}\,
\frac{\Gamma(\frac{n+1}{4}+\hf\,\langle
\xi,\,\frac{\partial}{\partial \xi}\rangle-
\frac{i\lambda}{2}+\frac{|\varepsilon-\delta|}{2})}
{\Gamma(\frac{1-n}{4}-\hf\,\langle \xi,\,\frac{\partial}{\partial
\xi}\rangle +\frac{i\lambda}{2}+\frac{|\varepsilon-\delta|}{2})}\,.
\end{align}
Then one has
\begin{equation}
f=D_1\,H\,\biggr|\!\biggr|_{{\mathcal X}_n^{\bullet}}\,,\qquad
\qquad h=D_2\,H\,\biggr|\!\biggr| _{{\mathcal X}_n^{\bullet}}\,:
\end{equation}
the double restriction bar indicates that one should first restrict
the function under consideration to the hypersurface of
\,$\R^{n+1}\times \R^{n+1}$\, of equation \ $\sc{x,\,\xi}=1$\,, next
use the invariance of the result under the usual action of
\,$\R^{\times}$\, to make it a function on the
corresponding quotient \,${\mathcal X}_n^{\bullet}$\, of this hypersurface.\\
\end{proposition}

\begin{proof}
On functions \,$H$\, with the parity \,$\delta$\, with respect to
\,$\xi$\,, one may write the function \ $(x,\xi)\mapsto H(x,t\xi)$\,
as \ $|t|_{\delta}^{\langle \xi,\,\frac{\partial}{\partial
\xi}\rangle}\,H$\,: after having inserted a factor
\,$e^{-\alpha|t|}$\, for convergence and letting \,$\alpha>0$\, go
to zero, one may use the Fourier transformation formula
\begin{equation}
\intR |t|_{\delta}^X\,e^{2i\pi
t}\,|t|_{\varepsilon}^{\frac{n-1}{2}+i\lambda}\,dt=
i^{|\varepsilon-\delta|}\,\pi^{-X-\frac{n}{2}-i\lambda} \
\frac{\Gamma(\frac{n+1}{4}+\frac{X}{2}+\frac{i\lambda}{2}+\frac{|\varepsilon-\delta|}{2})}
{\Gamma(\frac{1-n}{4}-\frac{X}{2}-\frac{i\lambda}{2}+\frac{|\varepsilon-\delta|}{2})}\,,
\end{equation}
which leads to the desired result, starting from (\ref{2f}) or (\ref{2h}).\\
\end{proof}

Functions \,$H=H(x,\,\xi)$\, on \,$\Omega$\, invariant under the
transformations \,$(x,\,\xi)\mapsto (tx,\,t^{-1}\xi)$\,, \,$t\neq
0$\,, can also be written as functions \ $h=h(s,\,\sigma;\,q)$\, on
\,$\R^n\times \R^n\times\R$\,, or on \,${\mathcal X}_n^{\bullet}
\times \R$\, ({\em cf.\/}\,(\ref{2ssigma})), under the
correspondence
\begin{equation}
H(x,\,\xi)=h\,(\frac{x_*}{x_{n+1}},\,\frac{\xi_*}{\xi_{n+1}};\,\langle
x,\,\xi\rangle)\,,
\end{equation}
with \ $x=(x_*,\,x_{n+1})$\, and \ $\xi=(\xi_*,\,\xi_{n+1})$\,.
Then, one has
\begin{multline}\label{2bigequ}
\langle
x,\,\xi\rangle\,\sum_{k=1}^{n+1}\,\frac{\partial^2\,H}{\partial
x_k\,\partial \xi_k}= (1+\langle
s,\,\sigma\rangle)\,\l[\,\sum_{j=1}^{n}\,\frac{\partial^2}{\partial
s_j\,\partial \sigma_j} +(\sum_{\ell=1}^n
s_{\ell}\,\frac{\partial}{\partial s_{\ell}})\,
(\sum_{m=1}^n \sigma_m\,\frac{\partial}{\partial \sigma_m})\,\r]\,h\\
+(n+1)\,q\,\frac{\partial h}{\partial q}+q^2\,
\frac{\partial^2 h}{\partial q^2}\,.\\
\end{multline}

This brings to light the operator
\begin{equation}\label{2Deltan}
\Delta_n= (1+\langle s,\sigma \rangle)\,\l[\,\sum_{j=1}^n
\frac{\partial^2}{\partial s_j\,
\partial \sigma_j}+(\sum s_{\ell}\,\frac{\partial}{\partial s_{\ell}})\,
(\sum \sigma_m\,\frac{\partial}{\partial \sigma_m})\,\r]\,,
\end{equation}
the fundamental invariant differential operator on the
(non--Riemannian) symmetric space \,${\mathcal
X}_n^{\bullet}=G_n/H_n$\,. One also sees that, if \,$H$\, lies in
the null space of the operator \
$\sum_{k=1}^{n+1}\frac{\partial^2}{\partial x_k\,\partial \xi_k}$\,,
and if \,$H$\, is, with respect to the variable \,$\xi$\, only,
homogeneous of degree \,$\rho$\,, its restriction \,$h$\, to
\,${\mathcal X}_n^{\bullet}$\, satisfies the eigenvalue equation
\begin{equation}\label{2eigen}
\Delta_n\, h=-\rho\,(n+\rho)\,h\,.
\end{equation}\\

An example is provided by the function \ $H(x,\,\xi)=|\langle
a,\,x\rangle\, \langle b,\,\xi\rangle|_{\delta}^{\rho}$\, on
\,$\Omega^{\bullet}$\, with \ $\langle a,\,b\rangle=0$\,, which
gives rise to the function \,$|\phi_{a,b}|_{\delta}^{\rho}$\, on
\,${\mathcal X}_n$\, with
\begin{equation}\label{2phiab}
\phi_{a,b}(s,\,\sigma)=\frac {(a_{n+1}+\langle
a_*,\,s\rangle)\,(b_{n+1}+\langle b_*,\,\sigma\rangle)} {1+\langle
s,\,\sigma\rangle}\,.
\end{equation}\\

As shown in \cite{van2,vdm}, the quasiregular representation of
\,$G_n$\, in \,$L^2({\mathcal X}_n^{\bullet})=L^2(G_n/H_n)$\,
decomposes into a continuous part and a discrete part, a fact
tantamount to the analogous statement regarding \,$\Delta_n$\,. We
shall take it, temporarily, for granted that --- a consequence of
the analysis to be developed  in the next section --- functions on
\,${\mathcal X}_n^{\bullet}$\, in the continuous part of the
decomposition can always be viewed (in many ways) as restrictions to
\,${\mathcal X}_n^{\bullet}$\, of \,$\R^{\times}$--invariant
functions satisfying in \,$\R^{n+1}\times \R^{n+1}$\,, in the
distribution sense, the equation
\begin{equation}\label{2J}
\square\,H\colon=\sum_{k=1}^{n+1} \frac{\partial^2\,H}{\partial
x_k\,\partial \xi_k}=0\,:
\end{equation}
this condition is remarkable from the point of view of
pseudodifferential analysis since it means that the operator
\,${\mathrm{Op}}_{\kappa}(H)$\,, as defined in (\ref{1Op}), does not
depend on \,$\kappa$\,. This follows from the equation ({\em cf.\/}
\cite[p.\,15]{aar} or do elementary manipulations using the Fourier
transformation) \ ${\mathrm{Op}}_{\kappa}(H)=
{\mathrm{Op}}_0\,\l(\exp\,(\frac{\kappa}{2i\pi}\, \sum_{k=1}^{n+1}
\frac{\partial^2}{\partial x_k\,\partial \xi_k})\,H\r)$\,: in
particular,
the operators with standard or antistandard symbol \,$H$\, are identical.\\

The following corollary is one half (that concerning the continuous
part \,$L^2_{\mathrm{cont}}(G_n/H_n)$\, of the decomposition of
\,$L^2(G_n/H_n)$) of the last result of the paper \cite{vdm} by
Molchanov and van Dijk.

\begin{corollary}
The operator \,$J_{i\lambda,\varepsilon}\colon
L^2_{\mathrm{cont}}(G_n/H_n)\to L^2_{\mathrm{cont}}(G_n/H_n)$\,
defined by the validity of the equation \
${\mathrm{Op}}_{i\lambda,\varepsilon}^{\vee}(f)=
{\mathrm{Op}}_{i\lambda,\varepsilon}(J_{i\lambda,\varepsilon}\,f)$\,
for every \,$f\in L^2(G_n/H_n)$\, is characterized by the following
property: on functions which are generalized eigenfunctions of
\,$\Delta_n$\, for the eigenvalue \,$-\rho\,(n+\rho)$\, and have,
with respect to \,$\xi$\, only, the parity characterized by
\,$\delta$\,, the operator \,$J_{i\lambda,\varepsilon}$\, coincides
with the scalar

\begin{equation}\label{2formdel0}
G_{i\lambda,\varepsilon}\,(\rho,\delta)=(-1)^{\delta}\,
\frac{\Gamma(\frac{n+1-\mu+\varepsilon}{2})}
{\Gamma(\frac{-n+\mu+\varepsilon}{2})}\,\frac{\Gamma(\frac{-n+\mu-\rho+|\varepsilon-\delta|}{2})}
{\Gamma(\frac{n+1-\mu+\rho+|\varepsilon-\delta|}{2})}\,\frac{\Gamma(\frac{1-\mu+\varepsilon}{2})}
{\Gamma(\frac{\mu+\varepsilon}{2})}\,\frac{\Gamma(\frac{\mu+\rho+|\varepsilon-\delta|}{2})}
{\Gamma(\frac{1-\mu-\rho+|\varepsilon-\delta|}{2})}\,,
\end{equation}
where  \,$\mu=\frac{n+1}{2}+i\lambda$\,.\\
\end{corollary}

\begin{proof}
Apply Proposition 2.2 with a symbol \,$H$\, such that \
$\square\,H=0$\,, so that \
${\mathrm{Op}}_0(H)={\mathrm{Op}}_1(H)$\, and \
$A_{i\lambda,\varepsilon}= B_{i\lambda,\varepsilon}$\, for every
pair \,$(\lambda,\varepsilon)$\,. It then suffices to consider the
expression of \,$D_1\,D_2^{-1}$\, obtained as a consequence of
Proposition 2.3, and to use the equation (\ref{2eigen}): of course,
it should be noted that the expression on the right--hand side of
(\ref{2formdel0}) is invariant under the change
\ $\rho\mapsto -n-\rho$\,, which makes it a function of \,$\rho\,(n+\rho)$\,.\\
\end{proof}

\noindent {\bf Remark:}  It will be seen in Section 5 that Corollary
2.4 remains valid for certain symbols --- important from the point
of view of pseudodifferential analysis --- which are far from lying
in \,$L^2(G_n/H_n)$\,.


\section{The square root method: the continuous part of the operator \,$\Delta_n-\frac{n^2}{4}$}

It is not our intention to give a complete exposition of the
decomposition of the space \,$L^2(G_n/H_n)$\, under the quasiregular
action of \,$G_n$\,, already made in the references just recalled.
The present section justifies the assertion, made just before
Corollary 2.4, concerning the possibility to realize functions in
the continuous part of the decomposition of \,$L^2(G_n/H_n)$\, with
the help of solutions in \,$\R^{n+1}\times \R^{n+1}$\, of the
equation \ $\square\,H=0$\,. At the same time, it introduces a new
construction of this continuous part, which, to our taste at least,
makes the whole picture very clear (Theorem 3.2). Rather than
continuing with a study of the discrete subspaces of the
decomposition of \,$L^2(G_n/H_n)$\,, we shall follow, in section 5,
with the study of some interesting special distributions on
\,${\mathcal X}_n^{\bullet}$\, related to symbols of operators in
the algebra generated by resolvents of certain infinitesimal
operators of the representation \,$\pi_{i\lambda,\varepsilon}$\,.\\

The equation (\ref{2eigen}) shows that generalized eigenvalues of
the operator in the title of this section present themselves in the
form \,$-(\rho+\frac{n}{2})^2$\,: more to the point, the spectral
theory of this operator --- as developed in \cite{van2,vdm} with the
help of \,$H_n$-- \,spherical distribution theory --- shows that the
spectrum of \,$\Delta_n$\, has a continuous part, consisting of all
numbers \,$\rho(-n-\rho)$\, with \,$\rho=-\frac{n}{2}+ir\in
-\frac{n}{2}+i\R$\,, and a discrete part consisting of the numbers
\,$(\frac{1-n}{2}+k)(\frac{-1-n}{2}-k)$\, with \,$k\in \N$\,. Now,
the non--negative integer \,$k$\, is certainly uniquely determined
by this latter rational number, but making the choice, for every
\,$r^2>0$\,, of one of the two numbers \,$r$\, and \,$-r$\,, would
be a very rough way of defining a square root of the continuous part
of the operator \,$\Delta_n-\frac{n^2}{4}$\,. A better solution
consists in making, under some transfer, the latter operator appear
as the square of some {\em differential\/} operator, defined not on
\,${\mathcal X}_n^{\bullet}$\, but on some other space
\,${\Sigma}_n^{\bullet}$\,, to wit the cone of equation \ $\langle
x,\,\xi\rangle =0$\,, divided by the equivalence \ $(x.\,\xi)\sim
(tx,\,t^{-1}\xi) (t\in \R^{\times})$\,, in such a way that the
generalized
eigenvalue \,$r^2$\, should split there as the pair \,$\pm r$\,. \\

The way to do this, introduced in \cite[section 18]{aumod} in the
case of the simplest Riemannian symmetric space
\,$SL(2,\R)/SO(2)$\,, is an alternative approach  to the spectral
theory of the invariant operator \,$\Delta$\, under consideration
with several advantages: in the situation already experienced, it
led to some renewed understanding of the Lax--Phillips scattering
theory \cite{lap} for the automorphic wave equation, and proved
useful
in modular form theory.\\

Let us start from the following analogue of \cite[p.\,11]{lap}:
under the change of variables \ $(x,\,\xi)\mapsto
(s,\sigma;\,\tau)$\, from
\,$\Omega_+^{\bullet}=\{(x,\,\xi)^{\bullet}\in
\Omega^{\bullet}\colon \langle x,\,\xi \rangle > 0\}$ \ to
\,${\mathcal X}_n^{\bullet}\,\times\,\R$\, defined by the pair of
equations (in which \ $x=(x_*,\,x_{n+1}),\
\xi=(\xi_*,\,\xi_{n+1})$\,)
\begin{equation}\label{3coord}
(s,\,\sigma)=(\frac{x_*}{x_{n+1}},\,\frac{\xi_*}{\xi_{n+1}})\,,\qquad
\tau=\log\,\langle x,\,\xi \rangle\,,
\end{equation}
and under the transformation \,$H\mapsto
H_1=e^{\frac{n\tau}{2}}\,H$\,, the equation \,$\square H=0$\, inside
\,$\Omega_+^{\bullet}$\, is equivalent to the wave equation
\begin{equation}\label{3wave}
\frac{\partial^2H_1}{\partial
\tau^2}+(\Delta_n-\frac{n^2}{4})\,H_1=0\,.
\end{equation}
Indeed, this follows from the equation (\ref{2bigequ}) if one writes
\ $(n+1)\,q\, \frac{\partial h}{\partial q}+q^2\,
\frac{\partial^2 h}{\partial q^2}=(\frac{d}{d\tau}+\frac{n}{2})^2h-\frac{n^2}{4}\,h$\,.\\

Of course, a solution of this equation in \,$\Omega_+^{\bullet}$\,
can be characterized by its first two traces on the hyperplane
\,$\tau=0$\,: however, only the first one is of interest to us in
the present context. It does not change if one replaces \,$H$\, by
its transform \,$\widetilde{H}$\, under the inversion map:
\begin{equation}\label{3inversion}
({\mathrm{Inv}}\,H)(x,\xi)=\widetilde{H}(x,\xi)=(\langle
x,\xi\rangle)^{-n}\, H\l(\frac{x}{\langle
x,\xi\rangle},\,\frac{\xi}{\langle x,\xi\rangle}\r)\,.
\end{equation}
This distribution lies in the nullspace of the operator
\,$\square$\, within \,$\Omega_+$\, if \,$H$\, does, and is also
invariant under the transformations \ $(x,\xi)\mapsto
(tx,\,t^{-1}\xi),\ t\in \R^{\times}$\,, if \,$H$\, is. In
particular, the function \,$H$\, defined just before (\ref{2phiab})
has the same restriction (to be followed by the passage to functions
on a quotient set) to \,${\mathcal X}_n^{\bullet}$\, as the function
\,$\widetilde{H}(x,\xi)=\frac{H(x,\xi)} {|\langle
x,\xi\rangle|^{2\rho+n}}$\,: note in this case that the degree of
homogeneity of \,$\widetilde{H}$\, with respect to \,$\xi$\, is
\,$-n-\rho$\, instead of \,$\rho$\,, which was to be expected in
view of (\ref{2eigen}). We base our present study of the continuous
part of the decomposition of \,$L^2({\mathcal
X}_n^{\bullet})=L^2(G_n/H_n)$\, on the construction of a certain map
\,$\Theta$\, from functions on \,$\Sigma_n^{\bullet}$\, to
\,$\R^{\times}$--invariant functions on \,$\R^{n+1}\times
\R^{n+1}$\, in the nullspace of \,$\square$\,, together with an
involution \,${\mathcal K}$\, on the first space of functions, such
that the knowledge of the first trace of \,$H=\Theta\,\Phi$\, on
\,${\mathcal X}_n^{\bullet}$\, should be equivalent to that of the
\,${\mathcal K}$--invariant part
of \,$\Phi$\, on \,$\Sigma_n^{\bullet}$\,.\\

The role of \,$\Sigma_n^{\bullet}$\, is also clear from the
representation--theoretic point of view. Indeed, the irreducible
unitary components of the quasiregular representation of \,$G_n$\,
in \,$L^2({\mathcal X}_n)$\, are sub--representations of
representations induced from a parabolic subgroup of parabolic rank
2. Such representations can be realized on spaces of functions on
the homogeneous space \,$G/MN$\,, where \,$N$\, is the Heisenberg
group of  dimension \,$2n-1$\,, a space
isomorphic to the cone under consideration.\\

Our first task is thus to give an efficient construction of all
solutions of the wave equation \ $\square \,H=0$\,. This can be done
in at least two different ways. The first one --- following
\cite[section 18]{aumod} --- is based on an extension of the theory
of Riesz operators \cite{rie} or, more properly said, distributions
\cite{sch}, to the case of the operator \,$\square$\,: this was our
first choice during the preparation of this paper, and it provides
more information  than the following one, based on the use of the
Fourier transformation;
the latter one has the advantage of being more concise.\\

On \,$\R^{n+1}\times \R^{n+1}$\,, we shall use throughout the
``symplectic'' Fourier transformation \,${\mathcal
F}_{\mathrm{symp}}$\, defined by the equation
\begin{equation}\label{3symplec}
({\mathcal F}_{\mathrm{symp}}\,{\mathcal S})(x,\,\xi)=
\langle\,{\mathcal S}\,,\,(y,\eta)\mapsto e^{2i\pi\,(\langle
x,\,\eta\rangle-\langle y,\,\xi\rangle)}\,\rangle\,:
\end{equation}
note that, when using the symplectic Fourier transformation, one
should not consider the ``dual'' variables as conceptually distinct
from the main ones. Then, the symplectic Fourier transform of a
distribution \,${\mathfrak S}$\, satisfying the distribution
equation \ $\sc{y,\,\eta}\,{\mathfrak S}=0$\, lies in the nullspace of \,$\square$\,.\\

Consider the following example. Given \,$\rho\in\C$\, and
\,$\delta=0$\, or $1$\, with \,$\rho+\delta\neq -1,-3,\dots$\, and
\,$\rho-\delta\neq 0,2,\dots$\, , finally given \,${a}\in
\R^{n+1}\backslash\{0\}$\,, introduce the distribution
\,$M^{a}_{\rho,\delta}$\, on \,$\R^{n+1}$\, defined by the equation
\begin{equation}\label{3Marho}
\langle\,M^{a}_{\rho,\delta},\,\psi\,\rangle=\intR
\psi(r{a})\,|r|_{\delta}^{-\rho-1}\,dr\,,\qquad\quad \psi\in
{\mathcal S}(\R^{n+1})\,:
\end{equation}
this is a measure supported by the line \,$\R {a}$\, in the case
when \ $\Re \rho<0$\,, and a well--defined distribution, homogeneous
of degree \,$-n-1-\rho$\,, whenever the pair \,$(\rho,\,\delta)$\,
satisfies the above conditions. If we take for \,${\mathfrak S}$\,
the distribution \ $
M^b_{\rho,\delta}\,\otimes\,M^{a}_{\rho,\delta}$\,, and if we assume
that \ $\sc{a,\,b}=0$\,, so that the equation \
$\sc{y,\,\eta}\,{\mathfrak S}=0$\, should hold, we obtain \
${\mathcal F}_{\mathrm{symp}}\,{\mathfrak S}= {\mathcal
F}^{-1}\,M^{a}_{\rho,\delta}\,\otimes\,{\mathcal
F}\,M^b_{\rho,\delta}$\,, {\em i.e.,\/}
\begin{equation}\label{3Fou}
({\mathcal F}_{\mathrm{symp}}\,{\mathfrak
S})\,(x,\xi)=\pi^{1+2\rho}\,
\l[\,\frac{\Gamma(\frac{-\rho+\delta}{2})}
{\Gamma(\frac{\rho+1+\delta}{2})}\,\r]^2 \ |\langle
{a},\,x\rangle\langle b,\,\xi\rangle|_{\delta}^{\rho}\,:
\end{equation}
we thus get back to the function \,$|\phi_{a,b}|_{\delta}^{\rho}$\,
considered just
before (\ref{2phiab}).\\

We now need to introduce a few geometric objects. The set
\,$\overline{\Sigma}_n$\, is the hypersurface of \,$\R^{n+1}\times
\R^{n+1}$\, of equation \ $\langle y,\,\eta\rangle =0$\,: let us
warn the reader that, when dealing with functions defined on this
space, it is sometimes necessary to switch from the variables
\,$(y,\eta)$\, to the variables \,$(x,\xi)$\,, since this space
plays a role ``on both sides'' of the (symplectic) Fourier
transformation. Let \,$\Sigma_n$\, denote the open dense subset of
\,$\overline{\Sigma}_n$\, characterized by the conditions \,$y\neq
0,\,\eta\neq 0$\,. The space \,${\Sigma}_n^{\bullet}$\, is the
quotient of \,${\Sigma}_n$\, by the group of transformations \
$(y,\eta)\mapsto (ty,t^{-1}\eta)$\, with \,$t\in \R^{\times}$\,:
since this action of \,$\R^{\times}$\, occurs in a consistent way,
such notions as \,$\R^{\times}$--invariant functions, or quotients
by \,$\R^{\times}$\,, will always make reference to this action. On
the cone \,${\Sigma}_n$\,, we may use the (singular) coordinates
\,$(y,\,\eta_*)$\,, since \ $\eta_{n+1}=-\frac{\langle
y_*,\,\eta_*\rangle}{y_{n+1}}$\,, and the \,$GL(n+1,\R)$--invariant
measure \ $dm(y,\,\eta_*)=|y_{n+1}|^{-1}\,dy\,d\eta_*$\,. On its
quotient \,${\Sigma}_n^{\bullet}$\,, we shall use the coordinates
\,$(y_*,\,\eta_*)$\, corresponding to the orbit, under the action of
\,$\R^{\times}$\,, of the point \ $(y_*,1;\,\eta_*,-\langle
y_*,\,\eta_*\rangle)$\,. In terms of these coordinates, we set
\,$dm^{\bullet}(y_*,\,\eta_*)=dy_*\,d\eta_*$\,: it is easy to show
that this measure is invariant under the action of \,$GL(n+1,\R)$\,
on \,${\Sigma}_n^{\bullet}$\, coming from the action
\,$g.\,(y,\eta)=(gy,\,{g'}^{-1}\eta)$\, of this group on
\,$\Sigma_n$\,. Only, recall that in
the case when \,$g\in SL(n+1,\R)$\,, it can be written as  \ $g=\l(\begin{smallmatrix}  M  & \vec{p} \\
\vec{q}^{\dag} &m \end{smallmatrix}\r)$\,, where \,$\vec{p}$\, is a
column vector and \,$\vec{q}^{\dag}$\, is the transpose of the
column vector \,$\vec{q}$\,: the computation of
\,$g.\,(y_*,\,\eta_*)$\, is trivial in the case when
\,$\vec{q}=0$\,, and it suffices to consider the case when
\,$\vec{p}=0$\, and \,$M=I$\,, in which one has \ $
g.\,(y_*,\,\eta_*)=\l(\frac{y_*}{1+\langle
\vec{q},\,y_*\rangle}\,;\,(1+\langle \vec{q},\,y_*\rangle)\,(\eta_*+
\langle y_*,\,\eta_*\rangle\,\vec{q})\r)$\,. Note that on
\,${\mathcal X}_n^{\bullet}$\,, too, we may still use the
coordinates \,$(x_*,\,\xi_*)$\, corresponding to the orbit, under
the action of \,$\R^{\times}$\,, of the point \
$(x_*,1;\,\xi_*,\,1-\langle x_*,\,\xi_*\rangle)$\, and that, in
these coordinates, the \,$G_n$--invariant measure \,$|1+\langle
s,\sigma\rangle|^{-n-1}\,ds\,d\sigma$\,  reduces again to \,$dx_*\,d\xi_*$\,.\\

Any \,$C^{\infty}$\, function \ $\Phi=\Phi(y_*,\,\eta_*)$\, on
\,${\Sigma}_n^{\bullet}$\, can be extended as a function
\,$\widetilde{\Phi}$\, on \,${\Sigma}_n$\, in a natural way, setting
\
$\widetilde{\Phi}\,(y,\,\eta_*)=\Phi\,(\frac{y_*}{y_{n+1}},\,y_{n+1}\,\eta_*)$\,.
Then, the distribution \ ${\mathfrak S}=\widetilde{\Phi}\,dm$\, is
supported in \,$\overline{\Sigma}_n$\,, is
\,$\R^{\times}$--invariant, and satisfies the distribution equation
\ $\sc{y,\,\eta}\,{\mathfrak S}=0$\,, so that its symplectic Fourier
transform lies in the nullspace of \,$\square$\,. We are interested
in two \,$\R^{\times}$--invariant functions, to wit the restrictions
of \ $\Theta\,\Phi={\mathcal
F}_{\mathrm{symp}}\,(\widetilde{\Phi}\,dm)$\, to the hypersurface of
equation \ $\sc{x,\,\xi}=1$\,, or to the cone \,$\Sigma_n$\, of
equation \ $\sc{x,\,\xi}=0$\,: actually, we immediately need to
consider the results of these restrictions as living on the quotient
\,${\mathcal X}_n^{\bullet}$\, or \,$\Sigma_n^{\bullet}$\, of the
corresponding hypersurface by \,$\R^{\times}$\, and, for clarity ---
a notation already used in Proposition 2.3 --- we denote by a double
bar the operation of restriction followed by the one of going to the
quotient set. With this convention, we set
\begin{equation}\label{3AHPhi}
{\mathcal A}\,\Phi={\mathcal
F}_{\mathrm{symp}}\,(\widetilde{\Phi}\,dm)\,\biggr|\!\biggr|
_{{\mathcal X}_n^{\bullet}}\qquad\qquad {\mathcal H}\,\Phi={\mathcal
F}_{\mathrm{symp}}\,(\widetilde{\Phi}\,dm)\,\biggr|\!\biggr|
_{\Sigma_n^{\bullet}}\ .
\end{equation}
We shall first study the operator \,${\mathcal H}$\,, which will
turn out to be closely related to the involution \,${\mathcal K}$\,
which we have in mind; then, the operator \,${\mathcal A}$\,,
restricted to \,${\mathcal K}$--invariant functions, will provide an
isomorphism with a dense subspace of the continuous part of
\,$L^2({\mathcal X}_n^{\bullet})$\,. Under the transfer by this
isomorphism, the operator \,$\Delta_n$\, will appear as the square
of an
Euler--type differential operator on \,$\Sigma_n^{\bullet}$\,.\\

\begin{theorem}
Let \,${\mathcal H}$\, be the operator from functions on
\,$\Sigma_n^{\bullet}$\, (say, \,$C^{\infty}$\, with compact
support) to functions on \,$\Sigma_n^{\bullet}$\, characterized by
the identity
\begin{equation}\label{3intrH}
{\mathcal H}\,\Phi=({\mathcal
F}_{\mathrm{symp}}\,(\widetilde{\Phi}\,dm))\,\biggr|\!\biggr|_{{\Sigma_n}^{\bullet}}\,.
\end{equation}
It extends as an unbounded self--adjoint operator on
\,$L^2(\Sigma_n^{\bullet})$\,: moreover, denoting as \,$r$\, the
self--adjoint operator defined by the equation (in the coordinates
\,$(x_*,\,\xi_*)$\, on \,$\Sigma_n^{\bullet}$\,) \ $\langle
\xi_*\,,\,\frac{\partial}{\partial \xi_*}\rangle=-\frac{n}{2}-ir$\,,
one has
\begin{equation}\label{3HstarH}
{\mathcal H}^2={\mathcal H}^*{\mathcal H}=\pi\,
\frac{\Gamma(ir)\,\Gamma(-ir)}{\Gamma(\hf+ir)\,\Gamma(\hf-ir)}\,.
\end{equation}
Set, assuming that \,$\Phi=\Phi(y_*,\,\eta_*)$\, has, with respect
to \,$\eta_*$\,, the parity associated with \,$\delta$\,,
\begin{equation}\label{3compukPhi}
{\mathcal
K}\Phi=(-1)^{\delta}\,\pi^{-\hf-2ir}\,\frac{\Gamma(\hf+ir)}{\Gamma(-ir)}\,
{\mathcal H}\,\Phi\,.
\end{equation}
Then, the operator \,${\mathcal K}$\, on
\,$L^2(\Sigma_n^{\bullet})$\,  is a unitary involution and, for
every function \,$\Phi\in C^{\infty}(\Sigma_n^{\bullet})$\, with
compact support, one has the identity
\begin{equation}\label{3FsympK}
{\mathcal F}_{\mathrm{symp}}\,(\widetilde{{\mathcal
K}\Phi}\,dm)={\mathrm{Inv}}\, ({\mathcal
F}_{\mathrm{symp}}\,(\widetilde{\Phi}\,dm))\,.
\end{equation}\\
\end{theorem}

\begin{proof}
The function \ $\Theta\,\Phi={\mathcal
F}_{\mathrm{symp}}\,(\widetilde{\Phi}\,dm)$\, is given by the
integral
\begin{equation}\label{3ThetaPhi}
(\Theta\,\Phi)(x,\,\xi)=\int
e^{2i\pi\,(\sc{x,\,\eta}-\sc{y,\,\xi})}\,\Phi(\frac{y_*}{y_{n+1}},\,\,y_{n+1}\,\eta_*)\
\frac{dy\,d\eta_*}{|y_{n+1}|}\,,
\end{equation}
in which \ $\eta_{n+1}=-\frac{\sc{y_*,\,\eta_*}}{y_{n+1}}$\,. In
order to make the operator \,${\mathcal H}$\, introduced in
(\ref{3intrH}) explicit, it suffices to set \
$(x;\,\xi)=(x_*,1;\,\xi_*,-\sc{x_*,\,\xi_*})$\,, since this has been
our choice of coordinates on \,$\Sigma_n^{\bullet}$\,. After one
renames as \,$t$\, the integration variable \,$y_{n+1}$\, and one
changes \,$(y_*,\,\eta_*)$\, to
\,$(y_{n+1}\,y_*,\,\frac{\eta_*}{y_{n+1}})$\,, one obtains
\begin{equation}\label{3HPhib}
({\mathcal H}\,\Phi)(x_*,\,\xi_*)=\int
e^{2i\pi\,(\frac{\sc{x_*-y_*,\,\eta_*}}{t}+
t\,\sc{x_*-y_*,\,\xi_*})}\ \Phi(y_*,\,\eta_*)\
\frac{dt}{|t|}\,dy_*\,d\eta_*\,.
\end{equation}
Consequently, the operator \,${\mathcal H}$\, is formally
self--adjoint as an operator on \,$L^2(\Sigma_n^{\bullet})$\,.
Though it resembles a Fourier transformation, it involves a
restriction to some hypersurface and a dual operation, and it is not
unitary.
Its Fourier--transformed expression is easier to manage.\\

Denote as \,${{\mathcal F}}_*$\, the \,$(2n)$--dimensional version
of the symplectic Fourier transformation, and remark that the parity
of \,$\Phi$\, with respect to the set of Greek variables is the same
as that of \,${\mathcal F}_*\,\Phi$\, with respect to the roman
ones. Starting from (\ref{3HPhib}), one easily obtains
\begin{align}\label{3Foutransf}
({\mathcal F}_*\,{\mathcal H}\,\Phi)(y_*,\,\eta_*)&=\intR
|t|^{n-1}\,e^{2i\pi t\sc{y_*,\,\eta_*}}
\ ({\mathcal F}_*\,\Phi)(-t^2\,y_*,\,\eta_*)\,dt\nonumber\\
&=\intR |t|^{n-1}\,e^{2i\pi t}\,|\sc{y_*,\,\eta_*}|^{-n} \
({\mathcal
F}_*\,\Phi)\l(-t^2\,\frac{y_*}{(\sc{y_*,\,\eta_*})^2},\,\eta_*\r)\
dt\,.
\end{align}
Set \ ${\mathcal F}_*\,r\,{\mathcal F}_*^{-1}=\hat{r}$\,, {\em
i.e.,\/} \ $\sc{y_*,\,\frac{\partial}{\partial y_*}}
=-\frac{n}{2}+i\,\hat{r}$\,. Define
\begin{equation}\label{3defM}
({\mathcal F}_*\,{\mathcal M}\,\Phi)(y_*,\,\eta_*)=\intR
|t|^{n-1}\,e^{2i\pi t} \ ({\mathcal
F}_*\,\Phi)\l(t^2\,y_*,\,\eta_*\r)\ dt\,,
\end{equation}
so that
\begin{equation}\label{3HM}
{\mathcal F}_*\,{\mathcal H}\,\Phi={\mathcal J}\,({\mathcal
F}_*\,{\mathcal M}\,\Phi)
\end{equation}
if we denote as \,${\mathcal J}$\, the involution characterized by
the equation
\begin{equation}\label{3defJ}
({\mathcal J}\,\Xi)(y_*,\,\eta_*)=|\sc{y_*,\,\eta_*}|^{-n} \
\Xi\l(-\frac{y_*}{(\sc{y_*,\,\eta_*})^2},\,\eta_*\r)\,:
\end{equation}
note that \,${\mathcal J}$\, anticommutes with \,$\hat{r}$\,. One
may rewrite (\ref{3defM}) as
\begin{equation}\label{3rewrM}
{\mathcal F}_*\,{\mathcal
M}\,\Phi=\pi^{\hf-2i\,\hat{r}}\,\frac{\Gamma(i\,\hat{r})}{\Gamma(\hf-i\,\hat{r})}
\ {\mathcal F}_*\,\Phi\,,
\end{equation}
a result which, when combined with (\ref{3HM}), leads to
\begin{equation}
{\mathcal F}_*\,{\mathcal H}^*\,{\mathcal H}\,{{\mathcal F}_*}^{-1}=
\pi\,\frac{\Gamma(i\,\hat{r})\,\Gamma(-i\,\hat{r})}{\Gamma(\hf+i\,\hat{r})\,
\Gamma(\hf-i\,\hat{r})}\,,
\end{equation}
then to (\ref{3HstarH}).\\

\noindent The definition (\ref{3compukPhi}) of \,${\mathcal
K}\,\Phi$\,, together with (\ref{3HM}), leads to
\begin{align}\label{3FouKPhi}
{\mathcal F}_*\,{\mathcal
K}\,\Phi&=(-1)^{\delta}\,\pi^{-\hf-2i\,\hat{r}}\,
\frac{\Gamma(\hf+i\,\hat{r})}{\Gamma(-i\,\hat{r})}\,{\mathcal
F}_*\,{\mathcal H}\,\Phi
\nonumber\\
&=(-1)^{\delta}\,{\mathcal J}\,\l(\,\pi^{-\hf+2i\,\hat{r}}\,
\frac{\Gamma(\hf-i\,\hat{r})}{\Gamma(i\,\hat{r})}\,{\mathcal
F}_*\,{\mathcal M}\,\Phi\,\r)
\nonumber\\
&=(-1)^{\delta}\,{\mathcal J}\,({\mathcal F}_*\,\Phi)
\end{align}
in view of (\ref{3rewrM}). In terms of the ($(2n)$--dimensional)
symplectic Fourier transform \,${\mathcal F}_*\,\Phi$\, of
\,$\Phi$\,, one may rewrite (\ref{3ThetaPhi}) as
\begin{multline}\label{ThetaPhiFou}
(\Theta\,\Phi)(x,\,\xi)=|x_{n+1}|^{-n}\,\int_{\R^n\times \R^n\times
\R}
|t|^{n-1}\ ({\mathcal F}_*\,\Phi)(q,\,p)\\
\exp\,\l(\frac{2i\pi}{x_{n+1}}\,\l(t^2\,\sc{q,\,\xi_*}-t\,(\sc{q,\,p}+\sc{x,\,\xi})
+\sc{x_*,\,p}\,\r)\r)\ dq\,dp\,dt\,,
\end{multline}
as seen after a perfectly elementary, if somewhat lengthy to write down, computation. \\

Next, the parity of \,${\mathcal F}_*\,\Phi$\, with respect to
\,$q$\, is the same as that of \,$\Phi$\, with respect to \,$y_*$\,,
{\em i.e.,\/} \,$\delta$\,. Substituting the result of
(\ref{3FouKPhi}) into (\ref{ThetaPhiFou}), one sees that, in order
to obtain \,$(\Theta\,{\mathcal K}\,\Phi)(x,\,\xi)$\,, it suffices
to take the right--hand side of (\ref{ThetaPhiFou}), replacing
\,$({\mathcal F}_*\,\Phi)(q,\,p)$\, by \ $|\sc{q,\,p}|^{-n} \
({\mathcal F}_*\,\Phi)(\frac{q}{(\sc{q,p})^2},\,p)$\,: we shall not
display the result, but shall immediately perform the change of
variables \ $(q,\,p)\mapsto (\frac{q}{(\sc{q,p})^2},\,p)$\,, which
leads to
\begin{multline}\label{ThetaPhiK}
(\Theta\,{\mathcal
K}\,\Phi)(x,\,\xi)=|x_{n+1}|^{-n}\,\int_{\R^n\times \R^n\times \R}
|t|^{n-1}\ |\sc{q,\,p}|^{-n}\,({\mathcal F}_*\,\Phi)(q,\,p)\\
\exp\,\l(\frac{2i\pi}{x_{n+1}}\,\l(\,\frac{t^2}{(\sc{q,\,p})^2}\,\sc{q,\,\xi_*}-
t\,(\frac{1}{\sc{q,\,p}}+\sc{x,\,\xi})+\sc{x_*,\,p}\,\r)\r)\
dq\,dp\,dt\,.
\end{multline}\\

On the other hand, changing in (\ref{ThetaPhiFou}) the integration
variable \,$t$\, for \,$s$\, one obtains
\begin{multline}\label{ThetaPhiFouInv}
({\mathrm{Inv}}\,(\Theta\,\Phi))(x,\,\xi)=|x_{n+1}|^{-n}\,\int_{\R^n\times
\R^n\times \R}
|s|^{n-1}\ ({\mathcal F}_*\,\Phi)(q,\,p)\\
\exp\,\l(\frac{2i\pi}{x_{n+1}}\,\l(s^2\,\sc{q,\,\xi_*}-s\,(\sc{q,\,p}
\sc{x,\,\xi}+1)+\sc{x_*,\,p}\,\r)\r)\ dq\,dp\,ds\,.
\end{multline}
Setting \ $s=\frac{t}{\sc{q,\,p}}$\, in the last integral, and
comparing the result to (\ref{ThetaPhiK}), one obtains the equation
\ ${\mathrm{Inv}}\,(\Theta\,\Phi)=
\Theta\,{\mathcal K}\,\Phi$\,, just the same as (\ref{3FsympK}).\\
\end{proof}

We can now state and prove the main result of this section.

\begin{theorem}
Recall from {\em (\ref{3AHPhi})\/} that
\begin{equation}
{\mathcal
H}\,\Phi=\Theta\,\Phi\,\biggr|\!\biggr|_{{\Sigma_n}^{\!\!\bullet}}\,,
\qquad\qquad {\mathcal
A}\,\Phi=\Theta\,\Phi\,\biggr|\!\biggr|_{{\mathcal
X}_n^{\bullet}}\,,
\end{equation}
and that the map \,${\mathcal H}$\, has been analyzed in Theorem
{\em 3.1\/}. Let \,$D_{\mathcal K}({\mathcal
H})\subset\,L^2(\Sigma_n^{\bullet})$\, be the space of \,${\mathcal
K}$--invariant functions in the domain of the self--adjoint operator
\,${\mathcal H}$\,. The linear map \,${\mathcal A}$\, is a linear
isomorphism from \,$D_{\mathcal K}({\mathcal H})$\, onto the
subspace of \,$L^2({\mathcal X}_n^{\bullet})=L^2(G_n/H_n)$\,
corresponding to the continuous part of the decomposition of this
latter space under the quasiregular action of \,$G_n$\,. One has the
identity
\begin{equation}
\Vert\,{\mathcal A}\,\Phi\,\Vert^2=2\ \Vert\,{\mathcal
H}\,\Phi\,\Vert^2
\end{equation}
for every \,$\Phi \in D_{\mathcal K}({\mathcal H})$\,. Finally, the
operator \ $(\frac{n}{2}+\sc{\eta_*,\,\frac{\partial}{\partial
\eta_*}})^2$\, commutes with \,${\mathcal K}$\, and, for every
function \,$\Phi\in D_{\mathcal K}({\mathcal H})$\, such that \
$(\frac{n}{2}+\sc{\eta_*,\,\frac{\partial}{\partial
\eta_*}})^2\,\Phi$\, lies in \,$D_{\mathcal K}({\mathcal H})$\, as
well, the identity
\begin{equation}\label{3diffsqu}
{\mathcal A}\,\l(\frac{n}{2}+\sc{\eta_*,\,\frac{\partial}{\partial
\eta_*}}\r)^2\,\Phi= (\Delta_n-\frac{n^2}{4})\,{\mathcal A}\,\Phi
\end{equation}
holds. It thus reduces the study of the continuous part of the
spectral decomposition of \,$\Delta_n$\, to the (trivial) spectral
theory of an Euler--type operator on the cone
\,$\Sigma_n^{\bullet}$\,.\\
\end{theorem}

\begin{proof}
On \,${\mathcal X}_n^{\bullet}$\, we still use the coordinates
\,$(x_*,\,\xi_*)$\,, corresponding this time to the point \,
$(x,\,\xi)=(x_*,1;\,\xi_*, 1-\sc{x_*,\xi_*})$\,: these coordinates
are related by the equation (\ref{2ssigma}) to the coordinates
\,$(s,\,\sigma)$\,, more useful when dealing with the symbolic
calculus; here, there are advantageous, since in particular the
\,$G_n$--invariant measure \,$|1+\langle
s,\sigma\rangle|^{-n-1}\,ds\,d\sigma$\, on \,${\mathcal
X}_n^{\bullet}$\, reduces again to \,$dx_*\,d\xi_*$\,. With the same
computation as the one in the beginning of the proof of Theorem 3.1,
we find the following equation, to be compared to (\ref{3HPhib}):
\begin{equation}
({\mathcal A}\,\Phi)(x_*,\,\xi_*)=\int
e^{2i\pi\,(\frac{\sc{x_*-y_*,\,\eta_*}}{t}+
t\,\sc{x_*-y_*,\,\xi_*})}\,e^{-2i\pi t}\ \Phi(y_*,\,\eta_*)\
\frac{dt}{|t|}\,dy_*\,d\eta_*\,,
\end{equation}
the Fourier--transformed version of which, to be compared to
(\ref{3Foutransf}), is
\begin{equation}\label{3FouAstar}
({\mathcal F}_*\,{\mathcal A}\,\Phi)(y_*,\,\eta_*)=\intR
|t|^{n-1}\,e^{2i\pi t\sc{y_*,\,\eta_*}} \,e^{-2i\pi\,t^{-1}}
\ ({\mathcal F}_*\,\Phi)(-t^2\,y_*,\,\eta_*)\,dt\,.\\
\end{equation}

The \,$n$--dimensional analogue of \,$\square$\,, to wit the
operator on \,$\R^n\times \R^n$\, defined by the equation \
$\square_*= \sum_{j=1}^n\,\frac{\partial^2}{\partial x_j\,\partial
\xi_j}$\,, can be interpreted in two different ways, since the
coordinates \,$(x_*,\,\xi_*)$\, can be used on \,${\mathcal
X}_n^{\bullet}$\, as well as on \,$\Sigma_n^{\bullet}$\,: of course,
the two operators obtained are conceptually different, and we shall
denote them as \,$\square_*^{\mathcal X}$\, and
\,$\square_*^{\Sigma}$\, respectively. Incidentally, one may prove
the equation
\begin{equation}
\Delta_n=\square_*^{\mathcal
X}-\,\sc{\xi_*,\,\frac{\partial}{\partial \xi_*}}
\,(n+\sc{\xi_*,\,\frac{\partial}{\partial \xi_*}})\,.
\end{equation}
Each of the two operators just defined is a self-adjoint operator
with continuous spectrum, and we denote as \,$P_+^{\mathcal X}$\,
and \,$P_-^{\mathcal X}$\, ({\em resp.\/} \,$P_+^{\Sigma}$\, and
\,$P_-^{\Sigma}$\,) the orthogonal projections onto the positive and
negative spaces of \,$\square_*^{\mathcal X}$\, ({\em resp.\/} \,$\square_*^{\Sigma}$\,).\\

We shall compute the norms of \,$P_{\pm}^{\mathcal X}\,{\mathcal
A}\,\Phi$\, in terms of \,$\Phi$\, separately: since \ ${\mathcal
F}_*\,\square_*\,{{\mathcal F}_*}^{-1}=4\pi^2\,\sc{y_*,\,\eta_*}$\,,
one has for instance
\begin{equation}
\Vert\,P_+^{\mathcal X}\,{\mathcal A}\,\Phi\,\Vert^2= \int
_{\sc{y_*,\,\eta_*}\,>0}\,|({\mathcal F}_* \,{\mathcal
A}\,\Phi)(y_*,\,\eta_*)|^2\ dy_*\,d\eta_*\,.
\end{equation}
With \ $y_*=(y_{**},\,y_n)$\, it is convenient to substitute for
\,$(y_*,\,\eta_*)$\, the (singular) coordinates
\,$(q,\,\eta_*;\,w)=(\frac{y_{**}}{y_n},\,\eta_*;\,\sc{y_*,\,\eta_*})\in
P_{n-1}(\R)\times \R^n\times\R$\,: then,
\begin{equation}
dy_*\,d\eta_*=\frac{|y_n|^n}{|\sc{y_*,\,\eta_*}|}\,dq\,d\eta_*\,dw=
\frac{|w|^{n-1}}{|\sc{q,\,\eta_{**}}+\eta_n|^n}\,dq\,d\eta_*\,dw\,,
\end{equation}
and \ $i\hat{r}=\frac{n}{2}+w\,\frac{\partial}{\partial w}$\,.\\

In terms of the coordinates just introduced, the operator \
${\mathcal B}={\mathcal F}_*\, {\mathcal A}\,{\mathcal F}_*^{-1}$\,,
as given by (\ref{3FouAstar}), really reduces to the operator on
functions of one variable only defined by
\begin{equation}\label{3h}
({\mathcal B}\,\phi)(w)=\intR |t|^{n-1}\,e^{2i\pi tw}\,e^{-2i\pi
t^{-1}}\,\phi(-t^2 w)\,dt\,:
\end{equation}
the coordinates \,$(q,\,\eta_*)$\, are now simple parameters and
have been omitted for
clarity; we are interested in \,${\mathcal B}$\, as an unbounded operator from the space \\
\,$L^2((-\infty,\,0);\,|w|^{n-1}\,dw)$\, to the space
\,$L^2((0,\,\infty);\,w^{n-1}\,dw)$\,. Using (\ref{3FouKPhi}) and
(\ref{3defJ}) again, one sees that the symmetry condition \
$\Phi={\mathcal K}\,\Phi$\, expresses itself in a way independent of
\,$\delta$\,: in the new coordinates, and still forgetting the ones
that are present as simple parameters, one has
\begin{equation}\label{3sym}
\phi(w)=|w|^{-n}\,\phi(\frac{1}{w})\,.
\end{equation}
Set, for real \,$x,\,y$\,,
\begin{equation}
\chi(y)=e^{ny}\,\phi(-e^{2y})\,,\qquad\qquad ({\mathcal
D}\,\chi)(x)=e^{-nx}\,({\mathcal B}\,\phi)(e^{-2x})\,:
\end{equation}
then, the function \,$\chi$\, is even. Also, one may then rewrite
the last equation as
\begin{align}\label{3mathcalD}
({\mathcal D}\,\chi)(x)&=2\,e^{-nx}\,\int_0^{\infty} t^{n-1}\
\cos\,(2\pi\,(t\,e^{-2x}-t^{-1}))
\ \phi(-t^2\,e^{-2x})\,dt\nonumber\\
&=2\,\intR \cos\,(4\pi\, e^{-x}\,\sinh y)\ \chi(y)\,dy\,.
\end{align}
Now, the inverse Fourier transform of the function \ $y\mapsto
e^{4i\pi e^{-x}\sinh y}$\,, evaluated at \,$s$\,, can be found from
\cite[p.\,86]{mos}, after one has inserted a factor \,$e^{-4\pi
\varepsilon\cosh y}$\, for convergence and let \,$\varepsilon$\, go
to zero: the result is the function \ $s\mapsto 2\,e^{-\pi^2
s}\,K_{2i\pi s}(4\pi e^{-x})$\,: denoting as \,$\hat{\chi}$\, the
usual Fourier transform of \,$\chi$\,, one finds
\begin{equation}
({\mathcal D}\,\chi)(x)=4\intR (\cosh\,\pi^2 s)\, K_{2i\pi s}(4\pi
e^{-x})\,\hat{\chi}(s)\,ds
\end{equation}
and
\begin{multline}\label{3bigfirst}
\Vert\,{\mathcal
D}\chi\,\Vert^2_{L^2(\R;\,dx)}=16\,\int_0^{\infty}\frac{dt}{t}
\int_{\R^2} K_{2i\pi s_1}(4\pi t)\,K_{2i\pi s_2}(4\pi t)\\
(\cosh\,\pi^2 s_1)\,(\cosh\,\pi^2 s_2)\
\hat{\chi}(s_1)\,\overline{\hat{\chi}(s_2)}\ ds_1\,ds_2\,.
\end{multline}
This integral can be computed with the help of \cite[p.\,46]{aumod},
leading to
\begin{align}\label{3dchisqu}
\Vert\,{\mathcal D}\,\chi\,\Vert^2_{L^2(\R;\,dx)}&= 8\intR
\Gamma(2i\pi s)\,\Gamma(-2i\pi s)\,(\cosh\,\pi^2 s)^2\
|\hat{\chi}(s)|^2\,ds
\nonumber\\
&=2\pi\,\intR \frac{\Gamma(i\pi s)\,\Gamma(-i\pi s)}
{\Gamma(\hf+i\pi s)\,\Gamma(\hf-i\pi s)}\ |\hat{\chi}(s)|^2\,ds\,.
\end{align}
Since the multiplication by \,$s$\, corresponds, under the Fourier
transformation, to \ $\frac{1}{2i\pi}\,\frac{d}{dy}$\,, and since
\begin{equation}
\frac{1}{2i\pi}\,\frac{d}{dy}\,\l(e^{ny}\,\phi(-e^{2y})\r)=\frac{1}{i\pi}\,
e^{ny}\,\l((\frac{n}{2}+w\,\frac{d}{dw})\,\phi\r)\,(w=-e^{2y})\,,
\end{equation}
finally using the equation \
$i\hat{r}=\frac{n}{2}+w\,\frac{\partial}{\partial w}$\, mentioned
above, one finds
\begin{equation}\label{3isom}
\Vert\,P_+^{\mathcal X}\,{\mathcal A}\,\Phi\,\Vert^2=2\pi\
\l(P_-^{\Sigma}\,\Phi\,\biggr|\,
\frac{\Gamma(ir)\,\Gamma(-ir)}{\Gamma(\hf+ir)\,\Gamma(\hf-ir)}\,P_-^{\Sigma}\,\Phi\r)\,.
\end{equation}\\

The study of \,$P_-^{\mathcal X}\,{\mathcal A}\,\Phi$\, calls this
time for that of \,${\mathcal B}$\, as an unbounded operator from
the space \,$L^2((0,\,\infty);\,w^{n-1}\,dw)$\, to the space
\,$L^2((-\infty,\,0);\,|w|^{n-1}\,dw)$\,. We then set
\begin{equation}
\psi(y)=e^{ny}\,\phi(e^{2y})\,,\qquad\qquad ({\mathcal
C}\,\psi)(x)=e^{-nx}\,({\mathcal B}\,\phi)(-e^{-2x})\,:
\end{equation}
the analogue of (\ref{3mathcalD}) is the equation
\begin{equation}
({\mathcal C}\,\psi)(x)=2\,\intR \cos\,(4\pi\, e^{-x}\,\cosh y)\
\psi(y)\,dy\,.
\end{equation}
Using \cite[p.\,86,\,67,\,66]{mos}, one can see that the inverse
Fourier transform of the function \ $y\mapsto e^{4i\pi e^{-x}\cosh
y}$\,, evaluated at \,$s$\,, is
\begin{align}
2\ {\mathrm{lim}}_{\varepsilon \to 0}\,K_{2i\pi
s}(4\pi\,(\varepsilon-i\,e^{-x}))&=
i\pi\,e^{-\pi^2s}\,H^{(1)}_{2i\pi s}(4\pi\,e^{-x})\nonumber\\
&=\pi\,\frac{e^{-\pi^2s}\,J_{-2i\pi
s}(4\pi\,e^{-x})-e^{\pi^2s}\,J_{2i\pi s}(4\pi\,e^{-x})}
{i\,\sinh\,2\pi^2 s}\,,
\end{align}
and that of the function \ $y\mapsto \cos\,(4i\pi e^{-x}\cosh y)$\,
is the real function
\begin{equation}
F_s(e^{-x})=\pi\,\frac{J_{-2i\pi s}(4\pi\,e^{-x})-J_{2i\pi
s}(4\pi\,e^{-x})} {2i\,\sinh\,\pi^2s}\,:
\end{equation}
one has
\begin{equation}
({\mathcal C}\,\psi)(x)=2\,\intR F_s(e^{-x})\,\hat{\psi}(s)\,ds\,,
\end{equation}
and we need to compute
\begin{align}\label{3needcomp}
\Vert\,{\mathcal
C}\psi\,\Vert^2_{L^2(\R)}&=4\,\int_0^{\infty}\frac{dt}{t}
\int_{\R^2} F_{s_1}(t)\,F_{s_2}(t)
\ \hat{\psi}(s_1)\,\overline{\hat{\psi}(s_2)}\ ds_1\,ds_2\nonumber\\
&={\mathrm{lim}}_{\rho\to 0} \,\int_{\R^2} K_{\rho}(s_1,\,s_2) \
\hat{\psi}(s_1)\,\overline{\hat{\psi}(s_2)}\ ds_1\,ds_2\,,
\end{align}
with
\begin{multline}
K_{\rho}(s_1,\,s_2)\\
=-\frac{\pi^2}{\sinh\,\pi^2s_1\,\sinh\,\pi^2s_2}
\sum_{\epsilon_1^2=\epsilon_2^2=1}
\epsilon_1\,\epsilon_2\,\int_0^{\infty} \, J_{2i\pi\,\epsilon_1
s_1}(4\pi\,t)\,J_{2i\pi\,\epsilon_2 s_2}(4\pi\,t) \
t^{\rho-1}\,dt\,.
\end{multline}
We set
\begin{equation}
\sc{\epsilon,\,s}=\epsilon_1\,s_1+\epsilon_2\,s_2\,,\qquad\qquad
\sc{\check{\epsilon},\,s}=-\epsilon_1\,s_1+\epsilon_2\,s_2
\end{equation}
and we use \cite[p.\,99]{mos}
\begin{multline}\label{3JJ}
\int_0^{\infty} J_{2i\pi \epsilon_1 s_1}(4\pi t)\,J_{2i\pi
\epsilon_2 s_2}(4\pi t)\ t^{\rho-1}\,dt
=\hf\,(2\pi)^{-\rho}\,\Gamma(1-\rho)\\
\times\ \frac{\Gamma(\frac{\rho}{2}+i\pi\,\sc{\epsilon,\,s})}
{\Gamma(1-\frac{\rho}{2}-i\pi\,
\sc{\check{\epsilon},\,s})\,\Gamma(1-\frac{\rho}{2}+i\pi\,\sc{\epsilon,\,s})\,
\Gamma(1-\frac{\rho}{2}+i\pi\,\sc{\check{\epsilon},\,s})}\,.
\end{multline}\\

\noindent Note that, for \,$s_1\neq \pm s_2$\,,
\,$K_{\rho}(s_1,\,s_2)$\, goes to zero as \,$\rho\to 0$\,, since
\begin{equation}
K_0(s_1,\,s_2)=-\frac{1}{4i\pi\,\sinh\,\pi^2s_1\,\sinh\,\pi^2s_2}
\,\sum_{\epsilon_1^2=\epsilon_2^2=1} \epsilon_1\,\epsilon_2\,
\frac{\sinh\,\pi^2\,\sc{\check{\epsilon},\,s}}
{\sc{\epsilon,\,s}\,\sc{\check{\epsilon},\,s}}
\end{equation}
changes to its negative under the change of parameters \
$(\epsilon_1,\,\epsilon_2) \mapsto (-\epsilon_1,\,-\epsilon_2)$\,.
Consequently, as a distribution on \,$\R^2$\,, \
$K_{\rho}(s_1,\,s_2)$\, has, as \,$\rho\to 0$\,, a limit supported
in the union of the two lines \ $s_1\pm s_2=0$\,. The calculation
can be completed in the same way as that in \cite[p.\,46]{aumod}.
The only contributions to \ $K_{\rho}(s_1,\,s_2)$\, which do not
vanish in the limit as \,$\rho \to 0$\, come from the pole at zero
of the Gamma factor on the top of the right--hand side of
(\ref{3JJ}), and we may replace this Gamma factor by its rational
equivalent. The terms with \ $\epsilon_1=\epsilon_2$\, add up to an
expression which has the same limit, as \,$\rho\to 0$\,, as
\begin{multline}
-\frac{\pi^2}{2\,\sinh\,\pi^2s_1\,\sinh\,\pi^2s_2}\,
[\,\Gamma(1-\frac{\rho}{2}-i\pi\,\sc{\check{\epsilon},\,s})
\Gamma(1-\frac{\rho}{2}+i\pi\,\sc{\check{\epsilon},\,s})\,]^{-1}\\
\times\,\l[\,\frac{1}{\frac{\rho}{2}+i\pi\,\sc{\epsilon,\,s}}
+\frac{1}{\frac{\rho}{2}-i\pi\,\sc{\epsilon,\,s}}\,\r]\,:
\end{multline}
the function on the second line goes to the distribution \
$2\,\delta(s_1+s_2)$\,, so that (remembering that \,$\hat{\psi}$\,
is an even function) the contribution to (\ref{3needcomp}) of the
terms with \,$\epsilon_1=\epsilon_2$\, is the integral of \
$|\hat{\psi}(s)|^2\,ds$\, against the coefficient
\begin{equation}
\frac{\pi^2}{(\sinh\,\pi^2s)^2}\,[\,\Gamma(1-2i\pi
s)\,\Gamma(1+2i\pi s)\,]^{-1}= \pi\,\frac{\Gamma(i\pi
s)\,\Gamma(-i\pi s)}{\Gamma(\hf+i\pi s)\,\Gamma(\hf-i\pi s)}\,;
\end{equation}
one obtains the same result from the consideration of the terms with
\ $\epsilon_1=-\epsilon_2$\,. Comparing this to (\ref{3dchisqu}),
the net result is that
\begin{equation}
\Vert\,P_-^{\mathcal X}\,{\mathcal A}\,\Phi\,\Vert^2=2\pi\
\l(P_+^{\Sigma}\,\Phi\,\biggr|\,
\frac{\Gamma(ir)\,\Gamma(-ir)}{\Gamma(\hf+ir)\,\Gamma(\hf-ir)}\,P_+^{\Sigma}\,\Phi\r)\,.
\end{equation}

The proof of Theorem 3.2 is now complete: for the first part, we
only need comparing (\ref{3isom}) to (\ref{3HstarH}); next, the
equation (\ref{3diffsqu}) is a consequence of
(\ref{2eigen}).\\
\end{proof}

Needless to say, the inversion problem, {i.e.,\/} the problem of
recovering \,$\Phi$\, in terms of \,${\mathcal A}\,\Phi$\,, under
the assumption that \,$\Phi$\, is \,${\mathcal K}$--invariant, is
easy. Since the computations can be made by following the same
transformations as above, let us satisfy ourselves with the result
of the computation:
\begin{equation}
({\mathcal
F}_*\,\Phi)(y_*,\,\eta_*)=\frac{|1-\sc{y_*,\,\eta_*}|}{2}\, \intR
|t|^n\,e^{2i\pi t\,(1+\sc{y_*,\,\eta_*})} \ ({\mathcal
F}_*\,{\mathcal A}\,\Phi)(-t^2\,y_*,\,\eta_*)\,dt\,.
\end{equation}\\
In view of Theorem 3.2, the preceding results, coupled with a Mellin
transformation on \,$\Sigma_n^{\bullet}$\,, provide a
diagonalization of the continuous part of the spectral
decomposition of \,$\Delta_n$\,.\\


\section{Composition: the soft approach}

The composition problem is a difficult one: we shall have a glimpse
of it in the last part of the present paper, in which we shall
consider it for some special symbols. In this short section, we
explain why the obvious approach to the composition problem does not
lead anywhere, in contradiction to what is the case with the
symbolic calculus on \,$\R^{n+1}$\,. It would be misleading to
believe that, in the quantization of symmetric spaces, the
composition \,$f_1\,\#\,f_2$\, of two symbols can be, even in a
rough way, described by means of a series of differential
expressions in the pair
of symbols under consideration; also, the integral formula is essentially worthless.\\

Definition 2.1 gives the integral kernels of the operators
\,${\mathrm{Op}}_{i\lambda,\varepsilon}(f)\,\theta_{i\lambda,\varepsilon}^{-1}$\,
and
\,$\theta_{i\lambda,\varepsilon}\,{\mathrm{Op}}^{\vee}_{i\lambda,\varepsilon}(h)$\,;
on the other hand, (\ref{2thetaugly}) gives the integral kernel of
the intertwining operator \,$\theta_{i\lambda,\varepsilon}$\,, or of
its inverse since it is unitary. It is thus immediate to obtain the
following integral formula,  analogous to the formula
\begin{equation}\label{4J}
(Jh)(x,\xi)=\int h(y,\eta)\,e^{2i\pi\,\langle
x-y,\eta-\xi\rangle}\,dy\,d\eta
\end{equation}
which, in the calculus on \,$\R^{n+1}$\,, makes it possible to link
the antistandard symbol \,$h$\, of some operator to its standard
symbol \,$Jh$\,.
\begin{proposition}
If \,$h\in L^2(G/H)$\,, one has
\begin{equation}
{\mathrm{Op}}^{\vee}_{i\lambda,\varepsilon}(h)=
{\mathrm{Op}}_{i\lambda,\varepsilon}(f)
\end{equation}
with (setting \ $d\mu(t,\tau)=|1+\langle t,\tau
\rangle|^{-n-1}\,dt\,d\tau$\,)
\begin{equation}
f(s,\sigma)=|C_{i\lambda,\varepsilon}|^2\,\int \biggr|\,
\frac{(1+\langle s,\sigma \rangle)\,(1+\langle t,\tau \rangle)}
{(1+\langle s,\tau \rangle)\,(1+\langle t,\sigma
\rangle)}\,\biggr|_{\varepsilon} ^{\frac{n+1}{2}+i\lambda}\
h(t,\tau)\ d\mu(t,\tau)\,.
\end{equation}\\
\end{proposition}
This formula should really be understood as the fact that the
function \,$(\sigma,s)\mapsto |1+\langle
s,\sigma\rangle|_{\varepsilon}^{-\frac{n+1}{2}-i\lambda}\,
f(s,\sigma)$ \ is
\,$C_{i\lambda,\varepsilon}\,\bar{C}_{i\lambda,\varepsilon}^{-1}$\,
times the image of the function \,$(t,\tau)\mapsto |1+\langle t,\tau
\rangle|_{\varepsilon}^{-\frac{n+1}{2}+i\lambda}\,f(t,\tau)$\, under
the operator \
$\theta_{-i\lambda,\varepsilon}\,\otimes\,\theta_{-i\lambda,\varepsilon}$\,.
One may interpret the following integral formula in a similar way:
\begin{proposition}
Let \,$f_1$\, and  \,$f_2$\, lie in \,$ L^2(G/H)$\,. One has
\begin{equation}
{\mathrm{Op}}_{i\lambda,\varepsilon}(f_1)\,{\mathrm{Op}}_{i\lambda,\varepsilon}(f_2)=
{\mathrm{Op}}_{i\lambda,\varepsilon}(f_1\,\#\,f_2)
\end{equation}
with
\begin{equation}
(f_1\,\#\,f_2)(s,\sigma)= |C_{i\lambda,\varepsilon}|^2\,\int
\biggr|\,\frac{(1+\langle s,\sigma \rangle)\,(1+\langle t,\tau
\rangle)} {(1+\langle s,\tau \rangle)\,(1+\langle t,\sigma
\rangle)}\,\biggr|_{\varepsilon} ^{\frac{n+1}{2}+i\lambda}\
f_1(s,\tau)\,f_2(t,\sigma)\ d\mu(t,\tau)\,.
\end{equation}\\
\end{proposition}
For the sake of comparison, we may here recall that the integral
formula giving the symbol \,$f_1\underset{\R^{n+1}}{\#} f_2$\, of
the two operators with standard symbols \,$f_1$\, and \,$f_2$\, is
\begin{equation}\label{4intcompeuc}
(f_1\underset{\R^{n+1}}{\#} f_2)(x,\xi)= \int
f_1(x,\eta)\,f_2(y,\xi)\,e^{2i\pi\,\langle
x-y,\eta-\xi\rangle}\,d\eta\,dy\,.
\end{equation}\\

We have already come across the operator \,$J$\, which occurs in
(\ref{4J}): indeed, it has been pointed out, right after (\ref{2J}),
that one has the equation \ $J=\exp\, \frac{\square}{2i\pi}$\,:
writing the exponential as a series, this immediately leads to the
expansions
\begin{equation}\label{4tr}
Jh\,\sim\,\sum_{k=0}^{\infty}
\frac{1}{k\,!}\,\l(\frac{\square}{2i\pi}\r)^k\,h
\end{equation}
and
\begin{equation}\label{4comp}
(f_1\underset{\R^{n+1}}{\#} f_2)(x,\xi)\,\sim\, \sum_{k=0}^{\infty}
\frac{1}{k\,!}\,\l(\frac{\square_{y,\eta}}{2i\pi}\r)^k\,
\l(\,f_1(x,\,\xi+\eta)\,f_2(x+y,\,\xi)\,\r)(y=\eta=0)\,.
\end{equation}
These formulas lie at the foundations of classical
pseudodifferential analysis: they do not define convergent series
except for special symbols --- such as those of differential
operators
--- but they are valid, as useful asymptotic expansions, for symbols, or pairs of symbols, lying
in large appropriate classes. Let us emphasize that what makes
pseudodifferential analysis such a useful tool in partial
differential equations is the easy way it makes it possible to
define, and handle, {\em auxiliary\/} operators: you certainly do
not need it to compose differential operators. Harmonic analysts may
find an extra, immediate, satisfaction in the fact that if one
introduces from the start, in the usual way, a Planck's constant in
the definition of the (say, standard) symbolic calculus, this
constant will appear, in (\ref{4tr}) or (\ref{4comp}), as a
coefficient in front of \,$\square$\,: thus, these two asymptotic
expansions may also be viewed as series expansions with respect to
Planck's
constant.\\

As will be shown on the example which is the subject of this paper,
this feature of pseudodifferential analysis in Euclidean space does
not survive in the quantization of symmetric spaces, whether you
wish to interpret the would--be analogues of (\ref{4tr}) and
(\ref{4comp}) as asymptotic expansions or as series in the parameter
\,$\lambda^{-1}$\,,
sometimes viewed as some kind of analogue of  Planck's constant.\\

We view the developments in the present section as necessary for a
good understanding of the nature of quantization: however, they are,
to a certain extent, of a negative nature, and we shall be as brief
as possible.
 Using the coordinates \,$(s,\,\sigma)$\, on the phase space
\,${\mathcal X}_n^{\bullet}$\, ({\em cf.\/} (\ref{2ssigma})), we
note that the point \,$(s,\,\sigma)$\, is the image of the point
\,$(0,0)$\, under the matrix \ $g_{s,\sigma}=\l(\begin{smallmatrix}
I & \frac{s}{1+\langle s,\sigma\rangle} \\ -\sigma^{\dag} &
\frac{1}{1+\langle s,\sigma\rangle}\end{smallmatrix}\r)$\,.
Concentrating on the composition formula --- since the operator
\,$h\mapsto f$\, from Proposition 4.1 has already been studied ---
we may use covariance to reduce the pointwise study of
\,$(f_1\,\#\,f_2)(s,\sigma)$\, to that of
\,$(f_1\,\#\,f_2)(0,0)$\,. In perfect analogy with the way
(\ref{4comp}) was derived from (\ref{4intcompeuc}), it suffices to
compute the Fourier transform of a certain distribution to obtain
the following:
\begin{proposition}
Denote as \,$F_{i\lambda,\varepsilon}$\, the function of one real
variable defined as follow:
\begin{equation}
F_{i\lambda,\,\varepsilon}(\beta)=\frac{2\,(2\pi)^{\frac{1-n}{2}-i\lambda}}{\Gamma(\frac{1-n}{2}-i\lambda)}
\ \beta^{\frac{1-n}{4}-\frac{i\lambda}{2}}\,
K_{\frac{n-1}{2}+i\lambda}(4\pi\,\beta^{\hf})
\end{equation}
if \,$\beta>0$\,, and
\begin{multline}
F_{i\lambda,\,\varepsilon}(\beta)\\
=(2\pi)^{\frac{1-n}{2}-i\lambda}\,\Gamma(\frac{n+1}{2}+i\lambda) \
|\beta|^{\frac{1-n}{4}-\frac{i\lambda}{2}}\ \l[
J_{\frac{n-1}{2}+i\lambda}(4\pi\,|\beta|^{\hf})
-(-1)^{\varepsilon}\,J_{\frac{1-n}{2}-i\lambda}(4\pi\,|\beta|^{\hf})\,\r]
\end{multline}
if \,$\beta<0$\,. Then one has
\begin{equation}\label{4f1sharpf2}
(f_1\,\#\,f_2)(0,0)=
\l[\,F_{\lambda}\,\l(-\frac{1}{4\pi^2}\,\sum\frac{\partial^2}{\partial
t_j\,\partial \tau_j}\r)\,
(f_1(0,\tau)\,f_2(t,0))\,\r]\,((t,\,\tau)=(0,\,0)).
\end{equation}
\end{proposition}

\begin{proof}
The right--hand side of the equation
\begin{equation}
(f_1\,\#\,f_2)(0,0)=|C_{i\lambda,\varepsilon}|^2\,\int |1+\langle
t,\,\tau \rangle|_{\varepsilon}
^{-\frac{n+1}{2}+i\lambda}\,f_1(0,\tau)\,f_2(t,0)\ dt\,d\tau
\end{equation}
can be interpreted as the value at \,$(0,0)$\, of a convolution
product on \,$\R^n\times \R^n$\,, to wit that of the function \
$G_{i\lambda,\,\varepsilon}(t,\,\tau)=|C_{i\lambda,\varepsilon}|^2
\,|1+\langle t,\,\tau \rangle|_{\varepsilon}
^{-\frac{n+1}{2}+i\lambda}$\, by the function \ $(t,\,\tau)\mapsto
f_1(0,\tau)\,f_2(t,0)$\,. In a more expressive way, one may write
\begin{equation}
(f_1\,\#\,f_2)(0,0)=\l[\,({\mathcal F}\,G_{i\lambda,\,\varepsilon})
\l(\frac{1}{2i\pi}\,\frac{\partial}{\partial t},\,\frac{1}{2i\pi}\,
\frac{\partial}{\partial
\tau}\r)\,(f_1(0,\tau)\,f_2(t,0))\,\r](0,0)\,.
\end{equation}\\

We thus compute the (usual, not symplectic, this time) Fourier
transform of the distribution \ $(t,\tau)\mapsto |1+\langle t,\,\tau
\rangle|_{\varepsilon} ^{-\frac{n+1}{2}+i\lambda}$\, denoting as
\,$(\rho,\,r)$\, a pair of variables dual of \,$(t,\tau)$\,. With \
$\nu=\frac{n+1}{2}-i\lambda$\,, one has for real \,$R$\, the
equation
\begin{equation}
|R|_{\varepsilon}^{-\nu}=\frac{(-i)^{\varepsilon}}{2\pi}\,\Gamma(\frac{1-\varepsilon+\nu}{2})\,
\Gamma(\frac{1+\varepsilon-\nu}{2})\,[\,(0-i\,R)^{-\nu}+(-1)^{\varepsilon}\,(0+i\,R)^{-\nu}\,]\,.
\end{equation}
For \,$a>0$\,, one has
\begin{equation}
[a+i\,(1+\langle t,\,\tau
\rangle)]^{-\nu}=\frac{(2\pi)^{\nu}}{\Gamma(\nu)}\, \int_0^{\infty}
e^{-\frac{2\pi}{h}\,[a+i\,(1+\langle t,\,\tau
\rangle)]}\,h^{-\nu-1}\,dh
\end{equation}
and, since
\begin{equation}
{\mathcal F}(e^{-\frac{2i\pi}{h}\langle t,\,\tau
\rangle})(\rho,\,r)= h^n\,e^{2i\pi h\,\langle r,\,\rho \rangle}\,,
\end{equation}
one obtains
\begin{equation}
{\mathcal F}\,[0+i\,(1+\langle t,\,\tau
\rangle)]^{-\nu}={\mathrm{lim}}_{a\to 0}
\frac{(2\pi)^{\nu}}{\Gamma(\nu)}\int_0^{\infty}
h^{n-\nu-1}\,e^{-\frac{2\pi}{h}(a+i)}\, e^{2i\pi h\,\langle r,\,\rho
\rangle}\,dh\,.
\end{equation}
This is an integral giving, classically \cite[p.\,85]{mos}, the
function \,$K_{n-\nu}$\,, but since one of the exponents, to wit \
$2i\pi h\,\langle r,\,\rho \rangle$\,, is pure imaginary, it must be
interpreted as \ $-(0-2i\pi h\,\langle r,\,\rho \rangle)$\,. Being
careful with phases, one then finds
\begin{multline}
{\mathcal F}\,[0+i\,(1+\langle t,\,\tau
\rangle)]^{-\nu}=2^{\nu+1}\,\frac{\pi^{\nu}}{\Gamma(\nu)}\, |\langle
r,\,\rho
\rangle|^{\frac{\nu-n}{2}}\,\exp\,\l(\frac{i\pi\,(n-\nu)}{4}\,(1+{\mathrm{sign}}\,
\langle r,\,\rho \rangle)\r)\,\\
K_{n-\nu}\,\l(4\pi\,|\langle r,\,\rho
\rangle|^{\hf}\,\exp\,\l(\frac{i\pi}{4}\,(1-{\mathrm{sign}}\,
\langle r,\,\rho \rangle)\r)\r)\,.
\end{multline}\\

To finish the computation, one notes \cite[p.\,67]{mos} that, for
every \,$\mu\in \C$\,, one has if \,$x>0$\, the equations
\begin{align}
K_{\mu}(x\,e^{-\frac{i\pi}{2}})+K_{\mu}(x\,e^{\frac{i\pi}{2}})&=
\hf\,\Gamma(\frac{\mu}{2})\,\Gamma(\frac{2-\mu}{2})\,(-J_{\mu}(x)+J_{-\mu}(x))\,,\nonumber\\
K_{\mu}(x\,e^{-\frac{i\pi}{2}})-K_{\mu}(x\,e^{\frac{i\pi}{2}})
&=\frac{i}{2}\,\Gamma(\frac{1+\mu}{2})\,\Gamma(\frac{1-\mu}{2})\,
(J_{\mu}(x)+J_{-\mu}(x))\,.
\end{align}\\

Hence, the Fourier transform of the distribution \ $(t,\tau)\mapsto
|1+\langle t,\,\tau \rangle|_{\varepsilon}
^{-\frac{n+1}{2}+i\lambda}$\, is a distribution in the variables
\,$(r,\rho)$\, which coincides, when \ $\langle r,\,\rho \rangle
\neq 0$\,, with the function \,$H_{i\lambda,\,\varepsilon}(\langle
r,\,\rho \rangle)$\, defined by
\begin{equation}
H_{i\lambda,\,\varepsilon}(\beta)=\frac{2\,(2\pi)^{\frac{n+1}{2}-i\lambda}}
{\Gamma(\frac{n+1}{2}-i\lambda)}\,
\frac{\Gamma(\frac{n+3}{4}-\frac{i\lambda}{2}-\frac{\varepsilon}{2})\,
\Gamma(\frac{1-n}{4}+\frac{i\lambda}{2}+\frac{\varepsilon}{2})}
{\Gamma(\frac{n+1}{4}+\frac{i\lambda}{2}+\frac{\varepsilon}{2})\,
\Gamma(\frac{3-n}{4}-\frac{i\lambda}{2}-\frac{\varepsilon}{2})}\,\beta^{\frac{1-n}{4}-\frac{i\lambda}{2}}\,
K_{\frac{n-1}{2}+i\lambda}(4\pi\,\beta^{\hf})
\end{equation}
if \,$\beta>0$\,, and
\begin{multline}
H_{i\lambda,\,\varepsilon}(\beta)=
\Gamma(\frac{n+3}{4}-\frac{i\lambda}{2}-\frac{\varepsilon}{2})\,
\Gamma(\frac{1-n}{4}+\frac{i\lambda}{2}+\frac{\varepsilon}{2})\,
\Gamma(\frac{n-1}{4}+\frac{i\lambda}{2}+\frac{\varepsilon}{2})\,
\Gamma(\frac{5-n}{4}-\frac{i\lambda}{2}-\frac{\varepsilon}{2})\,\times\\
\frac{(2\pi)^{\frac{n-1}{2}-i\lambda}}{\Gamma(\frac{n+1}{2}-i\lambda)}\,
|\beta|^{\frac{1-n}{4}-\frac{i\lambda}{2}}\,
\l[\,J_{\frac{1-n}{2}-i\lambda}(4\pi\,|\beta|^{\hf})-(-1)^{\varepsilon}\,
J_{\frac{n-1}{2}+i\lambda}(4\pi\,|\beta|^{\hf})\,\r]
\end{multline}
if \,$\beta<0$\,.\\

On the other hand, according to (\ref{2Cilambda}), one has
\begin{equation}
|C_{i\lambda,\varepsilon}|^2=\pi^{-n}\,
\frac{\Gamma(\frac{1+n}{4}-\frac{i\lambda}{2}+\frac{\varepsilon}{2})\,
\Gamma(\frac{1+n}{4}+\frac{i\lambda}{2}+\frac{\varepsilon}{2})}
{\Gamma(\frac{1-n}{4}+\frac{i\lambda}{2}+\frac{\varepsilon}{2})\,
\Gamma(\frac{1-n}{4}-\frac{i\lambda}{2}+\frac{\varepsilon}{2})}\,.
\end{equation}\\
Using the duplication formula
\begin{equation}
\Gamma(\frac{1\pm n}{4}-\frac{i\lambda}{2}+\frac{\varepsilon}{2})\,
\Gamma(\frac{3\pm n}{4}-\frac{i\lambda}{2}-\frac{\varepsilon}{2})
=(2\pi)^{\hf}\,2^{\frac{1\mp n}{2}+i\lambda}\, \Gamma(\frac{1\pm
n}{2}-i\lambda)\,,
\end{equation}
one finally obtains the function
\,$F_{i\lambda,\,\varepsilon}(\beta)$\, as characterized
in the proposition.\\
\end{proof}
{\bf Remark.} Using the series expansion of the Bessel functions,
one obtains in all cases the following expression (involving
Pochhammer's symbols)
\begin{multline}\label{4expans}
F_{i\lambda,\,\varepsilon}(\beta)=\sum_{m\geq 0}
\frac{(4\pi^2\,\beta)^m}{m\,!\
(\frac{n+1}{2}+i\lambda)_m}\\
+(4\pi^2)^{\frac{1-n}{2}-i\lambda}\,\frac{\Gamma(\frac{n-1}{2}+i\lambda)}
{\Gamma(\frac{1-n}{2}-i\lambda)}\
|\beta|_{\varepsilon}^{\frac{1-n}{2}-i\lambda}\, \sum_{m\geq 0}
\frac{(4\pi^2\,\beta)^m}{m\,!\ (\frac{3-n}{2}-i\lambda)_m}\,.
\end{multline}
Plugging this expansion into the equation (\ref{4f1sharpf2}), or in
its generalized version
\begin{multline}
(f_1\,\#\,f_2)(s,\,\sigma)=\\
\l[\,F_{i\lambda,\,\varepsilon}\,\l(-\frac{1}{4\pi^2}\,\sum\frac{\partial^2}
{\partial t_j\,\partial \tau_j}\r)\,
\l(f_1(s,[{g_{s,\sigma}'}^{-1}]\,\tau)\,f_2([g_{s,\sigma}]\,t,\sigma)\r)\,\r]\,
((t,\,\tau)=(0,\,0))\,,
\end{multline}
one sees that the series which is first term of (\ref{4expans})
contributes to \,$(f_1\,\#\,f_2)(s,\,\sigma)$\, a series of
differential operators applied to the tensor product
\,$f_1\,\otimes\,f_2$\,, evaluated at \,$(s,\sigma)$\,. However, the
second term of (\ref{4expans}) cannot be neglected. It is a ramified
function of \,$\beta$\, at the origin, which shows that, in reality,
no convergent, or simply asymptotic, series of differential
operators applied  to the tensor product \,$f_1\,\otimes\,f_2$\,
evaluated at \,$(s,\sigma)$\, can produce a satisfactory
approximation of the sharp product \,$f_1\,\#\,f_2$\,. This is of
course in contrast with what happens with the usual (standard,
antistandard or Weyl) symbolic calculi on \,$\R^{n+1}$\,.\\

A fully similar phenomenon appeared in \cite{U2}, in relation with
the quantization of the upper half--plane \ $SL(2,\R)/SO(2)$\,:
indeed, the function \,$E(z)$\, that occurs there, in Theorem 5.1,
has a comparable type of singularity at the origin, to wit a
ramified part. The main difference is that the exponent of
\,$|\beta|$\,, to with \,$-i\lambda$\, in the present context, is to
be replaced, in the former reference, by \,$\lambda$\, so that, as
\,$\lambda$\, increases, the ramified term is, in some sense, pushed
away. This explains why, indeed --- as was shown in \cite{U2} ---
the symbolic calculus developed there has better properties for
increasing values of \,$\lambda$\,, and the reason why a limit of
the calculus as \,$\lambda \to \infty$\, could be found (the Fuchs
calculus). Nothing of the sort can work in the present situation,
which is worse in this respect. So far as series expansions with
respect to \,$\lambda^{-1}$\,  are concerned, they never occur in
the quantization of symmetric spaces, since the true functions of
\,$\lambda$\, involved always have
an essential singularity at infinity.\\

To put on end to this section, let us observe that, assuming that
the symbols \,$H_1$\, and \,$H_2$\, of two pseudodifferential
operators \,$A_1$\, and \,$A_1$\, on \,$\R^{n+1}$\, both satisfy the
equation (\ref{2J}) (which means that the standard and antistandard
symbols of each of the two operators under consideration agree) does
not imply that the same holds for the composition \,$A_1\,A_2$\,: it
suffices to consider the two symbols \,$x_1\,\xi_2$\, and
\,$x_2\,\xi_3$\,, the sharp composition of which, in the standard or
antistandard calculus, is respectively \
$x_1\,x_2\,\xi_2\,\xi_3+\frac{1}{2i\pi}\,x_1\,\xi_3$\, or \
$x_1\,x_2\,\xi_2\,\xi_3$\,. This makes it impossible to reduce the
study of the composition
of symbols in the \,$\pi_{i\lambda,\varepsilon}$--\,calculus to the usual one on \,$\R^{n+1}$\,.\\

Another reason, even more decisive, why pseudodifferential analysis
on projective space cannot be fully reduced (despite Proposition
2.2) to the (standard, or Weyl) pseudodifferential analysis on
\,$\R^{n+1}$\,, as currently developed, has to do with more
technical aspects: even though, for their applications to partial
differential equations, miscellaneous classes of symbols and
associated operators have been considered, translations of the phase
space always play a role there, albeit a local one: in contrast, it
is only through its action by (local) conjugations that the group
\,$GL(n+1,\R)$\, or, more generally (in the case of the Weyl
calculus), the symplectic group \,${\mathrm{Sp}}(n+1,\R)$\,, plays a
role in the definition of such classes of symbols. To give but one
example, using the space \,${\mathcal S}(\R^{n+1})$\, of
\,$C^{\infty}$\, vectors of the Heisenberg representation and its
dual space \,${\mathcal S}'(\R^{n+1})$\,, one finds immediately a
very large class of symbols, to wit \,${\mathcal S}'(\R^{2n+2})$\,,
all of which give rise to meaningful operators; but giving a
characterization of, say, the standard symbols of linear operators
from the space of \,$C^{\infty}$\, vectors of the quasiregular
representation of \,$GL(n+1,\R)$\, in \,$L^2(\R^{n+1})$\, to the
dual space is a problem in harmonic analysis
--- possibly a not too difficult one --- which P.D.E. people would probably find no reason to
consider.


\section{Some special symbols}

In this last section, we try to familiarize ourselves with the
calculus by an analysis of the operators the symbols of which are
integral powers (the exponents can be of any sign, but negative ones
are more interesting) of functions of the species
\,$\sc{a,\,x}\sc{b,\,\xi}$\,, with \ $\sc{a,\,b}=0$\,: these
functions already appeared in (\ref{2phiab}) in the case when \
$a,\,b\in\R^{n+1}$\, but, here, they will be complex vectors. In the
case when \,$n=1$\,(the study of which was made in \cite{U3}), the
symbols under consideration generate the discrete spaces of the
decomposition of \,$L^2(G_1/H_1)$\,: moreover, the Hilbert sum of
(one half of) these spaces is closed under the sharp product of
symbols, in the
\,${\mathrm{Op}}_{i\lambda,\varepsilon}$\,--\,calculus, and the
composition formulas were made explicit with the help of the
so--called Rankin--Cohen brackets. In the case when \,$n\geq 2$\,,
these functions do not lie in \,$L^2(G_n/H_n)$\, any more. Our
interest in them lie in the fact that they are the symbols of
integral powers, with positive exponents, of resolvents of
certain infinitesimal operators of the representation \,$\pi_{i\lambda,\varepsilon}$\,.  \\

We need to introduce  the infinitesimal operators of the
representation \,$\pi_{i\lambda,\varepsilon}$\,: these are defined
by the equation
\begin{equation}
(d\pi_{i\lambda,\varepsilon}(X)\,u)(s)=\frac{d}{dt}_{|_{t=0}}\,
(\pi_{i\lambda,\varepsilon}(\exp\,tX)\,u)(s)\,, \qquad X\in
{\mathfrak g}_n\,,
\end{equation}
where \,${\mathfrak g}_n$\, is the Lie algebra of \,$G_n$\,. As a
linear basis of \,${\mathfrak g}_n$\,, we choose the set \
$(E_{jk})_{(j,k)\neq (n+1,n+1)}$\, defined as follows: if \,$j\neq
k$\,, \,$E_{jk}=e_j\otimes e_k^*$\, is the matrix such that
${(E_{jk})}_{\ell,m}=\delta^j_{\ell}\,\delta^k_m$\,; next, \
$E_{jj}$\, is the diagonal matrix with diagonal
\,$\{0,\dots,0,1,0,\dots,0,-1\}$\, where the \,$1$\, occupies the
\,$j^{\,\mathrm{th}}$\, place. A fixture of the developments to come
will be the operator \ $\sum_m s_m\frac{\partial}{\partial
s_m}+\frac{n+1}2+i\lambda$\,. It is convenient to set
\begin{equation}
\mu=\frac{n+1}2+i\lambda\,,\qquad D_{\mu}=\langle
s,\,\frac{\partial}{\partial s}\rangle +\mu\,.
\end{equation}
Applying (\ref{2transform})--(\ref{2defpi}), one finds the equations
\begin{align}\label{5infoprep}
&d\pi_{i\lambda,\varepsilon}(E_{j,n+1})=-\frac{\partial}{\partial
s_j}\,,\qquad\qquad\qquad
\ \qquad\qquad d\pi_{i\lambda,\varepsilon}(E_{n+1,k})=s_k\,D_{\mu}\,,\nonumber\\
&d\pi_{i\lambda,\varepsilon}(E_{jk})=-s_k\frac{\partial}{\partial
s_j}\quad(j,k,n+1\ {\mathrm{distinct}})\,, \,\quad
d\pi_{i\lambda,\varepsilon}(E_{jj})=-D_{\mu}-s_j\frac{\partial}{\partial
s_j}\,:
\end{align}
note that, so far as the {\em{formal}\/} infinitesimal operators
only are considered, there is no difference between the
representations associated with the same value of \,$\lambda$\, but
distinct values of \,$\varepsilon$\,: there is of course a
considerable difference when the self--adjoint extensions of the
operators under consideration are concerned. We also denote as
\,$d\pi_{i\lambda,\varepsilon}$\, the extension of this map to the
enveloping algebra
\,${\mathcal U}({\mathfrak g}_n)$\, of \,${\mathfrak g}_n$\,.\\

Under the assumption that not only some operator
\,${\mathrm{Op}}_{i\lambda,\varepsilon}(f)$\, but also the result of
its composition on the left by the image, under
\,$d\pi_{i\lambda,\varepsilon}$\,, of any element of \,${\mathcal
U}({\mathfrak g}_n)$\,,  is a Hilbert--Schmidt endomorphism of the
space \,$L^2(\R^n)$\,, one can, with the help of (\ref{2stand}) and
of the preceding equations, compute the symbol of the operator \
$d\pi_{i\lambda,\varepsilon}(X)\,{\mathrm{Op}}_{i\lambda,\varepsilon}(f)$\,
for any vector \,$X\in {\mathfrak g}$\,, getting after a trivial
computation the set of relations (in which \ $j,\,k\neq n+1$)
\begin{align}\label{5leftinfOp}
d\pi_{i\lambda,\varepsilon}(E_{j,n+1})\,{\mathrm{Op}}_{i\lambda,\varepsilon}(f)&=
{\mathrm{Op}}_{i\lambda,\varepsilon}\l(-\frac{\partial f}{\partial
s_j}
+\mu\,\frac{\sigma_j}{1+\langle s,\sigma \rangle}\,f\r)  \,,\nonumber\\
d\pi_{i\lambda,\varepsilon}(E_{n+1,k})\,{\mathrm{Op}}_{i\lambda,\varepsilon}(f)&=
{\mathrm{Op}}_{i\lambda,\varepsilon}\l(s_k\,\langle
s,\,\frac{\partial f}{\partial s}\rangle
+\mu\,\frac{s_k}{1+\langle s,\sigma \rangle}\,f\r)  \,,\nonumber\\
d\pi_{i\lambda,\varepsilon}(E_{jk})\,{\mathrm{Op}}_{i\lambda,\varepsilon}(f)&=
{\mathrm{Op}}_{i\lambda,\varepsilon}\l(-s_k\,\frac{\partial
f}{\partial s_j}
+\mu\,\frac{s_k\,\sigma_j}{1+\langle s,\sigma \rangle}\,f\r) \,,\nonumber\\
d\pi_{i\lambda,\varepsilon}(E_{jj})\,{\mathrm{Op}}_{i\lambda,\varepsilon}(f)&=
{\mathrm{Op}}_{i\lambda,\varepsilon}\l(-\langle s,\,\frac {\partial
f}{\partial s}\rangle-s_j\,\frac{\partial f}{\partial s_j}
+\mu\,\frac{s_j\,\sigma_j-1}{1+\langle s,\sigma \rangle}\,f\r) \,.
\end{align}\\

Our main concern, in this section, has to do with the operators the
symbols of which are integral powers of the function
\begin{equation}\label{5phialphabeta}
\phi_{a,b}((x,\,\xi)^{\bullet})=\langle a,\,x\rangle\,\langle
b,\,\xi\rangle\,,\qquad\quad (x,\,\xi)^{\bullet}\in {\mathcal
X}_n^{\bullet}
\end{equation}
or, in inhomogeneous coordinates, with \
$a=(a_1,\dots,a_{n+1})=(a_*,a_{n+1})$\, and \ $b=(b_*,b_{n+1})$\,,
\begin{equation}\label{5phialphabetainh}
\phi_{a,b}(s,\,\sigma)=\frac {(a_{n+1}+\langle
a_*,\,s\rangle)\,(b_{n+1}+\langle b_*,\,\sigma\rangle)} {1+\langle
s,\,\sigma\rangle}\,.
\end{equation}
It is assumed that \,$a$\, and \,$b$\, lie in \,$\C^{n+1}$\,: when
these two vectors are real, this function has already been
considered in (\ref{2phiab}). The case when \ $\sc{a,\,b}=0$\, will
be of special interest. Since the equations (\ref{5leftinfOp}) give
in particular (setting \,$f=1$\,) the symbols of the infinitesimal
operators of the representation \,$\pi_{i\lambda,\varepsilon}$\,,
one can verify that, {\em in this case\/}, the function
\,$\phi_{a,b}$\, is the symbol of the operator \
$\mu^{-1}\,d\pi_{i\lambda,\varepsilon}\, (X_{a,b})$\,, with
\begin{equation}\label{5infinitop}
X_{a,b}=\sum_{(j,\,k)\neq (n+1,n+1)} a_k\,b_j\,E_{jk}\,.
\end{equation}\\

We first make a quick study of the operator with symbol
\,$\phi_{a,b}^p$\, with \,$p=0,1,\dots$\,: we are more interested in
the same symbols with \,$p=-1,-2,\dots$\,, but this will require
some preparation. Even in the case when \,$p\in \N$\,, this symbol
is associated to the function \ $(x,\,\xi)\mapsto (\langle
a,\,x\rangle\,\langle b,\,\xi\rangle)^p$\, on \,$\R^{n+1}\times
\R^{n+1}$\,, certainly not a bounded function so that Proposition
2.2 does not apply, even though the application of differential
operators such as \,$x_j\,\frac{\partial}{\partial x_k}$\, or
\,$\xi_j\,\frac{\partial}{\partial \xi_k}$\, does not make the
symbol any worse. In this section, we shall sometimes extend the
meaning of \,${\mathrm{Op}}_{i\lambda,\varepsilon}$\, beyond the
domain in which full justifications have been carried, keeping in
mind that the following basic property of the calculus should
continue to hold: in the case when a symbol \,$f_1$\, depends only
on \,$s$\,, or when \,$f_2$\, depends only on \,$\sigma$\,, the
product \,$f_1\,f_2$\, must be the symbol of the composition
\,${\mathrm{Op}}_{i\lambda,\varepsilon}(f_1)\,
{\mathrm{Op}}_{i\lambda,\varepsilon}(f_2)$\,; of course, the
situation has to be reversed when dealing with the
\,${\mathrm{Op}}^{\vee}_{i\lambda,\varepsilon}$\,--\,calculus. Also,
we shall
take advantage of the equations (\ref{5leftinfOp}).\\

First, let us deal with powers of the symbol \ $1+\sc{s,\,\sigma}$\,.\\
\begin{lemma}
With \ $\mu=\frac{n+1}{2}+i\lambda$\,, consider the operator \
$D_{\mu}=\langle s,\, \frac{\partial}{\partial s}\rangle+\mu$\,, an
endomorphism of the space \,$H^{\infty}_{i\lambda,\varepsilon}$\, of
\,$C^{\infty}$\, vectors of the representation
\,$\pi_{i\lambda,\varepsilon}$\,, and set
\begin{equation}
\Lambda_p={\mathrm{Op}}_{i\lambda,\varepsilon}\,((1+\langle
s,\,\sigma \rangle)^{-p})\,.
\end{equation}
For \ $p=0,1,\dots$\,, one has
\begin{equation}\label{5lambdap}
\Lambda_{p}=\frac{D_{\mu}(D_{\mu}+1)\dots
(D_{\mu}+p-1)}{\mu(\mu+1)\dots (\mu+p-1)}\,.
\end{equation}\\
\end{lemma}

\begin{proof}
Set \ $D_0=\langle s,\,\frac{\partial}{\partial s}\rangle$\,. From
(\ref{5leftinfOp}), applied with \,$f=1$\,, one has
\begin{equation}\label{5dpiD0}
d\pi_{i\lambda,\varepsilon}\,(\sum E_{jj})=-(n+1)\,D_0-n\,\mu\,,
\end{equation}
so that \,$D_{\mu}$\, is indeed an endomorphism of the space
\,$H^{\infty}_{i\lambda,\varepsilon}$\,. On the other hand, by
(\ref{5leftinfOp}) again, given any symbol \,$f$\,, the symbol of
the operator \,$d\pi_{i\lambda,\varepsilon}\,(\sum
E_{jj})\,{\mathrm{Op}}_{i\lambda,\varepsilon}(f)$\, is the function
\begin{equation}\label{5symb}
-(n+1)\,\langle s,\,\frac{\partial f}{\partial s}\rangle
+\mu\,\frac{\langle s,\,\sigma\rangle-n}{1+\langle
s,\,\sigma\rangle}\,f\,:
\end{equation}
this leads to the equation
\begin{equation}
d\pi_{i\lambda,\varepsilon}\,(\sum\,E_{jj})\,\Lambda_p=(\mu+(n+1)\,p)\,\Lambda_p-(n+1)(p+\mu)\,
\Lambda_{p+1}\,,
\end{equation}
which can also be written, using (\ref{5dpiD0}), as
\begin{equation}\label{5lambdapp+1}
\l(\,\langle s,\,\frac{\partial f}{\partial
s}\rangle+\mu+p\,\r)\,\Lambda_p=(\mu+p)\, \Lambda_{p+1}\,,
\end{equation}
from which (\ref{5lambdap}) follows.\\
\end{proof}

We need to introduce some right inverse \,$D_{\mu}^{-1}$\, of
\,$D_{\mu}$\, and more generally, for later purposes, a resolvent
\,$(D_{\mu}-\rho)^{-1}$\, of this operator: we define it by the
equation
\begin{equation}\label{5invDmu}
((D_{\mu}-\rho)^{-1}u)(s)=\int_0^1 u(ts)\,t^{\mu-\rho-1}\,dt\,,
\end{equation}
and observe first (this is one of the so--called Hardy's
inequalities) that, when \,$\Re\rho<\hf$\,, it extends as a bounded
operator on \,$L^2(\R^n)$\,: indeed, it suffices to write
\begin{align}
(v\,|\,((D_{\mu}-\rho)^{-1}u)&=-\int_{\R^n}
\bar{v}(s)\,ds\,\int_1^{\infty}u(ts)\,t^{\mu-\rho-1}\,dt
\nonumber\\
&=-\int_1^{\infty} t^{\mu-\rho-1}\,dt\,\int_{\R^n}
\bar{v}(s)\,\,u(ts)\,ds\,,
\end{align}
where the last integral, by the Cauchy--Schwarz inequality, is bounded by\\
$t^{-\frac{n}{2}}\,\Vert v\Vert_{L^2(\R^n)}\,\Vert
u\Vert_{L^2(\R^n)}$\,. Note that, even for large values of
\,$\Re\rho$\,, the integral (\ref{5invDmu}) makes sense when \,$u$\,
is flat
enough at \,$s=0$\,.\\

\begin{lemma}
Given \,$a$\, and  \,$b\in \C^{n+1}$\,, the operators with symbols \
$(s,\sigma)\mapsto a_{n+1}+\langle a_*,\,s\rangle$\, and \
$(s,\sigma)\mapsto \mu\,\frac{b_{n+1}+\langle b_*,\,\sigma\rangle}
{1+\langle s,\,\sigma \rangle}$\, are respectively the operator
\,$S_{a}$\, of multiplication by the function \ $s\mapsto
a_{n+1}+\langle a_*,\,s\rangle$\, and the operator
\begin{equation}\label{5Tb}
T_{b}=b_{n+1}\,D_{\mu}-\langle b_*,\,\frac{\partial}{\partial
s}\rangle\,,
\end{equation}
where \ $D_{\mu}=\langle s,\,\frac{\partial}{\partial
s}\rangle+\mu$\,. In the case when \ $\langle a,\,b \rangle=0$\,,
the two operators under consideration
generate a Lie algebra isomorphic to that of the one--dimensional affine group.\\
\end{lemma}

\begin{proof}
From (\ref{5leftinfOp}), then (\ref{5infoprep}),
\begin{equation}
{\mathrm{Op}}_{i\lambda,\varepsilon}\,\l(\,\mu\,\frac{\langle
b_*,\,\sigma\rangle} {1+\langle s,\,\sigma \rangle}\,\r)=\sum_{j=
1}^n\, b_j\,d\pi_{i\lambda,\varepsilon}(E_{j,n+1})=-\langle
b_*\,,\,\frac{\partial}{\partial s}\rangle\,;
\end{equation}
on the other hand, it has been found in Lemma 5.1 that \ $
{\mathrm{Op}}_{i\lambda,\varepsilon}\,(\frac{\mu}{1+\langle
s,\,\sigma \rangle})=D_{\mu}$\,. Next, it is immediate to verify
that
\begin{align}\label{5commut}
[\,T_{b}\,,\,S_{a}\,]&=b_{n+1}\,\langle a_*,\,s\rangle-
\langle a_*,\,b_*\rangle\nonumber\\
&=b_{n+1}\,S_{a}-\langle a,\,b\rangle\,.
\end{align}\\
\end{proof}

We have recalled that the symbol of an operator such as
\,${\mathrm{Op}}_{i\lambda,\varepsilon}(f_1)\,{\mathrm{Op}}_{i\lambda,\varepsilon}(f_2)$\,
reduces to \,$f_1\,f_2$\, whenever the symbol \,$f_1$\, on the
left--hand side depends only on the variable \,$s$\,, or when the
symbol \,$f_2$\, on the right--hand side depends only on the
variable \,$\sigma$\,. This has the consequence that, with the
notation introduced in Lemma 5.3, one has
\begin{equation}
{\mathrm{Op}}_{i\lambda,\varepsilon}\,(\phi_{a,b})=\mu^{-1}\,S_{a}\,T_{b}\,.
\end{equation}
The operator \,$T_{b}$\, can be decomposed further as the product
\begin{equation}
T_{b}={\mathrm{Op}}_{i\lambda,\varepsilon}\,(\frac{\mu}{1+\langle
s,\,\sigma \rangle})\,
{\mathrm{Op}}_{i\lambda,\varepsilon}\,(b_{n+1}+\langle b_*,\,\sigma
\rangle)\,:
\end{equation}
with the help of Lemma 5.1, this leads to the equation
\begin{equation}
T_{b}=D_{\mu}\,
{\mathrm{Op}}_{i\lambda,\varepsilon}\,(b_{n+1}+\langle b_*,\,\sigma
\rangle)\,,
\end{equation}
which can be inverted as
\begin{equation}\label{5ssigma}
{\mathrm{Op}}_{i\lambda,\varepsilon}\,(b_{n+1}+\langle b_*,\,\sigma
\rangle)= D_{\mu}^{-1}\,T_{b}\,:
\end{equation}

In all that precedes, the condition \ $\langle a,\,b\rangle=0$\, was
not needed. It is, however, needed in the proposition follows, which
stresses the ``reproducing'' property of the symbols
\,$\phi_{a,b}$\, under consideration.\\

\begin{proposition}
Under the assumption that  \ $\langle a,\,b\rangle=0$\,, one has,
for \ $p=1,2,\dots$\,,
\begin{equation}
({\mathrm{Op}}_{i\lambda,\varepsilon}\,(\phi_{a,b}))^p=
\frac{\mu(\mu+1)\dots (\mu+p-1)}{\mu^p}\,
{\mathrm{Op}}_{i\lambda,\varepsilon}\,(\phi^p_{a,b})\,.
\end{equation}\\
\end{proposition}

\begin{proof}
We abbreviate in this proof \,$S_a$\, and \,$T_b$\, as \,$S$\, and
\,$T$\,. It is no loss of generality to assume that \,$b_{n+1}=1$\,
(in the case when \,$b_{n+1}=0$\,, one may use to that effect the
covariance of the calculus), so as to simplify the commutation
relation (\ref{5commut}). This immediately leads to
\begin{equation}
({\mathrm{Op}}_{i\lambda,\varepsilon}\,(\phi_{a,b}))^p=
(\mu^{-1}\,S\,T)^p=\mu^{-p}\,S^p\ T(T+1)\dots(T+p-1)\,.
\end{equation}
On the other hand, using the basic property of the calculus, Lemma
5.1 and (\ref{5ssigma}), one obtains
\begin{equation}
{\mathrm{Op}}_{i\lambda,\varepsilon}\,(\phi^p_{a,b})=S^p\,
\frac{D_{\mu}(D_{\mu}+1)\dots (D_{\mu}+p-1)}{\mu(\mu+1)\dots
(\mu+p-1)}\, [\,D_{\mu}^{-1}\,T\,]^p\,.
\end{equation}
The equation to be shown thus reduces to
\begin{equation}\label{5tobeshown}
T(T+1)\dots(T+p-1)=D_{\mu}(D_{\mu}+1)\dots
(D_{\mu}+p-1)\,[\,D_{\mu}^{-1}\,T\,]^p\,.
\end{equation}
One first verifies the commutation relation
\begin{equation}
[D_{\mu},\,T]=D_{\mu}-T\,,
\end{equation}
from which one gets
\begin{equation}
(T+1)\,D_{\mu}=(D_{\mu}+1)\,T
\end{equation}
and, by induction,
\begin{equation}\label{5step}
(T+k+1)\,(D_{\mu}+k)=(D_{\mu}+k+1)\,(T+k)
\end{equation}
for \,$k=0,1,\dots$\,. First simplifying the right--hand side of
(\ref{5tobeshown}) as
\begin{equation}
(D_{\mu}+p-1)\dots (D_{\mu}+1)\,T\ [\,D_{\mu}^{-1}\,T\,]^{p-1}\,,
\end{equation}
we show by induction on \,$k\ (0\leq k \leq p-1$\,) that it can also
be written as
\begin{equation}
(D_{\mu}+p-1)\dots (D_{\mu}+k+1)\,(T+k)\dots T\
[\,D_{\mu}^{-1}\,T\,]^{p-k-1}
\end{equation}
(an expression that reduces to the left--hand side of
(\ref{5tobeshown}) when \,$k=p-1$): the step from \,$k$\, to
\,$k+1$\, is managed with the help  of (\ref{5step}), writing
\begin{align}
(D_{\mu}+k+1)\,(T+k)\dots T&=(T+k+1)\,(D_{\mu}+k)\,(T+k-1)\dots T\nonumber\\
&=(T+k+1)(T+k)\,(D_{\mu}+k-1)\,(T+k-2)\dots T\nonumber\\
&=\dots \nonumber\\
&=(T+k+1)\dots (T+1)\,D_{\mu}\,.
\end{align}\\
\end{proof}

\begin{corollary}
Let \,$a,\,b\in\C^{n+1}$\, satisfy \ $\sc{a,\,b}=0$\,. When \
$f=\phi^p_{a,b}\,,\ p\in \N$\,, the equation \
${\mathrm{Op}}_{i\lambda,\varepsilon}
^{\vee}(f)={\mathrm{Op}}_{i\lambda,\varepsilon}(J_{i\lambda,\varepsilon}\,f)$\,
from Corollary {\em 2.4\/} extends, only replacing, in the
expression {\em(\ref{2formdel0})\/} of the function
\,$G_{i\lambda,\varepsilon}\,(\rho,\delta)$\,, \,$\rho$\, by \,$p$\,
and \,$\delta$\, by \,$p$\,
mod $2$\,.\\
\end{corollary}

\begin{proof}
Assuming without loss of generality that \,$b_{n+1}=1$\,, we first
note that \ $S_{\bar{a}}^*=S_a$\, and \ $T_{\bar{b}}^*=1-T_b$\,, and
that \ $(1-T_b)\,S_a=-S_a\,T_b$\,: as a consequence, starting from
the equation \
${\mathrm{Op}}_{i\lambda,\varepsilon}(\phi_{a,b})=\mu^{-1}\,S_a\,T_b$\,,
\begin{equation}
{\mathrm{Op}}_{i\lambda,\varepsilon}^{\vee}(\phi_{a,b})=
{\mathrm{Op}}_{i\lambda,\varepsilon}(\phi_{\bar{a},\bar{b}})^*=
\bar{\mu}^{-1}\,(1-T_b),S_a=-\bar{\mu}^{-1}\,S_a\,T_b\,,
\end{equation}
so that
\begin{equation}
{\mathrm{Op}}_{i\lambda,\varepsilon}^{\vee}(\phi_{a,b})=-\frac{\mu}{\bar{\mu}}\,
{\mathrm{Op}}_{i\lambda,\varepsilon}(\phi_{a,b})\,.
\end{equation}
We then obtain from Proposition 5.4 that
\begin{align}\label{5OpOpvee}
{\mathrm{Op}}_{i\lambda,\varepsilon}^{\vee}(\phi^p_{a,b})=
{\mathrm{Op}}_{i\lambda,\varepsilon}(\phi^p_{\bar{a},\bar{b}})^*&=
\frac{\bar{\mu}^p}{\bar{\mu}\dots(\bar{\mu}+p-1)}\,
[{\mathrm{Op}}_{i\lambda,\varepsilon}(\phi_{\bar{a},\bar{b}})^*]^p\nonumber\\
&=(-1)^p\,\frac{\mu(\mu+1)\dots(\mu+p-1)}{\bar{\mu}(\bar{\mu}+1)\dots(\bar{\mu}+p-1)}
\ {\mathrm{Op}}_{i\lambda,\varepsilon}(\phi^p_{a,b})\,.
\end{align}
On the other hand, going back to (\ref{2formdel0}), one finds
\begin{equation}
G_{i\lambda,\varepsilon}(p,\,p\ {\mathrm{mod}}\
2)=(-1)^p\,\frac{F(\mu)}{F(\bar{\mu})}
\end{equation}
with
\begin{equation}
F(\mu)=\frac{\Gamma(\frac{1-\mu+\varepsilon}{2})}{\Gamma(\frac{1-\mu-p+|\varepsilon-\delta|}{2})}
\,\times\,\frac{\Gamma(\frac{\mu+p+|\varepsilon-\delta|}{2})}{\Gamma(\frac{\mu+\varepsilon}{2})}\,:
\end{equation}
since both numbers \ $p\pm |\varepsilon-\delta|-\varepsilon$\, are
non--negative even numbers, one may interpret each of the two
factors above as a Pochhammer's symbol, which leads after a
case--by--case computation to the equation
\begin{equation}
F(\mu)=(-1)^p\,\frac{\mu(\mu+1)\dots (\mu+p-1)}{2^p}\,,
\end{equation}
and to the identification of the coefficient in front of the
right--hand side of (\ref{5OpOpvee})
with \ $G_{i\lambda,\varepsilon}(p,\,p\ {\mathrm{mod}} \ 2)$\,.\\
\end{proof}

As a preparation towards some calculations related to the
composition of operators with certain special symbols, we compute
the symbol of the resolvent operator \,$(D_{\mu}-\rho)^{-1}$\,.
Since, according the definition (\ref{5invDmu}), one has
\begin{equation}
(D_{\mu}-\rho)^{-1}=\int_0^1 t^{\sc{s,\,\frac{\partial}{\partial
s}}+\mu-\rho-1}\,dt\,,
\end{equation}
and since, from (\ref{2Op}), the symbol, in the
\,${\mathrm{Op}}_{i\lambda,\varepsilon}$\,--\,calculus, of the
operator \ $t^{\sc{s,\,\frac{\partial}{\partial s}}}$\,, is
immediately seen to be the function
\begin{equation}
f_t(s,\,\sigma)=\frac{|1+\sc{s,\,\sigma}|_{\varepsilon}^{\mu}}
{|1+t\,\sc{s,\,\sigma}|_{\varepsilon}^{\mu}}\,,
\end{equation}
the symbol of the operator \,$(D_{\mu}-\rho)^{-1}$\, is the function
\begin{equation}
h_{\rho}(s,\,\sigma)=|1+\sc{s,\,\sigma}|_{\varepsilon}^{\mu}\,
\int_0^1
 |1+t\,\sc{s,\,\sigma}|_{\varepsilon}^{-\mu}\ t^{\mu-\rho-1}\,dt\,:
\end{equation}
there is no need to display the elementary calculations, based on
the splitting of the integral into two parts in the case when \
$\sc{s,\,\sigma}<-1$\,, which lead to the explicit formula
\begin{equation}
h_{\rho}(s,\,\sigma)= (\mu-\rho)^{-1}\,(1+\sc{s,\,\sigma})^{\mu}\
{}_2\!F_1\,(\mu,\,\mu-\rho;\,\mu+1-\rho;\,- \sc{s,\,\sigma})
\end{equation}
when \ $\sc{s,\,\sigma}>-1$\,, and to
\begin{multline}
h_{\rho}(s,\,\sigma)=(-1)^{\varepsilon}\,\frac{\Gamma(\mu-\rho)\,
\Gamma(1-\mu)}{\Gamma(1-\rho)}\,|\sc{s,\,\sigma}|^{\rho-\mu}\,|1+\sc{s,\,\sigma}|^{\mu}\,\\
+(1-\mu)^{-1}\,|\sc{s,\,\sigma}|^{\rho-\mu}\,|1+\sc{s,\,\sigma}| \
\times \ {}_2\!F_1\,(\rho+1-\mu,\,1-\mu;\,2-\mu;\,1+\sc{s,\,\sigma})
\end{multline}
when \ $\sc{s,\,\sigma}<-1$\,. Note that, when \,$\rho=1,2,\dots$\,,
\,$h_{\rho}(s,\,\sigma)$\, reduces to a polynomial of degree
\,$\rho-1$\, in \,$\sc{s,\,\sigma}$\,, independent of
\,$\varepsilon$\,. One may mention the following formula (more
easily verified with the help of (\ref{5leftinfOp})): for
\,$p=0,1,\dots$\,,
\begin{equation}\label{5DmuPoch}
(D_{\mu}-1)^{-1}(D_{\mu}-2)^{-1}\dots(D_{\mu}-p)^{-1}=
{\mathrm{Op}}_{i\lambda,\varepsilon}\, \l(\,\frac{(1+\langle
s,\,\sigma \rangle)^{p}}{(\mu-1)(\mu-2)\dots(\mu-p)}\r)\,;
\end{equation}
the factors \,$(D_{\mu}-j)^{-1}$\, on the left--hand side do not act
within the space \,$L^2(\R^n)$\,, but their composition still makes
sense if interpreted (using a decomposition into simple elements) as
the sum \
\,$\sum_{j=1}^p\frac{(-1)^{p-j}}{(j-1)\,!\,(p-j)\,!}\,(D-j)^{-1}$\,.\\

The following lemma, in which we allow \,$i\lambda$\, to be replaced
by a complex number no longer pure imaginary, will be needed soon.
Note that, if \,$v\in {\mathcal S}(\R^n)$\,, the function \
$\rho\mapsto v^{\flat}_{\rho,\varepsilon}$\, extends as a
meromorphic function, valued in the space \,$C^{\infty}(\R^n)$\,, in
the whole complex plane, with simple poles only at points
\,$-(\frac{n+1}{2}+k)$\,, where \,$k=0,\,1,\dots$\, and \,$k\equiv
\varepsilon$\, mod $2$\,: this also makes it possible to define the
intertwining operator
\,$\theta_{\rho,\varepsilon}$\, in general.\\

\begin{lemma}
On functions in \,$\R^n$\, with a parity associated to \,$\delta$\,,
one has
\begin{equation}\label{5conntwoint}
\theta_{i\lambda,\varepsilon}=\theta_{i\lambda+1,\varepsilon}\,\times\,
\frac{1}{\pi}\,\frac{\Gamma(\frac{1-D_{\mu}+|\varepsilon-\delta|}{2})\,
\Gamma(\frac{1+D_{\mu}+|\varepsilon-\delta|}{2})}{\Gamma(\frac{D_{\mu}+|\varepsilon-\delta|}{2})\,
\Gamma(\frac{-D_{\mu}+|\varepsilon-\delta|}{2})}\,.
\end{equation}\\
\end{lemma}

\begin{proof}
We may pretend that we are testing both sides of the identity to be
proven on a given function \,$s\mapsto v^{\flat}(s)$\, homogeneous
of degree and parity \,$(-\frac{n}{2}-i\,\nu,\,\delta)$\,, keeping
in mind, however, that we really deal with {\em nice\/} integral
superpositions of such functions. Under the map (\ref{2charac}),
such a function \,$v^{\flat}$\, lifts to \,$\R^{n+1}$\, as the
function
\begin{equation}
v(x)=|x_{n+1}|_{|\varepsilon-\delta|}^{-\hf+i(\nu-\lambda)}\,v^{\flat}(x_*)\,,
\end{equation}
the Fourier transform of which is
\begin{multline}
({\mathcal
F}^{(n)}v)(x)=(-1)^{|\varepsilon-\delta|}\,\pi^{i(\lambda-\nu)}\,\frac
{\Gamma(\frac{1}{4}+\frac{i\,(\nu-\lambda)}{2}+\frac{|\varepsilon-\delta|}{2})}
{\Gamma(\frac{1}{4}-\frac{i\,(\nu-\lambda)}{2}+\frac{|\varepsilon-\delta|}{2})}
\ |x_{n+1}|_{|\varepsilon-\delta|}^{-\hf-i(\nu-\lambda)}\\
({\mathcal F}^{(n-1)}v^{\flat})(x_*)\,:
\end{multline}
it has been deemed prudent, here, to emphasize, as a superscript,
the dimension of the Fourier transform under consideration: the same
precaution will be used, presently, in connection with the
intertwining operators or quantizing maps \,${\mathrm{Op}}$\, to be
considered, as well as when using the constants
\,$C_{i\lambda,\varepsilon}$\, as defined in (\ref{2Cilambda}).
Then, (\ref{2deftheta}) yields
\begin{equation}
(\theta_{i\lambda,\varepsilon}\,v^{\flat})(\sigma)=(-1)^{|\varepsilon-\delta|}\,\pi^{i(\lambda-\nu)}\,\frac
{\Gamma(\frac{1}{4}+\frac{i\,(\nu-\lambda)}{2}+\frac{|\varepsilon-\delta|}{2})}
{\Gamma(\frac{1}{4}-\frac{i\,(\nu-\lambda)}{2}+\frac{|\varepsilon-\delta|}{2})}
\ ({\mathcal F}^{(n-1)}v^{\flat})(\sigma)\,,
\end{equation}
from which the equation (\ref{5conntwoint}) follows, once it has
been observed that, on a function such as \,$v^{\flat}$\,,
\,$D_{\mu}$\, reduces to the multiplication by \
$\hf+i\,(\lambda-\nu)$\,.\\
\end{proof}

We shall also need the following lemma, in which the variable
\,$s\in \R^n$\, is split as \,$s=(s_1,\,s_*)\in \R\times
\R^{n-1}$\,: note that the subscript \,$*$\, here concerns the
{\em last\/} \,$n-1$\, variables.\\

\begin{lemma}
If \ $u=u(s_1,\,s_*)$\, is homogeneous of degree and parity
\,$(\frac{1-n}{2}-i\nu,\,\delta)$\, with respect to the variables
\,$s_*$\,, one has
\begin{equation}
(\theta_{i\lambda,\varepsilon}^{(n)}\,u)(\sigma_1,\,\sigma_*)=
(-1)^{|\varepsilon-\delta|}\,C_{i(\nu-\lambda),|\varepsilon-\delta|}^{(1)}
\
\l(\,(\theta_{i(\lambda-\nu),|\varepsilon-\delta|}^{(1)}\,\otimes\,
\theta_{i\lambda+\hf,\varepsilon}^{(n-1)})\,u\,\r)\,(\sigma_1,\,\sigma_*)\,.
\end{equation}\\
\end{lemma}

\begin{proof}
Even though the genuine proof depends again on the lifting,
depending on \,$(i\lambda,\varepsilon)$\,, from functions on
\,$\R^n$\, to homogeneous functions on \,$\R^{n+1}$\,, we shall
satisfy ourselves with a shorter formal proof based on
(\ref{2thetaugly}). Starting from this equation, performing the
change of variables \ $s_*\mapsto (1+s_1\sigma_1)\,s_*$\, and using
the homogeneity, one obtains
\begin{multline}
(\theta_{i\lambda,\varepsilon}^{(n)}\,u)(\sigma_1,\,\sigma_*)=\\
C_{i\lambda,\varepsilon}^{(n)}\, \intR
|1+s_1\sigma_1|^{i(\lambda-\nu)-1}_{|\varepsilon-\delta|}\,ds_1\,\int_{\R^{n-1}}
|1+\sc{s_*,\,\sigma_*}|_{\varepsilon}^{-\frac{n+1}{2}+i\lambda}\,u(s_1,\,s_*)\,ds_*\\
=\frac{C_{i\lambda,\varepsilon}^{(n)}}{C^{(1)}_{i(\nu-\lambda),|\varepsilon-\delta|}\,
C_{i\lambda-\hf,\varepsilon}^{(n-1)}} \
\l(\,(\theta_{i(\lambda-\nu),|\varepsilon-\delta|}^{(1)}\,\otimes\,
\theta_{i\lambda-\hf,\varepsilon}^{(n-1)})\,u\,\r)\,(\sigma_1,\,\sigma_*)\,.
\end{multline}
From Lemma 5.5, we may substitute for \
$\theta_{i\lambda-\hf,\varepsilon}^{(n-1)}$\, the product of \
$\theta_{i\lambda+\hf,\varepsilon}^{(n-1)}$\, by the number
\begin{multline}
\frac{1}{\pi}\,\frac{\Gamma(\hf+\frac{i\,(\nu-\lambda)}{2}+\frac{|\varepsilon-\delta|}{2})
\,\Gamma(\hf-\frac{i\,(\nu-\lambda)}{2}+\frac{|\varepsilon-\delta|}{2})}
{\Gamma(-\frac{i\,(\nu-\lambda)}{2}+\frac{|\varepsilon-\delta|}{2})
\,\Gamma(\frac{i\,(\nu-\lambda)}{2}+\frac{|\varepsilon-\delta|}{2})}\\
=(-1)^{|\varepsilon-\delta|}\,C^{(1)}_{i(\lambda-\nu),|\varepsilon-\delta|}\,
C^{(1)}_{i(\nu-\lambda),|\varepsilon-\delta|}\,,
\end{multline}
which leads to the result indicated, if one notes also that
\begin{equation}\label{5idconst}
C^{(n)}_{i\lambda,\varepsilon}=C^{(n-1)}_{i\lambda-\hf,\varepsilon}\,.
\end{equation}\\
\end{proof}

Recalling our present interest in symbols such as
\,$\phi_{a,b}^p$\,, with \,$p\in \Z$\, and \ $\sc{a,\,b}=0$\,, we
first show how, using covariance, the analysis of operators with
such a kind
of symbols can be reduced to a seemingly specialized class.\\

\begin{lemma}
Let \,$a,\,b\in \C^{n+1}$\, be such that \ $\sc{a,\,b}=0\,,\
\sc{a,\,\bar{b}}\neq 0$\,. There exists a matrix \,$g\in G_n$\, such
that the vectors \,$g'\,a$\, and \,${g}^{-1}\,b$\, are both linear
combinations, with complex coefficients, of the vectors \,$e_1$\,
and \,$e_{n+1}$\,
from the canonical basis of \,$\R^{n+1}$\,.\\
\end{lemma}

\begin{proof}
There is no change in the statement if one substitutes for \,$b$\,
any multiple \,$\nu b$\, with \,$\nu\in\C^{\times}$\,, so that we
may assume, without loss of generality, that \
$\sc{a,\,\bar{b}}=2i$\,. Let us decompose the complex vectors
involved as \ $a=p+i\,q,\ b=r+i\,s$\,, so that
\begin{equation}
\sc{p,\,r}=\sc{q,\,s}=0\,,\quad \sc{q,\,r}=1\,,\ \sc{p,\,s}=-1\,.
\end{equation}
Since the matrix \ $\l(\begin{smallmatrix} \sc{q,\,r} & - \sc{p,\,r}\\
 \sc{q,\,s} & - \sc{p,\,s}\end{smallmatrix}\r)$\ is the identity matrix,  it is
possible to find a positive--definite symmetric
\,$(n+1)\times(n+1)$\, matrix \,$h$\, such that \ $h\,r= q\,,\ h\,s=
- p$\,. Let \,$h^{\hf}$\, be the positive--definite square--root of
\,$h$\,. As
\begin{align}
\Vert\,h^{-\hf}\,q\,\Vert^2&=\sc{h^{-1}\,q,\,q}=\,\sc{r,\,q}=1\,,\nonumber\\
\Vert\,h^{-\hf}\,p\,\Vert^2&=\sc{h^{-1}\,p,\,p}=-\,\sc{s,\,p}=1\,,\nonumber\\
\,\sc{h^{-\hf}\,q\,,\,h^{-\hf}\,p}&=\sc{h^{-1}\,q\,,\,p}=\sc{r,\,p}=0\,,
\end{align}
one can find \ $\omega \in {\mathrm{O}}(n+1)$\, such that \
$\omega\,e_1=\pm\,h^{-\hf}\,q\,, \ \omega_2=h^{-\hf}\,p$\,: assuming
\,$n+1\geq 3$\,  (if \,$n=1$\,, the lemma is trivial), one may take
for \,$\omega$\, a rotation matrix. Then, setting \
$g=h^{-\hf}\,\omega$\,, one has \ $p={g'}^{-1}\,e_1\,,\
q={g'}^{-1}\,e_{n+1}$\, and \ $r=h^{-1}\,q=g\,e_{n+1}$\,, finally \
$s=-h^{-1}\,p=-g\,e_1$\,, so that
\begin{equation}
g'\,a=e_1+i\,e_{n+1}\,,\qquad\qquad g^{-1}\,b=-i\,e_1+e_{n+1}\,.
\end{equation}\\
\end{proof}

As made possible by the lemma that precedes, we now specialize to
the case when the symbols \,$\phi_{a,b}^p\,,\ p\in \Z$\,, to be
considered together with their integral superpositions, all
correspond to the case when \,$a$\, and \,$b$\, are linear
combinations of \,$e_1$\, and \,$e_{n+1}$\,: in this way, the
situation is, up to some point, reduced to that obtained when
\,$n=1$\,. Not quite, though, in view of the ever--present
occurrence of the operator \
$D_{\mu}=\sc{s,\,\frac{\partial}{\partial
s}}+\frac{n+1}{2}+i\lambda$\,. However, in this case, setting \
$s=(s_1,\,s_*)$\, with \ $s_*=(s_2,\dots,s_n)$\,, the only operators
we shall have to deal with commute with the partial Euler operator \
$\sc{s_*,\,\frac{\partial}{\partial s_*}}+\frac{n-1}{2}$\, (the
extra constant makes \,$i$\, times this operator a self--adjoint
operator on \,$L^2(\R^{n-1})$): it is thus possible to decompose
functions \ $u=u(s)$\, as integrals of functions
\,$u_{*,\,i\nu,\,\delta}$\, homogeneous of degree and parity
\,$(\frac{1-n}{2}-i\,\nu,\,\delta)$\, with respect to the variables
\,$s_*$\, only. On such a function, the operator \,$D_{\mu}$\,
reduces to \ $s_1\,\frac{d}{ds_1}+1+i\,(\lambda-\nu)$\,: this is
just the analogue of the operator \,$D_{\mu}$\, in a
one--dimensional pseudodifferential calculus
\,${\mathrm{Op}}^{(1)}_{i\,(\lambda-\nu),\varepsilon}$\,. The recipe
for reducing our present analysis to the one--dimensional case thus
essentially calls for replacing \,$\lambda$\, by
\,$\lambda-\nu$\,.\\

To be more specific, let us recall some facts relative to the
discrete terms of the decomposition of \,$L^2(G_1/H_1)$\,. With the
help of the (singular) coordinates \,$(s_1,\,\sigma_1)$\, on
\,${\mathcal X}_1^{\bullet}$\, introduced in (\ref{2ssigma}), we
associate to each complex number \,$z\in \Pi$\,, the upper
half--plane, the function \,$\phi_z$\, such that
\begin{equation}\label{4phi}
\phi_z(s_1,\,\sigma_1)=\frac{(s_1-\bar{z})(1+\bar{z}\,\sigma_1)}{1+s_1\sigma_1}\,:
\end{equation}
an alternative expression, in terms of the homogeneous coordinates
\,$(x,\,\xi)$\, (with \ $\sc{x,\,\xi}=1$\,) of the point of
\,${\mathcal X}_1^{\bullet}$\, considered, is
\begin{equation}
\phi_z(s_1,\,\sigma_1)=\sc{a,\,x}\sc{b,\,\xi}\qquad
{\mathrm{with}}\quad a=\l(\begin{smallmatrix}1 \\
-\bar{z}\end{smallmatrix}\r)\,,\ \,b= \l(\begin{smallmatrix}\bar{z}
\\ 1\end{smallmatrix}\r)\,,
\end{equation}
an expression which may be compared to (\ref{3Fou}).\\

Given \,$k=0,1,\dots$\,, denote as \,$E_{k+1}$\, the closed subspace
of \ $L^2({\mathcal
X}_1^{\bullet})=L^2(\R^2;\,\frac{ds\,d\sigma}{(1+s\sigma)^2})$\,
generated by the functions \ $\phi_z^{-k-1}$\, with \,$z\in \Pi$\,:
this is an irreducible space of the quasiregular representation of
\,$G_1=SL(2,\R)$\, in \,$L^2({\mathcal X}_1^{\bullet})$\,. It makes
up half the eigenspace of \,$\Delta_1$\, for the eigenvalue \
$-k(k+1)$\,: the other half is obtained with the help of the similar
functions related to the lower half--plane. On the other hand, let
us recall that the representation \,$\pi_{2k+2}$\, taken from the
holomorphic discrete series of \,$G_1$\, can be realized in the
space \,${\mathcal D}_{2k+2}$\, consisting of all holomorphic
functions \,$f$\, on \,$\Pi$\, such that
\begin{equation}
\Vert\,f\,\Vert_{2k+2}^2=\int_{\Pi} |f(z)|^2\,(\Im
z)^{2k+2}\,d\mu(z)\,<\infty\,,
\end{equation}
with \ $d\mu(z)=(\Im z)^{-2}\,d\,\Re z\,\wedge\,d\,\Im z$\,. One has
\begin{equation}\label{4pik}
\l(\pi_{2k+2}\l({\generic}\r)f\r)\,(z)=(-cz+a)^{-2k-2}\,f(\frac{dz-b}{-cz+a})\,.
\end{equation}
One may then recall \cite[prop.2.2]{U3} the following. Set \
$\alpha_{k+1}=2^{-2k}\,\l(\begin{smallmatrix}2k \\
k\end{smallmatrix}\r)\,\pi^2$\,, and define the operator
\,$T_{k+1}$\, by
\begin{equation}\label{4Tk+1}
(T_{k+1}\,h)(z)=\alpha_{k+1}^{-1}\,\int_{{\mathcal X}_1^{\bullet}}
h(s_1,\,\sigma_1)\,\overline{\phi}_z^{-k-1}(s_1,\,\sigma_1)\
\frac{ds_1\,d\sigma_1}{(1+s_1\sigma_1)^2}
\end{equation}
for every \ $h\in L^2({\mathcal X}_1^{\bullet})$\, and \,$z\in
\Pi$\,. Then, the operator \
$(\frac{(2k+1)\,\alpha_{k+1}}{4\pi})^{\hf}\,T_{k+1}$\, is an
isometry from \,$E_{k+1}$\, onto \,${\mathcal D}_{2k+2}$\,. It acts
as an intertwining operator between the quasiregular action of
\,$G_1$\, in \,$E_{k+1}$\, and the representation \,$\pi_{2k+2}$\,
of \,$G_1$\, in \,${\mathcal D}_{2k+2}$\,. Its inverse is given by
the formula
\begin{equation}
h(s_1,\,\sigma_1)=\frac{2k+1}{4\pi}\,\int_{\Pi} (T_{k+1}\,h)(z)\
\phi_z^{-k-1}(s_1,\,\sigma_1) \ (\Im z)^{2k+2}\,d\mu(z)\,.
\end{equation}
It has been shown in ({\em loc.cit.\/}) that the Hilbert sum of the
spaces \,$E_{k+1}$\, is an algebra for the sharp composition of
symbols, the sharp products expressing themselves in terms of
Rankin--Cohen brackets of the \,$T_{k+1}$\,--\,transforms of the
terms from the
decompositions of the two symbols under consideration.\\

If a symbol \,$f$\, lies in \,$E_{k+1}$\,, so that it is an integral
superposition of symbols \ $(\sc{a,\,x}\sc{b,\,\xi})^{-k-1}$\, with
\ $ a=\l(\begin{smallmatrix}1 \\ -\bar{z}\end{smallmatrix}\r)\,,\
\,b= \l(\begin{smallmatrix}\bar{z} \\ 1\end{smallmatrix}\r)$\,, and
where  \ $ x=\l(\begin{smallmatrix}x_1 \\
x_2\end{smallmatrix}\r)\,,\ \,\xi= \l(\begin{smallmatrix}\xi_1 \\
\xi_2\end{smallmatrix}\r)$\,, we can turn it to a symbol
\,$\tilde{f}$\, in the
\,${\mathrm{Op}}_{i\lambda,\varepsilon}$\,--\,calculus in \,$n$\,
variables, substituting for the two--dimensional vectors above the
\,$(n+1)$--\,dimensional ones
\begin{equation}
a=\l(\begin{smallmatrix} 1\\ 0 \\ . \\ .\\ .\\0 \\
-\bar{z}\end{smallmatrix}\r)\,,\quad b=\l(\begin{smallmatrix}
\bar{z} \\ 0\\ . \\ . \\ .\\0\\1\end{smallmatrix}\r)\,,\quad
x=\l(\begin{smallmatrix} x_1\\. \\ . \\.\\
x_{n+1}\end{smallmatrix}\r)\,,\quad \xi =\l(\begin{smallmatrix}
\xi_1 \\. \\ . \\.\\ \xi_{n+1}\end{smallmatrix}\r)\,\,:
\end{equation}
taking integral superpositions, with respect to \,$z\in \Pi$\,, of
such symbols, we obtain symbols which can be written, in the
\,$(s,\,\sigma)$\,--\,coordinates on \,${\mathcal X}_n$\,, as
\begin{equation}\label{5tildef}
\tilde{f}(s,\,\sigma)=\l(\frac{1+\sc{s,\,\sigma}}{1+s_1\sigma_1}\r)^{k+1}\,f(s_1,\,\sigma_1)\,.
\end{equation}

As a final topic in this paper, we analyse the operator with symbol
\,$\tilde{f}$\,: in view of the different kind of dependence of the
latter with respect to the two groups of variables involved, one may
start from a decomposition of the function \,$u=u(s_1,\,s_*)$\, to
which the
operator is applied into homogeneous components.\\

\begin{proposition}
Let \,$n\geq 2$\,, and assume that \,$\tilde{f}$\, is given by
{\em(\ref{5tildef})\/}. On functions of \,$s=(s_1,\,s_*)$\,
homogeneous of degree and parity \,$(\frac{1-n}{2}-i\nu,\,\delta)$\,
with respect to the variables \,$s_*$\,, one has
\begin{multline}\label{5decOphom}
({\mathrm{Op}}^{(n)}_{i\lambda,\varepsilon}\,(\tilde{f})\,u)(s)=\frac
{\Gamma(\frac{n+1}{2}+i\lambda)}{\Gamma(\frac{n-1}{2}+i\lambda-k)}\,
\frac{\Gamma(i\,(\lambda-\nu)-k)}{\Gamma(i\,(\lambda-\nu)+1)}\\
\times\
{\mathrm{Op}}^{(1)}_{i(\lambda-\nu),|\varepsilon-\delta|}(f)\,(s_1\mapsto
u(s_1,\,s_*))\,.
\end{multline}\\
\end{proposition}

\begin{proof}
Changing \ $\sigma_*=(\sigma_2,\dots,\sigma_n)$\, to
\,$(1+s_1\,\sigma_1)\,\sigma_*$\, in the integral (\ref{2Op}), and
using the fact that the function \
$(\theta^{(n)}_{i\lambda,\varepsilon}\,u)(\sigma_1,\,\sigma_*)$\, is
homogeneous of degree and parity \,$(\frac{1-n}{2}+i\nu,\,\delta)$\,
with respect to \,$\sigma_*$\,, we obtain
\begin{multline}
({\mathrm{Op}}^{(n)}_{i\lambda,\varepsilon}(\tilde{f})\,u)(s)=(-1)^{\varepsilon}\,
C^{(n)}_{-i\lambda,\varepsilon}\,\int f(s_1,\,\sigma_1)\
|1+s_1\,\sigma_1|_{|\varepsilon-\delta|}^{-1+i\,(\nu-\lambda)}\\
[1+\sc{s_*,\,\sigma_*}]^{k+1}\,|1+\sc{s_*,\,\sigma_*}|_{\varepsilon}^{-\frac{n+1}{2}-i\lambda}
\ (\theta^{(n)}_{i\lambda,\varepsilon}\,u)(\sigma_1,\,\sigma_*)\
d\sigma_1\,d\sigma_*\,:
\end{multline}
expressing \
$(\theta^{(n)}_{i\lambda,\varepsilon}\,u)(\sigma_1,\,\sigma_*)$\,
with the help of Lemma 5.6, one may interpret this as
\begin{equation}
\frac{(-1)^{\varepsilon}\,C^{(n)}_{-i\lambda,\varepsilon}}
{(-1)^{|\varepsilon-\delta|}\,C^{(1)}_{-i(\lambda-\nu),|\varepsilon-\delta|}\,\times\,
(-1)^{\varepsilon}\,C^{(n-1)}_{-i\lambda-\hf,\varepsilon}} \ \times\
(-1)^{|\varepsilon-\delta|}\,C^{(1)}_{i(\nu-\lambda),|\varepsilon-\delta|}
\end{equation}
times
\begin{equation}
\l(\,({\mathrm{Op}}^{(1)}_{i(\lambda-\nu),|\varepsilon-\delta|}(f)\,\otimes\,
({\mathrm{Op}}^{(n-1)}_{i\lambda+\hf,\varepsilon}\,([1+\sc{s_*,\,\sigma_*}]^{k+1}))\,u\,\r)
\,(s_1,\,s_*)\,:
\end{equation}
now, the constant above reduces to $1$ in view of (\ref{5idconst}).
On the other hand, the operator with symbol
\,$[1+\sc{s_*,\,\sigma_*}]^{k+1}$\, has been made explicit in
(\ref{5DmuPoch}): note that \,$\mu=\frac{n+1}{2}+i\lambda$\, does
not change if the pair \,$(n,\,i\lambda)$\, is replaced by
\,$(n-1,\,i\lambda+\hf)$\, and that, in our case,
\,$\sc{s_*,\,\frac{\partial}{\partial s_*}}+\mu$\, reduces to
\,$1+i\,(\lambda-\nu)$\,, which
leads to the result indicated.\\
\end{proof}

{\bf Remark.} The operator under consideration is not bounded in
\,$L^2(\R^n)$\, in view of the pole at \,$\nu=\lambda$\, of the
second Gamma factor on top of the first line of the right--hand side
of (\ref{5decOphom}): but it becomes bounded when composed with the
spectral projection, relative to the self--adjoint operator \
$i\,(\sc{s_*,\,\frac{\partial}{\partial s_*}}+\frac{n-1}{2})$\,,
corresponding to the
complementary, in the real line, of any neighborhood of the point \,$\lambda$\,.\\

The following proposition extends Proposition 5.3 to negative
integral exponents: let us warn the reader that, though a formal
proof, shorter than the one developed below, can be obtained as a
consequence of the equations (\ref{5leftinfOp}), it is only as an
application of Proposition 5.8 that a meaning is given to the
operator with symbol \,$\phi_{a,b}^{-p}$\,, and that it would
be just as much work to extend to this case the validity of the quoted equations.\\

\begin{proposition}
Assume \,$n\geq 2$\,. Let \,$a,\,b\in \C^{n+1}$\, be such that \
$\sc{a,\,b}=0\,, \ \sc{a,\,\bar{b}}\neq 0$\,. Recalling that
\,$\phi_{a,b}$\, has been defined in
{\em(\ref{5phialphabetainh})\/}, one has for \,$p=0,\,1,\dots$\, the
equation
\begin{equation}
({\mathrm{Op}}_{i\lambda,\varepsilon}(\phi_{a,b}))^{-p}=\frac{\mu^p}{(\mu-1)\dots
(\mu-p)}\, {\mathrm{Op}}_{i\lambda,\varepsilon}(\phi_{a,b}^{-p})\,.
\end{equation}\\
\end{proposition}

\begin{proof}
According to Lemma 5.7, it is no loss of generality to assume that \
$a=e_1+i\,e_{n+1},\ b=-i\,e_1+e_{n+1}$\,, in which case, with the
notation in (\ref{5tildef}), one has \
$\phi_{a,b}^{-p}=\tilde{f}_p$\, if one sets
\begin{equation}
f_p(s_1,\,\sigma_1)=\l(\frac{(s_1+i)(1-i\,\sigma_1)}{1+s_1\sigma_1}\r)^{-p}\,.
\end{equation}
Our aim is to prove the equation
\begin{equation}
{\mathrm{Op}}_{i\lambda,\varepsilon}(\phi_{a,b}^{-p-1})=\frac{\mu-p-1}{\mu}
\ {\mathrm{Op}}_{i\lambda,\varepsilon}(\phi_{a,b}^{-1})
\,{\mathrm{Op}}_{i\lambda,\varepsilon}(\phi_{a,b}^{-p})\,,
\end{equation}
using the equations (from Proposition 5.8)
\begin{multline}
({\mathrm{Op}}^{(n)}_{i\lambda,\varepsilon}\,(\phi_{a,b}^{-p})\,u)(s)=\frac
{\Gamma(\frac{n+1}{2}+i\lambda)}{\Gamma(\frac{n+1}{2}+i\lambda-p)}\,
\frac{\Gamma(i\,(\lambda-\nu)-p+1)}{\Gamma(i\,(\lambda-\nu)+1)}\\
\times\
{\mathrm{Op}}^{(1)}_{i(\lambda-\nu),|\varepsilon-\delta|}(f_p)\,(s_1\mapsto
u(s_1,\,s_*))\,,
\end{multline}
valid when applied to functions \,$u=u(s_1,\,s_*)$\, which are
homogeneous of degree and parity \,$(\frac{1-n}{2}-i\nu,\,\delta)$\,
with respect to the variables \,$s_*$\,: the formula reduces to a
formula in the one--dimensional
\,${\mathrm{Op}}^{(1)}_{i(\lambda-\nu),|\varepsilon-\delta|}$\,--
\,calculus, to wit
\begin{equation}\label{5compn1}
f_p\,\#\,f_1=\frac{i\,(\lambda-\nu)}{i\,(\lambda-\nu)-p}\,f_{p+1}\,.
\end{equation}
Of course, when \,$n=1$\,, symbols such as \,$f_p$\, with
\,$p=1,\,2,\dots$\, are square--integrable, so that the composition
is easier to analyze. A detailed proof of (\ref{5compn1}) is to be
found in \cite[Prop.\,4.1]{U3}, but here is some help towards
sorting--out the notation: in ({\em loc.cit.\/}),
\,$(s,\,\sigma)$\, was denoted as \,$(s,\,-t^{-1})$\, so that
\,$f_p$\, would have been denoted as \,$(-1)^p\,g_i^p$\, there;
finally, only the case when \ $ |\varepsilon-\delta|=0$\, was
explicitly considered in this reference, but no change whatsoever
occurs when dealing only with symbols such as \,$f_p$\,, taken from
the discrete spaces of the
decomposition of \,$L^2(G_1/H_1)$\,.\\
\end{proof}

{\bf Remark.} More generally, with the help of the results of ({\em
loc.cit.\/}), together with Proposition 5.8, one can make a
composition such as\,$\tilde{f}_1\,\#\,\tilde{f}_2$\,, with
\,$f_1\in E_{k_1+1}$\, and \,$f_2\in E_{k_2+1}$\,, fully explicit.
We may come back to the more general composition problem at some
later time. Let us just mention, without (the lengthy) proof, the
following result, an analogue of the last proposition, concerned
this time with symbols that
occur in the continuous part of the decomposition of \,$L^2(G_n/H_n)$\,.\\

\begin{proposition}
Let \,$a,\,b\in \R^{n+1}$\, satisfy \ $\sc{a,\,b}=0$\,. Set \
$R_{a,b}=i\mu\,{\mathrm{Op}}_{i\lambda,\varepsilon}(\phi_{a,b})$\,:
this is an (unbounded) self--adjoint operator in \,$L^2(\R^n)$\,
with a purely continuous spectrum, to wit the real line. Denote as
\,$(\Pi_{a,b})_{\pm}$\, the projection operators corresponding to
the positive and negative parts of the spectrum of \,$R_{a,b}$\,,
and set, with \,$\rho\in \C,\ \Re \rho=-\frac{n}{2}\,,\ \delta=0$\,
or $1$\,,
\begin{equation}
(R_{a,b})_{\pm}=\pm\,R_{a,b}\,(\Pi_{a,b})_{\pm}\,,\qquad\quad
|R_{a,b}|_{\delta}^{\rho}=(R_{a,b})_+^{\rho}+(-1)^{\delta}\,(R_{a,b})_-^{\rho}\,.
\end{equation}
Then, one has
\begin{equation}
{\mathrm{Op}}_{i\lambda,\varepsilon}(|\phi_{a,b}|^{\rho}_{\delta})=
(-1)^{\varepsilon}\,i^{\delta}\,2^{-\rho}\,\frac{\Gamma(\frac{\mu+\varepsilon}{2})}
{\Gamma(\frac{1-\mu+\varepsilon}{2})}\,\frac{\Gamma(\frac{1-\mu-\rho+|\varepsilon-\delta|}{2})}
{\Gamma(\frac{\mu+\rho+|\varepsilon-\delta|}{2})} \
|R_{a,b}|_{\delta}^{\rho}\,.
\end{equation}
\end{proposition}
\vspace{1cm}


\end{document}